\documentclass[11pt]{article}
\usepackage{latexsym,amssymb}
\usepackage{authblk}
\usepackage{enumitem}
\usepackage{amsmath}
\usepackage{mathtools}
\usepackage{xcolor}
\usepackage{tikz}

\renewcommand{\baselinestretch}{1.20}
\newcommand{\C}[1]{{\protect\cal #1}}
\newcommand{\B}[1]{{\bf #1}}
\newcommand{\I}[1]{{\mathbb #1}}
\newcommand{\OO}[1]{\overline{#1}}
\newcommand{\e}{\varepsilon}
\renewcommand{\mid}{:}

\definecolor{uuuuuu}{rgb}{0.27,0.27,0.27}
\definecolor{sqsqsq}{rgb}{0.1255,0.1255,0.1255}

\newcommand{\req}[1]{\textrm{(\ref{eq:#1})}}

\newtheorem{theorem}{Theorem}
\newtheorem{lemma}[theorem]{Lemma}

\newtheorem{corollary}[theorem]{Corollary}

\newtheorem{observation}[theorem]{Observation}
\newtheorem{definition}[theorem]{Definition}
\newtheorem{fact}[theorem]{Fact}

\newtheorem{cla}{Claim}[theorem]
\newcommand{\bcl}[2]{\begin{cla}\label{cl:#1}#2\end{cla}}

\newcommand{\bpf}[1][Proof.]{\smallskip\noindent{\it #1} }
\newcommand{\qed}{\nolinebreak\mbox{\hspace{5 true pt}%
  \rule[-0.85 true pt]{3.9 true pt}{8.1 true pt}}}
\newcommand{\cqed}{\nolinebreak\mbox{\hspace{5 true pt}%
  \rule[-0.85 true pt]{2.0 true pt}{8.1 true pt}}}
\newcommand{\epf}{\qed \medskip}
\newcommand{\bcpf}{\bpf[Proof of Claim.]}
\newcommand{\ecpf}{\cqed \medskip}

\newcommand{\IS}{\I S} 
\newcommand{\ex}{\mathrm{ex}}

\newcommand{\branch}[2]{\mathrm{br}_{#1}(#2)}
\newcommand{\blow}[2]{#1(\!(#2)\!)}
\newcommand{\multiset}[1]{\{\hspace{-0.25em}\{\hspace{0.1em}#1\hspace{0.1em}\}\hspace{-0.25em}\}}

\newcommand{\rep}[2]{{#1}^{(#2)}}
\newcommand{\hide}[1]{}
\def\multisets#1#2{\ensuremath{\left(\kern-.3em\left(\genfrac{}{}{0pt}{}{#1}{#2}\right)\kern-.3em\right)}}

  \newcommand{\proj}[1]{{\rm proj\,}#1}

\newcommand{\me}{\mathrm{e}}
\newcommand{\STS}{\mathcal{ST\!S}}

\newcommand{\BV}[1]{\B V_{#1}}
\newcommand{\TruncatedLevel}[1]{^{\mathrm{level}\le #1}}
\newcommand{\SizeTruncated}[1]{^{\mathrm{size}\ge #1}}

\begin{document}

%
\title{\bf\Large Finite Hypergraph Families with Rich Extremal Tur\'an Constructions via Mixing Patterns}

\date{\today}
\author{Xizhi Liu\thanks{Research was supported by ERC Advanced Grant 101020255. Email: xizhi.liu@warwick.ac.uk} }
\author{Oleg Pikhurko\thanks{Research was supported by ERC Advanced Grant 101020255 and
            Leverhulme Research Project Grant RPG-2018-424. Email: o.pikhurko@warwick.ac.uk} }
\affil{Mathematics Institute and DIMAP,
            University of Warwick,
            Coventry, CV4 7AL, UK}
\maketitle
\begin{abstract}
 We prove that, for any finite set of minimal $r$-graph patterns, there is a finite family $\mathcal F$ of forbidden $r$-graphs such that the extremal Tur\'an constructions for $\mathcal F$ are  precisely the maximum $r$-graphs obtainable from mixing the given patterns in any way via blowups and recursion. This extends the result by the second author~\cite{PI14}, where the above statement was established for a single pattern.

We present two applications of this result. First, we construct a finite family $\mathcal F$ of $3$-graphs such that there are exponentially many maximum $\mathcal F$-free $3$-graphs of each large order $n$ and, moreover, the corresponding Tur\'an problem is not finitely stable. Second,  we show that there exists a finite family $\mathcal{F}$ of $3$-graphs
whose feasible region function attains its maximum on a Cantor-type set of positive Hausdorff dimension.

\medskip

\end{abstract}

\section{Introduction}\label{Intro}

\subsection{Tur\'{a}n problem}\label{SUBSEC:Turan-number-stability}

For an integer $r\ge 2$, an \textbf{$r$-uniform hypergraph} (henceforth an \textbf{$r$-graph}) $H$ is a collection of $r$-subsets of some finite set $V$.
Given a family $\mathcal{F}$ of $r$-graphs, we say that $H$ is \textbf{$\mathcal{F}$-free}
if it does not contain any member of $\mathcal{F}$ as a subgraph.
The \textbf{Tur\'{a}n number} $\mathrm{ex}(n,\mathcal{F})$ of $\mathcal{F}$ is the maximum
number of edges in an $\mathcal{F}$-free $r$-graph on $n$ vertices.
The \textbf{Tur\'{a}n density} $\pi(\mathcal{F} )$ of $\mathcal{F}$ is defined as
$\pi(\mathcal{F})\coloneqq \lim_{n\to \infty}\mathrm{ex}(n,\mathcal{F})/\binom{n}{r}$; the
existence of the limit was established in~\cite{KatonaNemetzSimonovits64}.
The study of $\mathrm{ex}(n,\mathcal{F})$ is one of the central topic in extremal graph and hypergraph theory.

Much is known about $\mathrm{ex}(n,\mathcal{F})$ for graphs, that is, when $r=2$.
For example, Tur\'{a}n~\cite{Turan41} determined $\mathrm{ex}(n,K_{\ell}^2)$ for all $n> \ell> 2$ (where, more generally, $K^r_\ell$ denotes the complete $r$-graph on $\ell$ vertices).
Also, the Erd\H{o}s--Stone--Simonovits theorem~\cite{ErdosSimonovits66,ErdosStone46} determines the Tur\'{a}n density for every family $\mathcal{F}$ of graphs; namely, it holds that $\pi(\mathcal{F})=\min\{1-1/\chi(F)\mid F\in\mathcal{F}\}$, where $\chi(F)$ denotes the chromatic number of the graph~$F$.

For $r\ge 3$, determining $\pi(\mathcal{F})$ for a given family $\mathcal{F}$ of $r$-graphs
seems to be extremely difficult in general.
For example, the problem of determining $\pi(K_{\ell}^{r})$ raised already in the 1941 paper by Tur\'{a}n~\cite{Turan41}
 is still open for all $\ell>r\ge 3$;
 thus the $\$ 500$ prize of Erd\H{o}s (see e.g.\ \cite[Section III.1]{Erdos81c}) for determining $\pi(K_{\ell}^{r})$ for at least one pair $(\ell,r)$ with $\ell > r \ge 3$ remain unclaimed.


The ``inverse'' problem of understanding the sets
\begin{align*}
\Pi_{\mathrm{fin}}^{(r)} & \coloneqq  \left\{\pi(\mathcal{F}) \colon \text{$\mathcal{F}$ is a finite family of $r$-graphs} \right\}, \quad\text{and} \\
\Pi_{\infty}^{(r)} & \coloneqq  \left\{\pi(\mathcal{F}) \colon \text{$\mathcal{F}$ is a (possibly infinite) family of $r$-graphs} \right\}
\end{align*}
of possible $r$-graph Tur\'an densities is also very difficult for $r\ge 3$. (For $r=2$, we have by the Erd\H{o}s--Stone--Simonovits theorem~\cite{ErdosSimonovits66,ErdosStone46} that
$
\Pi_{\infty}^{(2)} = \Pi_{\mathrm{fin}}^{(2)}=  \{1\}\cup \left\{1-{1}/{k} \colon \mbox{integer $k\ge 1$}\right\}$.)

One of the earliest results on this direction is the theorem of Erd\H{o}s~\cite{E64} from the 1960s that $\Pi_{\infty}^{(r)} \cap (0, r!/r^r) = \emptyset$ for every integer $r\ge 3$. However, our understanding of the locations and the lengths of other maximal intervals avoiding $r$-graph Tur\'an densities and the right accumulation points of $\Pi_{\infty}^{(r)}$ (the so-called \textbf{jump problem}) is very limited; for some results in this direction see e.g.~\cite{BT11,FranklPengRodlTalbot07,FR84,Pikhurko15,YanPeng21}.

It is known that the set $\Pi_{\infty}^{(r)}$ is the topological closure of $\Pi_{\mathrm{fin}}^{(r)}$
(see \cite[Proposition~1]{PI14}) and thus the former set is determined by the latter.
In order to show that the set $\Pi_{\mathrm{fin}}^{(r)}\subseteq [0,1]$ has rich structure for each $r\ge 3$, 
the second author proved in~\cite[Theorem~3]{PI14} that, for every minimal $r$-graph pattern $P$, there is a finite family $\mathcal{F}$ of $r$-graphs such that the maximum $\mathcal{F}$-free graphs are precisely the maximum $r$-graphs that can be obtained by taking blowups of $P$
and using recursion. 
(See Section~\ref{SUBSEC:Pattern} for all formal definitions.)
In particular, the maximum asymptotic edge density
obtainable this way from the pattern $P$ is an element of~$\Pi_{\mathrm{fin}}^{(r)}$.

Another factor that makes the hypergraph Tur\'an problem difficult is  that some families may have many rather different (almost) extremal configurations.
A series of recent papers~\cite{HLLMZ22,LM22,LMR3,LMR1} (discussed in more detail in Sections~\ref{se:ExpMany}  and~\ref{SUBSEC:feasible-region}) concentrated on exhibiting examples for which the richness of extremal configurations can be proved.

Our paper contributes further to this line of research. The new results proved here are, informally speaking,  as follows.
Our main result, on which all new constructions are based, is Theorem~\ref{THM:mixed-pattern}. It extends~\cite[Theorem~3]{PI14}
to the case when there is a finite set $\{P_i\mid i\in I\}$ of minimal patterns (instead of just one) and we can mix them in any way when using recursion. 
By applying Theorem~\ref{THM:mixed-pattern}, we present two
examples of a $3$-graph family $\C F$ with a rich set of (almost) extremal configurations. The first one (given by Theorem~\ref{th:ExpMany})
has the property that the set of  maximum $\C F$-free $3$-graphs on $[n]$ has exponentially many in $n$ non-isomorphic hypergraphs and, moreover, the Tur\'an problem for $\C F$ is not finitely stable, that is, roughly speaking, there is no bounded number of constructions such that every almost maximum $\C F$-free $3$-graph is close to one of them. The second  finite family $\C F$ of $3$-graphs
(given by Corollary~\ref{CORO:feasible-function-Cantor})
satisfies the property that the limit set of possible densities of the shadows of asymptotically maximum $\C F$-free $3$-graphs is a Cantor-like set of positive Hausdorff dimension.

Let us now present the formal statements of the new results (together with some further definitions and background).

\subsection{Patterns}\label{SUBSEC:Pattern}

In order to state our main result (Theorem~\ref{THM:mixed-pattern}), we need to give a number of definitions.

Let an \textbf{$r$-multiset} mean an unordered collection
of $r$ elements with repetitions allowed. Let $E$ be  a
collection of $r$-multisets on $[m]\coloneqq \{1,\dots,m\}$, let $V_1,\dots,V_m$
be disjoint sets and set $V\coloneqq V_1\cup\dots\cup V_m$. The \textbf{profile}
of an $r$-set $X\subseteq V$ (with respect to $V_1,\dots,V_m$) is
the $r$-multiset on $[m]$ that contains $i\in [m]$
with multiplicity $|X\cap V_i|$. For an $r$-multiset $Y\subseteq [m]$, let
$\blow{Y}{V_1,\dots,V_m}$ consist of all $r$-subsets of $V$ whose profile is
$Y$. We call this $r$-graph the \textbf{blowup} of $Y$ (with respect to $V_1,\dots,V_m$) and
the $r$-graph
 \begin{align*}
 \blow{E}{V_1,\dots,V_m}\coloneqq \bigcup_{Y\in E} \blow{Y}{V_1,\dots,V_m}
 \end{align*}
 is called the \textbf{blowup} of $E$ (with respect to $V_1,\dots,V_m$).

An \textbf{($r$-graph) pattern} is a triple $P=(m,E,R)$ where $m$ is a positive integer, $E$ is a
collection of $r$-multisets on $[m]$,
and $R$ is a subset of $[m]$ (we allow $R$ to be the empty set).
As a convention, we require that $R$ does not contain $i$ if $E$ contains the multiset consisting of $r$ copies of~$i$.  
Suppose that $I$ is a non-empty index set and
 \begin{equation}\label{eq:PI}
 P_I\coloneqq \{P_i  \colon i\in I\},\quad\mbox{where $P_i= (m_i, E_i, R_i)$ for $i\in I$},
 \end{equation}
 is a collection of $r$-graph patterns indexed by $I$.

\begin{definition}[$P_I$-Mixed Constructions]
For $P_I$ as in~\eqref{eq:PI}, a \textbf{$P_I$-mixing construction} on a set $V$ is any $r$-graph with vertex set $V$ which has no edges or can be recursively constructed as follows. Pick some $i\in I$ and take any partition $V=V_1\cup\cdots\cup V_{m_i}$ such that $V_j\not=V$ for each $j\in R_i$. Let $G$ be obtained from the blowup $\blow{E_i}{V_1,\ldots,V_{m_i}}$ by adding, for each $j\in R_i$, an arbitrary $P_I$-mixing construction on $V_j$.
\end{definition}

Informally speaking,
we can start with a blowup $\blow{E_i}{V_1,\dots,V_{m_i}}$ for some $i\in I$, then put a blowup of some $E_{i'}$ inside each \textbf{recursive part} (that is, $V_j$ with $j\in R_i$), then put blowups into the new recursive parts, and so on.
Note that there is no restriction on the choice of the index at any step. For example, on the second level of the recursion, different recursive parts $V_j$, $j\in R_i$, may choose different indices. See Section~\ref{example} for two examples illustrating this construction.

The family of all $P_I$-mixing constructions will be denoted by~$\Sigma P_{I}$.
We say $G$ is a \textbf{$P_I$-mixing subconstruction} if it is a subgraph of some $P_{I}$-mixing construction on $V(G)$.

Let $\Lambda_{P_I}(n)$ be the maximum number of edges that an $r$-graph in~$\Sigma P_{I}$ with $n$ vertices can have:
\begin{align}\label{eq:LambdaP}
\Lambda_{P_I}(n)\coloneqq \max\big\{\,|G|\mid \mbox{$G$ is a $P_I$-mixing construction on $[n]$}\,\big\}.
\end{align}
 Using a simple averaging argument (see Lemma~\ref{lm:monotone}), one can show
that the ratio $\Lambda_{P_I}(n)/\binom{n}{r}$ is non-increasing and therefore
tends to a limit which
we denote by $\lambda_{P_I}$ and call it the \textbf{Lagrangian}
of $P_I$:
\begin{align}
\lambda_{P_I}\coloneqq \lim_{n\to\infty} \frac{\Lambda_{P_I}(n)}{\binom{n}{r}}.
\end{align}

If $P_I=\{P\}$ consists of a single pattern $P$ then we always have to use this pattern $P$ and the definition of a $P_{I}$-mixing construction coincides with the definition of a \textbf{$P$-construction} from~\cite{PI14}.
For brevity, we abbreviate $\Lambda_{P}\coloneqq \Lambda_{\{P\}}$, $\lambda_{P}\coloneqq \lambda_{\{P\}}$, etc.

For example, if $r=2$ and $P=(2,\{\,\{1,2\}\,\},\emptyset)$ then $P$-constructions (that is, $\{P\}$-mixing constructions) are exactly complete bipartite graphs, $P$-subconstructions are exactly graphs with chromatic number at most $2$, $\Lambda_{P}(n)=\lfloor n^2/4\rfloor$ for every integer $n\ge 0$, and $\lambda_P=1/2$.

For a pattern $P=(m,E,R)$ and $j\in [m]$, let $P-j$ be the pattern obtained
from $P$ by \textbf{removing index $j$}, that is, we remove $j$ from $[m]$ and delete
all multisets containing $j$ from $E$ (and relabel the
remaining indices to form the set $[m-1]$). In other words,
$(P-j)$-constructions are precisely those $P$-constructions where we
always let the $j$-th part be empty.
We call $P$ \textbf{minimal} if $\lambda_{P-j}$ is strictly smaller than $\lambda_{P}$ for every
$j\in [m]$. For example, the 2-graph pattern
$P\coloneqq (3,\{\,\{1,2\},\{1,3\}\,\},\emptyset)$ is not minimal as
$\lambda_{P}=\lambda_{P-3}=1/2$.

Let $\C F_\infty
$ be the family consisting of those $r$-graphs that are not $P_I$-mixing subconstructions, that is,
\begin{align}\label{eq:CFInfty}
    \C F_\infty\coloneqq \{\mbox{$r$-graph $F$}\mid
\mbox{every  $P_{I}$-mixing construction $G$ is $F$-free}\},
\end{align}
and for every $M\in \mathbb{N}$ let $\C F_M$ be the collection of members in $\C F_\infty$ with at most $M$ vertices, that is, for $v(F)\coloneqq |V(F)|$, we have
\begin{align}\label{eq:CFn}
    \C F_M\coloneqq \{F\in\C F_\infty\mid v(F)\le M\}.
\end{align}


Our main result is as follows.

\begin{theorem}\label{THM:mixed-pattern}
Let $r\ge 3$ and let $P_I=\{P_i\mid i\in I\}$ be an arbitrary collection of minimal $r$-graph patterns, where the index set $I$ is finite.
Then there exists $M\in \mathbb{N}$ such that the following statements hold.
\begin{enumerate}[label=(\alph*)]
\item\label{it:a}
For every positive integer $n$ we have
$\mathrm{ex}(n,\mathcal{F}_M) = \max\left\{|G|\colon v(G) = n \text{ and }G\in \Sigma P_{I}\right\}$.
Moreover, every  maximum $n$-vertex $\C F_M$-free $r$-graph is a member in $\Sigma P_{I}$.

\item\label{it:b}
For every $\varepsilon>0$ there exist $\delta>0$ and $N_0$ such that for every $\C F_M$-free $r$-graph $G$ on $n\ge N_0$ vertices with $|G|\ge (1-\delta)\mathrm{ex}(n,\mathcal{F}_M)$,
there exists an $r$-graph $H\in \Sigma P_{I}$ on $V(G)$ such that $|G\triangle H|\le \varepsilon n^r$.
\end{enumerate}
\end{theorem}

Note that Part~\ref{it:a} of Theorem~\ref{THM:mixed-pattern} gives that the family of maximum $\C F_M$-free $r$-graphs is exactly the family of maximum $P_I$-mixing constructions.

In the case of a single pattern (when $|I| = 1$), Theorem~\ref{THM:mixed-pattern} gives \cite[Theorem~3]{PI14}.
As we will see in Lemma~\ref{LEMMA:lambda-sigma}, it holds that $\lambda_{P_I}=\max\{\lambda_{P_i}\mid i\in I\}$ and thus we do not increase the set of obtainable Tur\'an densities by mixing patterns. The main purpose of Theorem~\ref{THM:mixed-pattern} is to show that some finite Tur\'an problems have rich sets of (almost) extremal graphs. In this paper, we present two applications of Theorem~\ref{THM:mixed-pattern} of this kind as follows.

\subsection{Finite families with exponentially many extremal Tur{\'a}n $r$-graphs}\label{se:ExpMany}

Let us call a family $\C F$ of $r$-graphs \textbf{$t$-stable} if for every $n\in\I N$ there are $r$-graphs $G_{1}(n),\dots,G_t(n)$ on $[n]$ such that for every $\e>0$ there are $\delta>0$ and $n_0$ so that if $G$ is an $\C F$-free $r$-graph with $n\ge n_0$ vertices and least $(1-\delta)\ex(n,\C F)$ edges then $G$ is within edit distance $\e n^r$ to some $G_i(n)$. The \textbf{stability number} $\xi(\C F)$ of $\C F$ is the smallest $t\in\I N$ such that $\C F$ is $t$-stable; we set $\xi(\C F)\coloneqq \infty$ if no such $t$ exists. We call $\C F$ \textbf{stable} if $\xi(\C F)=1$ (that is, if $\C F$ is 1-stable). According to our definition, every family $\C F$ of $r$-graphs with $\pi(\C F)=0$ is stable: take $G_1(n)$ to be the edgeless $r$-graph on~$[n]$.

The first stability theorem, which says that $K_{\ell}^2$ is stable for all integers $\ell\ge 3$,
was established independently by Erd\H{o}s~\cite{Erdos67a} and Simonovits $\cite{SI68}$.
In fact, the classical Erd\H{o}s--Stone--Simonovits theorem~\cite{ErdosStone46,ErdosSimonovits66} and Erd\H{o}s--Simonovits stability theorem~\cite{Erdos67a,SI68}
imply that every  family of graphs is stable.

For hypergraphs, there are some conjectures on the Tur\'an density of various concrete families which, if true, imply that these families are not stable.
One of the most famous examples in this regard is the tetrahedron $K_{4}^3$ whose conjectured Tur\'an density is $5/9$. If this conjecture is true then the constructions by Brown~\cite{BR83} (see also \cite{Fonderflaass88,Froh08,KO82}) show that $\xi(K_4^3)=\infty$. A similar statement applies to some other complete $3$-graphs $K_\ell^3$;
we refer the reader to~\cite{KE11,Sido95} for details.
Another natural example of conjectured infinite stability number
is the Erd\H{o}s--S\'{o}s Conjecture on triple systems with bipartite links; we refer the reader to~\cite{FF84} for details.

Despite these old conjectures, no finite family with more than one asymptotic Tur\'an extremal construction was known until recently.
In~\cite{LM22}, Mubayi and the first author constructed the first finite non-stable
family $\C F$ of $3$-graphs; in fact, their family satisfies $\xi(\C F)=2$.
Further, in~\cite{LMR1}, Mubayi, Reiher, and the first author found, for every integer $t\ge 3$, a finite family $\C F_t$ of $3$-graphs with $\xi(\C F_t)= t$.
Their construction was extended to families of $r$-graphs for every integer $r \ge 4$ in~\cite{LMR3}.
In~\cite{HLLMZ22}, Hou, Li, Mubayi, Zhang, and the first author constructed
a finite family $\C F$ of $3$-graphs such that $\xi(\C F)=\infty$.

Note that it is possible that $\xi(\C F)=1$ but there are many maximum $\C F$-free $r$-graphs of order~$n$. For example, if $k\ge 5$ is odd and we forbid the star $K_{1,k}^2$ (where, more generally, $K_{k_1,\ldots,k_\ell}^2$ denotes
the complete $\ell$-partite graph with part sizes $k_1,\ldots,k_\ell$), then the extremal graphs on $n\ge r$ vertices are precisely $(k-1)$-regular graphs and, as it is easy to see, there are exponentially many in $n$ such non-isomorphic graphs. For $r$-graphs with $r\ge 3$, a similar conclusion for an infinite sequence of $n$ can be achieved by forbidding e.g.\ two $r$-edges intersecting in $r-1$ vertices: if a sufficiently large $n$ satisfies the obvious divisibility conditions then by the result of Keevash~\cite{Keevash18} there are $\exp(\Omega(n^{r-1}\log n))$ extremal $r$-graphs on $[n]$, namely designs
where each $(r-1)$-set is covered exactly once. While the above Tur\'an problems are \textbf{degenerate} (i.e.\ have the Tur\'an density 0), a non-degenerate example for graphs can be obtained by invoking a result of Simonovits~\cite{SI68}, a special case of which is that every maximum graph of order $n\to\infty$ without $K_{1,t,t}^2$ can be obtained from a complete bipartite graph $K_{a,n-a}$ with $a=(1/2+o(1))n$ by adding a maximum $K_{1,t}^2$-free graph into each part. Very recently, Balogh, Clemen, and Luo~\cite{BaloghClemenLuo} found a single $3$-graph $F$ with $\pi(F)>0$ and with $\exp(\Omega(n^{2}\log n))$ non-isomorphic extremal constructions on $n$ vertices for an infinite sequence of~$n$. Note that all families in this paragraph are 1-stable.

In the other direction, the relation $\xi(\C F)=\infty$ does not generally imply that there are many maximum $\C F$-free $r$-graphs, since one of asymptotically optimal constructions may produce strictly better bounds (in lower order terms) on $\ex(n,\C F)$ than any other. So it is of interest to produce $\C F$ with many extremal graphs and with $\xi(\C F)=\infty$.
The above mentioned paper~\cite{HLLMZ22} made a substantial progress towards this problem, by giving a finite familly $\C F$ of $3$-graphs such that, in addition to $\xi(\C F)=\infty$, there are $\Omega(n)$ non-isomorphic maximum $\C F$-free for infinitely many $n$ (e.g.\ for all large $n$ divisible by $3$).

As an application of Theorem~\ref{THM:mixed-pattern}, we provide a finite family of $3$-graphs with infinite stability number and exponentially many extremal constructions for \textbf{all} large integers~$n$.

\begin{theorem}\label{th:ExpMany}
There is a finite family $\C F$ of $3$-graphs such that, for some $C>0$ and for every $n\in\I N$, the number of non-isomorphic maximum $\C F$-free $3$-graphs on $n$ vertices is at least $\mathrm{exp}(Cn)$. Additionally, $\xi(\C F)=\infty$.
\end{theorem}

\subsection{Feasible region}\label{SUBSEC:feasible-region}
Theorem~\ref{THM:mixed-pattern} has also applications to the so-called feasible region problem of hypergraphs.
To state our results we need more definitions.

Given an $r$-graph $G$, its \textbf{$s$-shadow} is defined as
\begin{align}
\partial_s G
\coloneqq 
\left\{A\in \binom{V(G)}{r-s}\colon \exists\, B\in G \text{ such that }
	A\subseteq B\right\}. \notag
\end{align}
We use $\partial G$ to represent $\partial_1 G$ and call it the \textbf{shadow} of $G$.
The \textbf{edge density} of $G$ is defined as $\rho(G)\coloneqq  |G|/\binom{v(G)}{r}$, and
the \textbf{shadow density} of $G$ is defined as $\rho(\partial G)\coloneqq  |\partial G|/\binom{v(G)}{r-1}$.

For a family $\mathcal{F}$, the \textbf{feasible region} $\Omega(\mathcal{F})$ of $\mathcal{F}$
is the set of points $(x,y)\in [0,1]^2$ such that there exists a sequence of $\mathcal{F}$-free $r$-graphs
$\left( G_{n}\right)_{n=1}^{\infty}$ with
\begin{align*}\lim_{n \to \infty}v(G_{n}) = \infty,\quad
\lim_{n \to \infty}\rho(\partial G_{n}) = x, \quad\text{and}\quad \lim_{n \to \infty}\rho(G_{n}) = y.\end{align*}
The feasible region unifies and generalizes asymptotic versions of some classical problems
such as the  Kruskal--Katona theorem~\cite{KA68,KR63} and the Tur\'{a}n problem.
It was introduced in \cite{LM1} to understand the extremal properties of $\mathcal{F}$-free hypergraphs
beyond just the determination of $\pi(\mathcal{F})$.

The following general results about $\Omega(\mathcal{F})$ were proved
in~\cite{LM1}. The projection of $\Omega(\mathcal{F})$ to the first
coordinate,
\[
	\proj{\Omega(\mathcal{F})}
	\coloneqq 
	\left\{ x \colon  \text{there is $y \in [0,1]$ such that $(x,y) \in \Omega(\mathcal{F})$} \right\},
\]
is an interval $[0, c(\mathcal{F})]$ for some $c(\mathcal{F})\in [0, 1]$.
Moreover, there is a left-continuous almost everywhere differentiable function
$g(\mathcal{F})\colon \proj{\Omega(\mathcal{F})} \to [0,1]$ such that
\[
	\Omega(\mathcal{F})
	=
	\bigl\{(x, y)\in [0, c(\mathcal{F})]\times [0, 1]\colon 0\le y\le g(\mathcal{F})(x)\bigr\}\,.
\]
The function $g(\mathcal{F})$ is called the \textbf{feasible region function} of $\mathcal{F}$.
It was showed in~\cite{LM1} that $g(\mathcal{F})$ is not necessarily continuous,
and it was showed in~\cite{Liu20a} that $g(\mathcal{F})$ can
have infinitely many local maxima even for some simple and natural families~$\mathcal{F}$.

For every $r$-graph family $\mathcal{F}$, define the set $M(\mathcal{F}) \subseteq \proj{\Omega(\mathcal{F})}$ as the collection of $x$ such that $g(\mathcal{F})$ attains its global maximum at $x$, that is,
\begin{align}
M(\mathcal{F}) \coloneqq  \{x\in \proj{\Omega(\mathcal{F})}\colon g(\mathcal{F})(x) = \pi(\mathcal{F})\}. \notag
\end{align}

For most families $\mathcal{F}$ that were studied before, the set $M(\mathcal{F})$ has size one,
i.e.\ the function $g(\mathcal{F})$ attains its maximum at only one point.
In general, the set $M(\mathcal{F})$ is not necessarily a single point,
and, in fact, trying to understand how complicated $M(\mathcal{F})$ can be is one of the motivations for constructions in~\cite{HLLMZ22,LM22,LMR1}.
Indeed, the construction in~\cite{LM22} shows that there exists a finite family $\mathcal{F}$ of $3$-graphs for which
$M(\mathcal{F})$ has size exactly two, the
constructions in~\cite{LMR1} show that for every positive integer $t$
there exists a finite family $\mathcal{F}$ of $3$-graphs for which $M(\mathcal{F})$ has size exactly $t$,
and the constructions in~\cite{HLLMZ22} show that there exists a finite family $\mathcal{F}$ of $3$-graphs for which
$M(\mathcal{F})$ is a non-trivial interval.

We show that $M(\mathcal{F})$ can be even more complicated than the examples above.
More specifically, we show that there exists a finite family $\mathcal{F}$ of $3$-graphs for which $M(\mathcal{F})$ is a \textbf{Cantor-type set}, i.e.\ is topologically homeomorphic to the standard Cantor set $\left\{\sum_{i=1}^\infty t_i 3^{-i}\colon \forall\,i\ t_i\in\{0,2\}\right\}$. For this, we need some further preliminaries.

Suppose that $P_i = (m_i, E_i, R_i)$ for $i\in I$ is a collection of patterns with the same Lagrangian, say $\lambda$ for some $\lambda\in [0,1]$. Recall that $\C F_\infty$ is defined by~\eqref{eq:CFInfty}; thus $\C F_\infty$-free graphs are exactly subgraphs of $P_I$-mixing constructions.
It follows that
$M(\C F_\infty)$ can be equivalently described as the set of all points $x\in [0,1]$
such that there exists a sequence $(H_n)_{n=1}^{\infty}$ of $r$-graphs such that $H_n$ is a $P_I$-mixing construction for all $n\ge 1$ and
\begin{align}
\lim_{n\to \infty}v(H_n) = \infty, \quad \lim_{n\to \infty}\rho(H_n) = \lambda, \quad\text{and}\quad \lim_{n\to \infty}\rho(\partial H_n) = x. \notag
\end{align}

Using a standard diagonalization argument in analysis we obtain the following observation.

\begin{observation}\label{OBS:shadow-density-mixing-construction}
The set $M(\C F_\infty)$ is a closed subset of $[0,1]$.
\end{observation}

Using Theorem~\ref{THM:mixed-pattern} and some further arguments we obtain the following result for $3$-graphs.

\begin{theorem}\label{THM:feasibe-region}
Suppose that $I$ is a finite set, $P_i = (m_i, E_i, R_i)$ is a minimal $3$-graph pattern for each $i\in I$, and there exists $\lambda \in (0,1)$ such that $\lambda_{P_i} = \lambda$ for all $i\in I$.
Then there exists a finite family $\mathcal{F}\subseteq \C F_\infty$ such that $M(\C F) = M(\C F_\infty)$.
\end{theorem}

Later, we will give specific patterns $P_1$ and $P_2$ with the same Lagrangian such that, for $P_I=\{P_1,P_2\}$, the set $M(\C F_\infty)$ is a
Cantor-type set with Hausdorff dimension ${\log 2}/{\log \left(4 \sqrt{7}+11\right)} \approx 0.225641$.
Thus, by Theorem~\ref{THM:feasibe-region}, we obtain the following corollary.

\begin{corollary}\label{CORO:feasible-function-Cantor}
There exists a finite family $\mathcal{F}$ of $3$-graphs such that the set $M(\C F)$ is a Cantor-type set with Hausdorff dimension
${\log 2}/{\log \left(4 \sqrt{7}+11\right)}$.
\end{corollary}

\medskip
\noindent \textbf{Organisation of the paper.} The rest of the paper is organized as follows.
In Section~\ref{notation}, we define some further notation and present two examples.
Theorems~\ref{THM:mixed-pattern}, \ref{th:ExpMany} and~\ref{THM:feasibe-region} are proved in 
Sections~\ref{SEC:Proof-Turan-number},  \ref{se:rich} and~\ref{SEC:feasible-region}, respectively.
Corollary~\ref{CORO:feasible-function-Cantor} is proved in Section~\ref{SEC:examples}.
Some concluding remarks are contained in Section~\ref{se:conclusion}.

\section{Notation}\label{notation}

Let us introduce some further notation complementing and expanding
that from the Introduction. Some other (infrequently used) definitions
are given shortly before they are needed for the first time in this
paper.

Let $\mathbb{N}\coloneqq \{1,2,\ldots\}$ be the set of all positive integers.
 Also, recall that $[m]$ denotes the set $\{1,\dots,m\}$.

Recall that an $r$-\textbf{multiset} $D$ is an unordered
collection of $r$ elements $x_1,\dots,x_r$ with repetitions
allowed. Let us denote this as $D=\multiset{x_1,\dots,x_r}$. The \textbf{multiplicity $D(x)$} of $x$ in $D$ is the number of
times that $x$ appears. If the underlying set is understood to be
$[m]$, then we can represent $D$ as the ordered
$m$-tuple $(D(1),\dots,D(m)) \in \{0,\ldots, r\}^{m}$ of multiplicities. Thus, for example, the profile
of $X\subseteq V_1\cup\dots\cup V_m$ is the multiset on $[m]$ whose
multiplicities are $(|X\cap V_1|,\dots,|X\cap V_m|)$. Also, let
$\rep{x}{r}$ denote the sequence consisting of $r$ copies of $x$; thus
the multiset consisting of $r$ copies of $x$ is denoted by $\multiset{\rep{x}{r}}$.
 If we need to
emphasise that a multiset is in fact a set (that is, no element
has multiplicity more than 1), we call it a \textbf{simple set}.

For $D\subseteq [m]$ and sets $U_1,\dots,U_m$, denote $U_D\coloneqq \bigcup_{i\in D} U_i$.
Given a set $X$ let $\binom{X}{r}$ and $\multisets{X}{r}$ denote the collections of all $r$-subsets of $X$ and all $r$-multisets of $X$, respectively.

\subsection{Hypergraphs}\label{SUBSEC:Def-Hypergraphs}
Let $G$ be an $r$-graph. The \textbf{complement} of $G$
is $\OO G\coloneqq \{E\subseteq V(G)\mid |E|=r, \ E\not\in G\}$.
For a vertex $v\in V(G)$, its \textbf{link} is the $(r-1)$-hypergraph
\begin{align*}
    L_{G}(v)\coloneqq \{E\subseteq V(G)\mid v\not\in E,\ E\cup\{v\}\in G\}.
\end{align*}
For $U\subseteq V(G)$,
its \textbf{induced subgraph} is $G[U]\coloneqq \{E\in G\mid E\subseteq U\}$.
The vertex sets of $\OO G$, $L_{G}(v)$, and $G[U]$ are by default $V(G)$,
$V(G)\setminus\{v\}$, and $U$ respectively.
The \textbf{degree} of $v\in V(G)$ is $d_G(v)\coloneqq |L_{G}(v)|$.
Let $\Delta(G)$ and $\delta(G)$ denote respectively the maximum and minimum
degrees of the $r$-graph $G$.

An \textbf{embedding} of an
$r$-graph $F$ into $G$ is an injection $f:V(F)\to V(G)$ such that
$f(E)\in G$ for every $E\in F$. Thus $F$ is a subgraph of $G$ if and only if $F$ admits an embedding into $G$.
An embedding is \textbf{induced} if non-edges are mapped to non-edges.

The \textbf{edit distance} between two $r$-graphs $G$ and $H$ with the same number of vertices is the smallest number  of \textbf{edits} (i.e.\ removal and addition of edges) that have to be applied to $G$ to make it isomorphic to $H$; in other words, this is the minimum of $|G\bigtriangleup \sigma(H)|$ over all bijections $\sigma:V(H)\to V(G)$. We say that $G$ and $H$ are \textbf{$s$-close} if they are at  edit distance at most~$s$.

\subsection{Further definitions and results for a single pattern}\label{SUBSEC:Def-patterns}

In this section we focus on a single pattern $P = (m, E, R)$. Let the \textbf{Lagrange polynomial} of $E$ be
\begin{align}\label{eq:LP}
     \lambda_E(x_1,\dots,x_m)\coloneqq r!\,\sum_{D\in E}\; \prod_{i=1}^m\;
        \frac{x_i^{D(i)}}{D(i)!}.
\end{align}
The special case of~\req{LP} when $E$ is an $r$-graph (i.e.\
$E$ consists of simple sets)
is a generalisation
of the well-known \textbf{hypergraph Lagrangian} (see e.g.~\cite{BT11,FR84}) that
has been successfully applied to Tur\'an-type problems, with
the basic idea going back to
Motzkin and Straus~\cite{MS65}.

For $i\in [m]$, let the \textbf{link} $L_{E}(i)$ consist of all $(r-1)$-multisets $A$
such that if we increase the multiplicity of $i$ in $A$ by one, then the
obtained $r$-multiset belongs to $E$.

The following simple fact follows easily from the definitions.
\begin{fact}\label{FACT:Lagrangian}
    The following statements hold for every $r$-graph pattern $P = (m, E, R)$. 
    \begin{enumerate}[label=(\alph*)]
        \item\label{FACT:Lagrangian-1} For every partition $[n]=V_1\cup\dots\cup V_m$, we have that
        \begin{align}\label{eq:lambdaE}
            \lambda_E\left(\frac{|V_1|}n,\dots,\frac{|V_m|}n\right)
            = \rho(\blow{E}{V_1,\dots,V_m}) +o(1),\qquad  \mbox{as $n\to\infty$}.
         \end{align}
         \item\label{FACT:Lagrangian-2} For every $i\in [m]$, we have  $\frac{\partial \lambda_{E}}{\partial_i}(\B x) = r\cdot \lambda_{L_{E}(i)}(\B x)$.
    \end{enumerate}
\end{fact}
See also Lemma~\ref{lm:Omega1} that
relates $\lambda_E$ and $\lambda_{P}$.

We call a pattern $P$ \textbf{proper} if it is minimal and $0<\lambda_{P}<1$.
Trivially, every minimal pattern $P=(m,E,R)$ satisfies that
\begin{align}\label{eq:MinLinks}
    L_{E}(i)\not=\emptyset,\qquad \mbox{for every $i\in [m]$.}
\end{align}

\begin{lemma}[{\cite[Lemma~16]{PI14}}]
\label{lm:DNs}
Let $P=(m,E,R)$ be a minimal pattern.
If distinct $j,k \in[m]$ satisfy $L_{E}(j)\subseteq L_{E}(k)$, then
$j\in R$, $k\not\in R$, and $L_E(j) \neq L_E(k)$.
In particular, no two vertices in $[m]$ have the same links in $E$.
\end{lemma}
We will also need the following result from \cite{PI14}, which characterizes the patterns whose Lagrangian is $1$.

\begin{lemma}[{\cite[Lemma~12]{PI14}}]
\label{lm:density=1}
An $r$-graph pattern $P = (m,E,R)$ satisfies $\lambda_{P}=1$ if and only if at least one of the following holds:
\begin{enumerate}[label=(\alph*)]
\item there is $i\in [m]$ such that $\multiset{\rep{i}{r}}\in E$, or
\item there are $i\in R$ and $j\in[m]\setminus\{i\}$ such that
$\multiset{\rep{i}{r-1},j}\in E$.
\end{enumerate}
\end{lemma}

The \textbf{standard $(m-1)$-dimensional simplex} is
 \begin{align}\label{eq:Sm}
 \IS_m\coloneqq \{\B x\in \I R^m\mid x_1+\dots+x_m=1,\ \forall\, i\in[m]\ x_i\ge 0\}.
 \end{align}
Also, let
\begin{align*}
 \IS_m^*\coloneqq \{\B x\in\I R^m\mid x_1+\dots+x_m=1,\ \forall\,i\in[m]\ 0\le x_i<1\}
 \end{align*}
 be obtained from the simplex $\IS_m$
 by excluding all its vertices (i.e.\ the standard basis
vectors).

For a pattern $P=(m,E,R)$, we say that a vector $\B x\in\I R^{m}$ is \textbf{$P$-optimal} if $\B x\in \IS_{m}^*$  and
\begin{align}\label{eq:opt}
   \lambda_{E}(\B x) + \lambda_{P} \sum_{j\in R} x_j^r = \lambda_{P}.
\end{align}
Note that when we define $P$-optimal vectors we restrict ourselves
to $\IS_{m}^*$, i.e.\ we do not allow any standard basis vector to be included.

We will need the following result from~\cite{PI14}, which extends some classical results (see e.g.~\cite[Theorem~2.1]{FR84}) about the Lagrangian of hypergraphs.

\begin{lemma}[{\cite[Lemma~14]{PI14}}]\label{lm:Omega1}
Let $P=(m,E,R)$ be a proper $r$-graph pattern and let
\begin{align*}f(\B x)\coloneqq  \lambda_{E}(\B x) + \lambda_{P} \sum_{j\in R} x_j^r
\end{align*} be
the left-hand side of~\req{opt}. Let $\C X$ be the set of all $P$-optimal vectors. Then the following statements hold.
\begin{enumerate}[label=(\alph*)]
\item\label{it:part1} The set $\C X$ is non-empty.
\item\label{it:f} We have $f(\B x)\le \lambda_{P}$ for all $\B x\in \IS_{m}$.
\item\label{it:boundary} The set $\C X$ does not intersect the boundary of $\IS_{m}$, i.e.\ no vector in $\C X$ has a zero coordinate.
\item\label{it:f'} For every $\B x\in\C X$ and $j\in [m]$ we have $\frac{\partial f}{\partial_j}(\B x) = r\cdot \lambda_{P}$.
\item\label{it:closed}  The set $\C X$, as a subset of  $\IS_{m}$, is closed. (In particular, no sequence of $P$-optimal vectors can converge to a standard basis vector.)
\item\label{it:almost} For every $\e>0$ there is $\alpha>0$ such that for every
$\B y\in\IS_{m}$ with $\max(y_1,\dots,y_{m})\le 1-\e$ and $f(\B y)\ge \lambda_{P}-\alpha$ there
is $\B x\in\C X$ with $\|\B x-\B y\|_\infty\le \e$.
\item\label{it:separated} There is a constant $\beta>0$ such that for every $\B x\in \C X$ and every $j\in[m]$ we have $x_j\ge \beta$.
\end{enumerate}
\end{lemma}

Observe that, in the above lemma, Parts~\ref{it:part1} and \ref{it:f} imply that
$\C X$ is precisely the set of elements in $\IS_{m}^*$ that maximise $\lambda_{E}(\B x) + \lambda_{P} \sum_{j\in R} x_j^r$. Also, note that \ref{it:boundary} is a consequence of~\ref{it:separated}.


We will also need the following easy inequality.

\begin{lemma}\label{LEMMA:Taylor-inequ}
Suppose that $E$ is a collection of $r$-multisets on $[m]$ and
$\B u, \B x \in \mathbb{S}_m$.
Then for every $j\in [m]$ we have
\begin{align*}
\left|\frac{\partial\lambda_{E}}{\partial_{j}}(\B u) -\frac{\partial\lambda_{E}}{\partial_{j}}(\B x)\right|
\le r^2 m \cdot \|\B u - \B x\|_\infty.
\end{align*}
\end{lemma}
\bpf
Fix $j_{\ast} \in [m]$, by Fact~\ref{FACT:Lagrangian}\ref{FACT:Lagrangian-2}, we have $\frac{\partial\lambda_{E}}{\partial_{j_{\ast}}}(\B y) = r\cdot \lambda_{L_{E}(j_{\ast})}(\B y)$ for every $\B y \in \mathbb{S}_{m}$. 
Let $\hat{E} \coloneqq L_{E}(j_{\ast})$. 
Similarly, for every $i \in [m]$ and $\B y \in \mathbb{S}_{m}$, we have $\frac{\partial\lambda_{\hat{E}}}{\partial_{i}}(\B y) = (r-1)\cdot \lambda_{L_{\hat{E}}(i)}(\B y)$, and hence, 
\begin{align*}
    \left|\frac{\partial^2\lambda_{E}}{\partial_{i}\partial_{j_{\ast}}}(\B y)\right|
    = r\cdot \frac{\partial\lambda_{\hat{E}}}{\partial_{i}}(\B y) 
    & = r(r-1)\cdot \lambda_{L_{\hat{E}}(i)}(\B y) \\
    & \le r^2 \sum_{S\in \multisets{[m]}{r-2}}\prod_{i\in S} y_i
    = r^2 \left(\sum_{i=1}^{m} y_i\right)^{r-2}
    =  r^2. 
\end{align*}
Given $\B u, \B x \in \mathbb{S}_m$, it follows from the Mean Value Theorem that there exists $\B y = \alpha \B u+ (1-\alpha) \B x \in \mathbb{S}_{m}$ for some $\alpha \in [0,1]$ such that 
\begin{align*}
    \left| \frac{\partial\lambda_{E}}{\partial_{j_{\ast}}}(\B u) -\frac{\partial\lambda_{E}}{\partial_{j_{\ast}}}(\B x)\right|
     = \left| \sum_{i \in [m]}  \frac{\partial^2\lambda_{E}}{\partial_{i}\partial_{j_{\ast}}}(\B y) \cdot (u_i - x_i)\right|  
    & \le \max_{i\in [m]}\left\{ \left|\frac{\partial^2\lambda_{E}}{\partial_{i}\partial_{j_{\ast}}}(\B y)\right| \right\} \cdot \sum_{i=1}^{m}\left|u_i-x_i\right| \\
    & \le  r^2\sum_{i=1}^{m}\left|u_i-x_i\right|
    \le r^2m \cdot \|\B u -\B x\|_\infty. 
\end{align*}
This proves Lemma~\ref{LEMMA:Taylor-inequ}.
\epf

\subsection{Mixing patterns}\label{se:MixingPatterns}

Let~\eqref{eq:PI} apply, that is, $P_I=\{P_i\mid i \in I\}$, where $P_{i} = (m_i, E_i, R_i)$ is an $r$-graph pattern for each~$i\in I$.

First, to each $P_I$-mixing construction $G$ we will associate its \textbf{base value}~$b(G)$. For this, we need to designate some special element not in $I$, say, $\emptyset\not\in I$. If $G$ has no edges, then we set $b(G)\coloneqq \emptyset$; otherwise let $b(G)$ be the element of $I$ such that the first step when making the $P_I$-mixing construction $G$ was to take a blowup of~$E_j$. It may in principle be possible that some different choices of $j$ lead to another representation of the same $r$-graph $G$ as a $P_I$-mixing construction. We fix one choice for $b(G)$, using it consistently. We say that $G$ has \textbf{base} $P_{b(G)}$ and call $b$ the \textbf{base function}.

Next, to every $G\in\Sigma P_I$ we will associate a pair $(\BV{G},\B T_G)$ which encodes
how the $P_{I}$-mixing construction $G$ is built. In brief, the \textbf{tree} $\B T_G$ of $G$ records the information how the patterns are mixed while the \emph{partition structure} $\BV{G}$ specifies which vertex partitions were used when making each blowup.

Let us formally define $\BV{G}$ and $\B T_G$, also providing some related terminology. We start with $\BV{G}$ having only one set $V_{\emptyset}\coloneqq V(G)$ and with $\B T_{G}$ consisting of the single root $((), \emptyset)$. If $G$ has no edges, then this finishes the definition of $(\BV{G},\B T_G)$.
So suppose that $|G|>0$. By the definition of $j\coloneqq b(G)$,
there exists a partition $V(G) = V_1 \cup \dots \cup V_{m_{j}}$ with $V_i\not=V$ for $i\in R_{j}$
such that $G\setminus \left(\bigcup_{k\in R_{j}}G[V_k]\right) = \blow{E_{j}}{V_1, \ldots, V_{m_{j}}}$.
(Again it may be possible that some different choices of the partition lead to another representation of $G\in\Sigma P_I$, so fix one choice and use it consistently.) This initial partition $V(G)=V_1\cup\dots\cup V_{m_{j}}$ is called the \textbf{level-$1$} partition
and $V_1, \ldots, V_{m_{j}}$ are called the \textbf{level-$1$} (or \textbf{bottom}) parts.  We add
the parts $V_1,\ldots,V_{m_{j}}$ to $\BV{G}$ and add $((1), j), \ldots, ((m_{j}), j)$ as the children of the root $((),\emptyset)$ of $\B T_G$. Note that we add all $m_{j}$ children, even if some of the corresponding parts are empty.   Next, for every $i\in R_{j}$ with $|G[V_i]|>0$, we apply recursion to $G[V_i]$ to define the descendants of $((i),j)$ as follows. Let $P_t$ be the base of $G[V_i]$ (that is, $t\coloneqq b(G[V_i])$). Add $((i,1),t),\dots,(i,m_t),t)$ as the children of $((i), j)$ in~$\B T_G$ and add to $\BV{G}$ the sets  $V_{i,1},\dots, V_{i,m_t}$ (called \textbf{level-$2$} parts) forming  the partition of $V_i$  that was used for the blowup of~$E_t$.
In other words, the set of level-$2$ parts of $G$ is the union over $i\in R_{j}$
of level-1 parts of~$G[V_i]$. We repeat this process. Namely, a node $(\B i,j)$ in $\B T_G$, where $\B i=(i_1,\ldots,i_s)$, has children if and only if  $|G[V_{\B i}]|>0$ (then necessarily $i_s\in R_j$), in which case the children are $((\B i,k),t)$ for $k\in [m_t]$, where $t\coloneqq b(G[V_{\B i}])$; then we also add the \textbf{level-$(s+1)$} parts $V_{\B i,1},\dots,V_{\B i,m_t}$ to $\BV{G}$ which we used for the blowup of~$E_t$. Note that we are slightly sloppy with our bracket notation, with $\B i,k$ meaning $i_1,\dots,i_s,k$ (and with $\emptyset$ sometimes meaning the empty sequence).
Thus level-$(s+1)$ parts of $G$ can be defined as the level-$s$ parts of $G[V_i]$ for all $i\in R_{b(G)}$. This finishes the definition of $(\BV{G},\B T_G)$.  For some examples, see Section~\ref{example}.

Clearly, the pair $(\BV{G},\B T_G)$ determines $G$ (although this representation has some redundancies). We often
work with just $\BV{G}$ (without explicitly referring to $\B T_G$). When $G$ is understood, we may write $\BV{}$ for~$\BV{G}$.

A sequence $(i_1,\dots,i_s)$ is called \textbf{legal} (with respect to $G$) if $V_{i_1,\ldots,i_s}$ is defined during the above process. This includes the empty sequence, which is always legal.
We view $\BV{G}$ and $\B T_G$ as vertical (see e.g.\ Figure~\ref{fig:ordered-tree})
with the index sequence
growing as we go up in level. This motivates calling the level-$1$ parts \textbf{bottom} (even though we regard $V_{\emptyset}$ as level-$0$). By default, the profile of $X\subseteq V(G)$
is taken with respect to the bottom parts, that is,
its multiplicities are $(|X\cap V_1|,\dots,|X\cap V_m|)$. The \textbf{height} of $G$ is the maximum length of a legal sequence (equivalently, the maximum number of edges in a directed path in the tree $\B T_G$).
Call a part $V_{i_1,\ldots,i_s}$ \textbf{recursive} if $s=0$ or the (unique) node $((i_1,\ldots,i_s),t)$ of $\B T_G$ satisfies $i_s\in R_t$. (Note that we do not require that $|G[V_{i_1,\dots,i_s}]|>0$ in this definition.) Otherwise, we call $V_{i_1,\ldots,i_s}$ \textbf{non-recursive}.
The \textbf{branch $\branch{G}{v}$} of a vertex $v\in V(G)$ is the (unique) maximal sequence $\B i$ such that $v\in V_{\B i}$. If $b(G)=\emptyset$, then every branch is the empty sequence; otherwise every branch is non-empty. Let $\B T_{G}\TruncatedLevel{\ell}$ be obtained from $\B T_{G}$ by removing all nodes at level higher than~$\ell$.



Note that the values of the base function $b$ are incorporated into the tree $\B T_G$ by ``shifting'' them one level up: namely,  if a part $V_{\B i}$ is recursive and  $j\coloneqq b(G[V_{\B i}])$ is not $\emptyset$, then $j$ appears as the second coordinate on all children of the unique node of $\B T_G$  with the first coordinate~$\B i$. This is notationally convenient: for example, if $G'$ is obtained from $G\in\Sigma P_I$ by deleting all edges inside some part $V_{\B i}$ then $\B T_{G'}$ is a \textbf{subtree} of $\B T_G$, that is, it can be obtained from $\B T_G$ by iteratively deleting leaves. (In fact, the subtree $\B T_{G'}$ is \textbf{full}, that is, every non-leaf node of $\B T_{G'}$ has the same set of children in $\B T_{G'}$ as in $\B T_{G}$.)

We say that $\B T$ is a \textbf{feasible tree} if there exists a
$P_{I}$-mixing construction $G$ with~$\B T_{G}=\B T$.
We say that a feasible  tree $\B T$ is \textbf{non-extendable} if it has positive height and every leaf $v=((i_1,\dots,i_s),t)$ of $\B T$ satisfies $i_s\not\in R_t$.
Equivalently, $\B T$ is non-extendable if it is not a subtree  of a strictly larger feasible tree.
Otherwise, we say that $\B T$ is \textbf{extendable}.%
\hide{
We say $\B V$ is a
\textbf{feasible partition structure} if there exists a
$ P_{I}$-mixing construction $G$ whose partition structure is $\B V$.
We say a feasible partition structure $\B V$ is \textbf{non-extendable} if for every $P_I$-mixing construction $G$ with partition structure $\B V$ and for every leaf $(V_{\B i}, *)$ the set $V_{\B i}$ is non-recursive in $G$.
Otherwise, we say $\B V$ is \textbf{extendable}.
}


\hide{By the \textbf{height}, \textbf{tree}, etc, of $G\in\Sigma P_I$, we  mean the height, tree, etc, of $\B V_G$.
For example, the height of any edgeless $G$ is $0$ while, if $R_i=\emptyset$ for each $i\in I$, then every $P_I$-mixing construction $G$ with $|G|>0$ has height exactly~$1$. If $G$ is understood, then we may omit the subscript $G$ in $\B V_{G}$.}

Given a family $P_I=\{P_i \colon i\in I\}$ of $r$-graph patterns, recall that $\C F_\infty$ consists of those $r$-graphs $F$ that do not embed into a $P_{I}$-mixing construction.
For an integer $s$, let $\C F_s$ consist of all members of $\C F_\infty$ with at most $s$ vertices.

Note that we do not require here (nor in Theorem~\ref{THM:mixed-pattern}) that all patterns $P_i$ have the same Lagrangian (although this will be the case in all concrete applications of Theorem~\ref{THM:mixed-pattern} that we present in this paper). As we will see in Lemma~\ref{LEMMA:lambda-sigma}, patterns $P_i$ with $\lambda_{P_i}< \max\{\lambda_{P_j}\mid j\in I\}$ will not affect the largest asymptotic density of $P_I$-mixing constructions (however, they may appear in maximum constructions at the recursion level when we consider parts of bounded size).

\subsection{Examples}\label{example}

In this section we give some examples to illustrate some of the above definitions. The set $P_I$ of the first example  will be later used in our proof of Theorem~\ref{THM:feasibe-region}.

Recall that $K_{5}^3$ is the complete $3$-graph on $5$ vertices (let us assume that its vertex set is $[5]$).
Let $B_{5,3}$ be the $3$-graph on $7$ vertices (let us assume that its vertex set is $[7]$) whose edge set is the union of all triples that have at least two vertices in $\{1,2,3\}$ and all triples  (with their ordering ignored) from the following set:
\begin{align}
\Big(\{1\}\times \{4,5\}\times \{6,7\}\Big) \cup \Big(\{2\}\times \{4,6\}\times \{5,7\}\Big) \cup \Big(\{3\}\times \{4,7\}\times \{5,6\}\Big). \notag
\end{align}
In particular, for $i\in \{1,2,3\}$ the induced subgraph of $L_{B_{5,3}}(i)$ on $\{4,5,6,7\}$ is a copy of $K_{2,2}$ (where $K_{m,n}$ denotes the complete bipartite graph with parts of size $m$ and $n$), and their union covers each pair in $\{4,5,6,7\}$ exactly twice. See Figure~\ref{fig:B53-links} for an illustration.

The motivation for defining $B_{5,3}$ comes from the so-called \textbf{crossed blowup} defined in~\cite{HLLMZ22} so, in fact,  $B_{5,3}$ is a 3-crossed blowup of $K_{5}^{3}$.

Let $P_1 \coloneqq  (5,K_{5}^{3}, \{1\})$ and $P_2 \coloneqq  (7,B_{5,3}, \{1\})$.
Suppose that $G$ is a $\{P_1, P_2\}$-mixing construction with exactly three levels:
the base for $G$ is $P_1$,
the base for $G[V_1]$ is $P_2$, and
the base for $G[V_{1,1}]$ is $P_1$ (see Figure~\ref{fig:mixing-construction}).
It is clear to see that every feasible tree (see e.g.\ Figure~\ref{fig:ordered-tree}) is extendable since every pattern $P_i$ has non-empty~$R_i$.

For the second example, we take $P_1 = (3, \{123\}, \{1\})$ and $P_2 = (4, K_4^3, \emptyset)$.
Let $G$ be a $\{P_1, P_2\}$-mixing construction with three levels:
the base for $G$ is $P_1$,
the base for $G[V_1]$ is $P_1$,
and the base for $G[V_{1,1}]$ is $P_2$, see Figure~\ref{fig:mixing-construction-2}.
Note that the tree $\B T_G$ of $G$ (see Figure~\ref{fig:ordered-tree-2}) of $G$ is non-extendable, since every leaf in $\B T_G$ is non-recursive.

\begin{figure}[htbp]
\centering
\tikzset{every picture/.style={line width=0.75pt}} 
\begin{tikzpicture}[x=0.75pt,y=0.75pt,yscale=-1,xscale=1]
\draw[line width=1pt,color=sqsqsq]   (128,49) -- (198,49) -- (128,105) -- (198,105) -- cycle ;
\draw[line width=1pt,color=sqsqsq]  (385.22,44.18) -- (384.78,114.18) -- (329.22,43.82) -- (328.78,113.82) -- cycle ;
\draw[line width=0.5pt,dash pattern=on 1pt off 1.2pt]  (231,25) -- (230,134) ;
\draw[line width=0.5pt,dash pattern=on 1pt off 1.2pt]  (431,27) -- (430,136) ;
\draw[line width=1pt,color=sqsqsq]   (520,47) -- (584,47) -- (584,111) -- (520,111) -- cycle ;
\draw [fill=uuuuuu] (63,49) circle (1.5pt);
\draw (63,49+10) node {$1$};
\draw [fill=uuuuuu] (266,85) circle (1.5pt);
\draw (266,85+10) node {$2$};
\draw [fill=uuuuuu] (457,113.82) circle (1.5pt);
\draw (457,113.82+10) node {$3$};
\draw [fill=uuuuuu] (128,49) circle (1.5pt);
\draw (128-10,49) node {$4$};
\draw [fill=uuuuuu] (128,105) circle (1.5pt);
\draw (128-10,105) node {$5$};
\draw [fill=uuuuuu] (198,49) circle (1.5pt);
\draw (198+10,49) node {$6$};
\draw [fill=uuuuuu] (198,105) circle (1.5pt);
\draw  (198+10,105) node {$7$};
\draw [fill=uuuuuu] (520,47) circle (1.5pt);
\draw (520-10,47) node {$4$};
\draw [fill=uuuuuu] (329.22,43.82) circle (1.5pt);
\draw (329.22-10,43.82) node {$4$};
\draw [fill=uuuuuu] (520,111) circle (1.5pt);
\draw (520-10,111) node {$5$};
\draw [fill=uuuuuu] (328.78,113.82) circle (1.5pt);
\draw (328.78-10,113.82) node {$5$};
\draw [fill=uuuuuu] (584,47) circle (1.5pt);
\draw (584+10,47) node {$6$};
\draw [fill=uuuuuu] (385.22,44.18) circle (1.5pt);
\draw (385.22+10,44.18) node {$6$};
\draw [fill=uuuuuu] (584,111) circle (1.5pt);
\draw (584+10,111) node {$7$};
\draw [fill=uuuuuu] (384.78,114.18) circle (1.5pt);
\draw (384.78+10,114.18) node {$7$};
\end{tikzpicture}
\caption{The induced subgraph of $L_{B_{5,3}}(i)$ on vertex set $\{4,5,6,7\}$ is a copy of $K_{2,2}$ for $i\in \{1,2,3\}$.}
\label{fig:B53-links}
\end{figure}
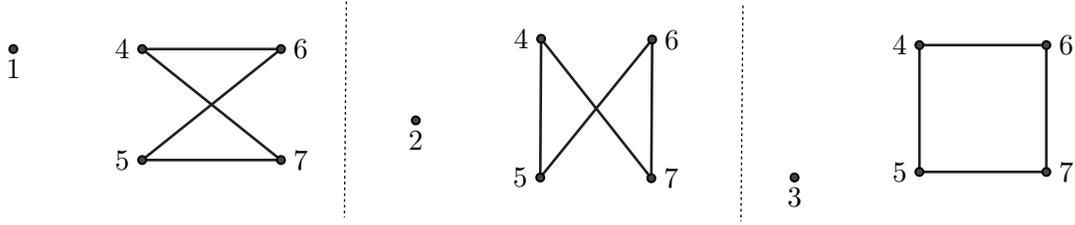

\begin{figure}[htbp]
\centering
\tikzset{every picture/.style={line width=0.75pt}} 

\begin{tikzpicture}[x=0.75pt,y=0.75pt,yscale=-1,xscale=1]

\draw[line width=1pt,color=sqsqsq]   (241,93.5) .. controls (241,50.7) and (275.7,16) .. (318.5,16) .. controls (361.3,16) and (396,50.7) .. (396,93.5) .. controls (396,136.3) and (361.3,171) .. (318.5,171) .. controls (275.7,171) and (241,136.3) .. (241,93.5) -- cycle ;
\draw[line width=1pt,color=sqsqsq]   (281,64.5) .. controls (281,43.79) and (297.79,27) .. (318.5,27) .. controls (339.21,27) and (356,43.79) .. (356,64.5) .. controls (356,85.21) and (339.21,102) .. (318.5,102) .. controls (297.79,102) and (281,85.21) .. (281,64.5) -- cycle ;
\draw[line width=1pt,color=sqsqsq]   (312,156.5) .. controls (312,151.25) and (316.25,147) .. (321.5,147) .. controls (326.75,147) and (331,151.25) .. (331,156.5) .. controls (331,161.75) and (326.75,166) .. (321.5,166) .. controls (316.25,166) and (312,161.75) .. (312,156.5) -- cycle ;
\draw[line width=1pt,color=sqsqsq]   (379,195) .. controls (379,167.94) and (400.94,146) .. (428,146) .. controls (455.06,146) and (477,167.94) .. (477,195) .. controls (477,222.06) and (455.06,244) .. (428,244) .. controls (400.94,244) and (379,222.06) .. (379,195) -- cycle ;
\draw[line width=1pt,color=sqsqsq]   (155,190) .. controls (155,162.94) and (176.94,141) .. (204,141) .. controls (231.06,141) and (253,162.94) .. (253,190) .. controls (253,217.06) and (231.06,239) .. (204,239) .. controls (176.94,239) and (155,217.06) .. (155,190) -- cycle ;
\draw[line width=1pt,color=sqsqsq]   (206,290) .. controls (206,262.94) and (227.94,241) .. (255,241) .. controls (282.06,241) and (304,262.94) .. (304,290) .. controls (304,317.06) and (282.06,339) .. (255,339) .. controls (227.94,339) and (206,317.06) .. (206,290) -- cycle ;
\draw[line width=1pt,color=sqsqsq]   (326,294) .. controls (326,266.94) and (347.94,245) .. (375,245) .. controls (402.06,245) and (424,266.94) .. (424,294) .. controls (424,321.06) and (402.06,343) .. (375,343) .. controls (347.94,343) and (326,321.06) .. (326,294) -- cycle ;
\draw[line width=1pt,color=sqsqsq]   (308,41.5) .. controls (308,36.25) and (312.25,32) .. (317.5,32) .. controls (322.75,32) and (327,36.25) .. (327,41.5) .. controls (327,46.75) and (322.75,51) .. (317.5,51) .. controls (312.25,51) and (308,46.75) .. (308,41.5) -- cycle ;
\draw[line width=1pt,color=sqsqsq]   (362,123.5) .. controls (362,118.25) and (366.25,114) .. (371.5,114) .. controls (376.75,114) and (381,118.25) .. (381,123.5) .. controls (381,128.75) and (376.75,133) .. (371.5,133) .. controls (366.25,133) and (362,128.75) .. (362,123.5) -- cycle ;
\draw[line width=1pt,color=sqsqsq]   (340,145.5) .. controls (340,140.25) and (344.25,136) .. (349.5,136) .. controls (354.75,136) and (359,140.25) .. (359,145.5) .. controls (359,150.75) and (354.75,155) .. (349.5,155) .. controls (344.25,155) and (340,150.75) .. (340,145.5) -- cycle ;
\draw[line width=1pt,color=sqsqsq]   (252,107.5) .. controls (252,98.39) and (259.39,91) .. (268.5,91) .. controls (277.61,91) and (285,98.39) .. (285,107.5) .. controls (285,116.61) and (277.61,124) .. (268.5,124) .. controls (259.39,124) and (252,116.61) .. (252,107.5) -- cycle ;
\draw[line width=1pt,color=sqsqsq]   (274,140.5) .. controls (274,131.39) and (281.39,124) .. (290.5,124) .. controls (299.61,124) and (307,131.39) .. (307,140.5) .. controls (307,149.61) and (299.61,157) .. (290.5,157) .. controls (281.39,157) and (274,149.61) .. (274,140.5) -- cycle ;
\draw[line width=1pt,color=sqsqsq]   (368,93.5) .. controls (368,88.25) and (372.25,84) .. (377.5,84) .. controls (382.75,84) and (387,88.25) .. (387,93.5) .. controls (387,98.75) and (382.75,103) .. (377.5,103) .. controls (372.25,103) and (368,98.75) .. (368,93.5) -- cycle ;
\draw[line width=1pt,color=sqsqsq]   (288,58.5) .. controls (288,53.25) and (292.25,49) .. (297.5,49) .. controls (302.75,49) and (307,53.25) .. (307,58.5) .. controls (307,63.75) and (302.75,68) .. (297.5,68) .. controls (292.25,68) and (288,63.75) .. (288,58.5) -- cycle ;
\draw[line width=1pt,color=sqsqsq]   (297,83.5) .. controls (297,78.25) and (301.25,74) .. (306.5,74) .. controls (311.75,74) and (316,78.25) .. (316,83.5) .. controls (316,88.75) and (311.75,93) .. (306.5,93) .. controls (301.25,93) and (297,88.75) .. (297,83.5) -- cycle ;
\draw[line width=1pt,color=sqsqsq]   (325,83.5) .. controls (325,78.25) and (329.25,74) .. (334.5,74) .. controls (339.75,74) and (344,78.25) .. (344,83.5) .. controls (344,88.75) and (339.75,93) .. (334.5,93) .. controls (329.25,93) and (325,88.75) .. (325,83.5) -- cycle ;
\draw[line width=1pt,color=sqsqsq]   (331,59.5) .. controls (331,54.25) and (335.25,50) .. (340.5,50) .. controls (345.75,50) and (350,54.25) .. (350,59.5) .. controls (350,64.75) and (345.75,69) .. (340.5,69) .. controls (335.25,69) and (331,64.75) .. (331,59.5) -- cycle ;

\draw (307,178) node [anchor=north west][inner sep=0.75pt]   [align=left] {$V_1$};
\draw (191,182) node [anchor=north west][inner sep=0.75pt]   [align=left] {$V_2$};
\draw (242,281) node [anchor=north west][inner sep=0.75pt]   [align=left] {$V_3$};
\draw (371,289) node [anchor=north west][inner sep=0.75pt]   [align=left] {$V_4$};
\draw (422,185) node [anchor=north west][inner sep=0.75pt]   [align=left] {$V_5$};
\draw (310,105) node [anchor=north west][inner sep=0.75pt]   [align=left] {$V_{1,1}$};
\draw (260-5,99) node [anchor=north west][inner sep=0.75pt]   [align=left] {$V_{1,2}$};
\draw (282-5,132) node [anchor=north west][inner sep=0.75pt]   [align=left] {$V_{1,3}$};
\end{tikzpicture}
\caption{The partition structure of a $\{P_1, P_2\}$-mixing construction $G$ with exactly three levels: the base for level-1 is $P_1$, while
the bases for the (unique) recursive parts at levels 2 and 3 are respectively $P_2$ and $P_1$.}
\label{fig:mixing-construction}
\end{figure}
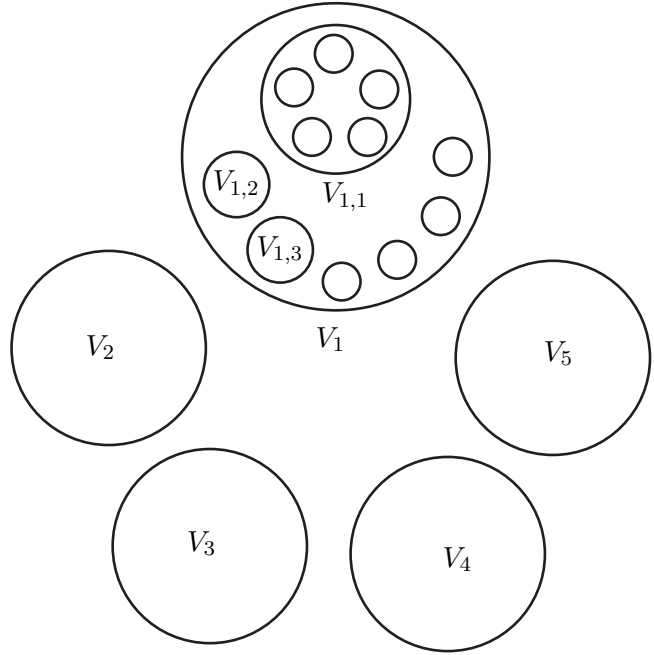

\begin{figure}[htbp]
\centering
\tikzset{every picture/.style={line width=0.75pt}} 
\begin{tikzpicture}[x=0.75pt,y=0.75pt,yscale=-1,xscale=1]

\draw[line width=1pt,color=sqsqsq]    (515,-38) -- (383,-90) ;
\draw[line width=1pt,color=sqsqsq]    (515,-38) -- (574,-90) ;
\draw[line width=1pt,color=sqsqsq]    (515,-38) -- (635,-90) ;
\draw[line width=1pt,color=sqsqsq]    (515,-38) -- (461,-90) ;
\draw[line width=1pt,color=sqsqsq]    (515,-38) -- (515,-90) ;
\draw[line width=1pt,color=sqsqsq]    (383,-90) -- (253,-160) ;
\draw[line width=1pt,color=sqsqsq]    (383,-90) -- (491,-160) ;
\draw[line width=1pt,color=sqsqsq]    (383,-90) -- (460,-160) ;
\draw[line width=1pt,color=sqsqsq]    (383,-90) -- (306,-160) ;
\draw[line width=1pt,color=sqsqsq]    (383,-90) -- (338,-160) ;
\draw[line width=1pt,color=sqsqsq]    (383,-90) -- (376,-160) ;
\draw[line width=1pt,color=sqsqsq]    (383,-90) -- (416,-160) ;
\draw[line width=1pt,color=sqsqsq]    (253,-160) -- (212,-230) ;
\draw[line width=1pt,color=sqsqsq]    (253,-160) -- (254,-230) ;
\draw[line width=1pt,color=sqsqsq]    (253,-160) -- (292,-230) ;
\draw[line width=1pt,color=sqsqsq]    (253,-160) -- (331,-230) ;
\draw[line width=1pt,color=sqsqsq]    (253,-160) -- (163,-230) ;

\draw (80,-38) node {Level-$0\colon$};
\draw (80,-90) node {Level-$1\colon$};
\draw (80,-160) node {Level-$2\colon$};
\draw (80,-230) node {Level-$3\colon$};
\draw [fill=uuuuuu] (515,-38) circle (1.5pt);
\draw (515,-38+12) node {\footnotesize{$((), \emptyset)$}};
\draw [fill=uuuuuu] (383,-90) circle (1.5pt);
\draw (383-7,-90+12) node {\footnotesize{$((1),1)$}};
\draw [fill=uuuuuu] (461,-90) circle (1.5pt);
\draw (461,-90-10) node {\footnotesize{$((2),1)$}};
\draw [fill=uuuuuu] (515,-90) circle (1.5pt);
\draw (515,-90-10) node {\footnotesize{$((3),1)$}};
\draw [fill=uuuuuu] (574,-90) circle (1.5pt);
\draw (574,-90-10) node {\footnotesize{$((4),1)$}};
\draw [fill=uuuuuu] (635,-90) circle (1.5pt);
\draw (635,-90-10) node {\footnotesize{$((5),1)$}};
\draw [fill=uuuuuu] (253,-160) circle (1.5pt);
\draw (253-15,-160+12) node {\footnotesize{$((1,1),2)$}};
\draw [fill=uuuuuu] (306,-160) circle (1.5pt);
\draw (306,-160-12) node {\footnotesize{$((1,2),2)$}};
\draw [fill=uuuuuu] (338,-160) circle (1.5pt);
\draw [fill=uuuuuu] (376,-160) circle (1.5pt);
\draw (376,-160-12) node {$\cdots$};
\draw [fill=uuuuuu] (416,-160) circle (1.5pt);
\draw (416,-160-12) node {$\cdots$};
\draw [fill=uuuuuu] (460,-160) circle (1.5pt);
\draw [fill=uuuuuu] (491,-160) circle (1.5pt);
\draw (491,-160-12) node {\footnotesize{$((1,7),2)$}};
\draw [fill=uuuuuu] (163,-230) circle (1.5pt);
\draw (163,-230-12) node {\footnotesize{$((1,1,1),1)$}};
\draw [fill=uuuuuu] (212,-230) circle (1.5pt);
\draw [fill=uuuuuu] (254,-230) circle (1.5pt);
\draw (254,-230-12) node {$\cdots$};
\draw [fill=uuuuuu] (292,-230) circle (1.5pt);
\draw [fill=uuuuuu] (331,-230) circle (1.5pt);
\draw (331,-230-12) node {\footnotesize{$((1,1,5),1)$}};;
\end{tikzpicture}
\caption{The tree $\B T_G$ of $G$.}
\label{fig:ordered-tree}
\end{figure}

\begin{figure}[htbp]
\centering
\tikzset{every picture/.style={line width=0.75pt}} 

\begin{tikzpicture}[x=0.75pt,y=0.75pt,yscale=-1,xscale=1]

\draw[line width=1pt,color=sqsqsq]   (224,110) .. controls (224,56.98) and (266.98,14) .. (320,14) .. controls (373.02,14) and (416,56.98) .. (416,110) .. controls (416,163.02) and (373.02,206) .. (320,206) .. controls (266.98,206) and (224,163.02) .. (224,110) -- cycle ;
\draw[line width=1pt,color=sqsqsq]   (149,223.5) .. controls (149,198.37) and (169.37,178) .. (194.5,178) .. controls (219.63,178) and (240,198.37) .. (240,223.5) .. controls (240,248.63) and (219.63,269) .. (194.5,269) .. controls (169.37,269) and (149,248.63) .. (149,223.5) -- cycle ;
\draw[line width=1pt,color=sqsqsq]   (402,222.5) .. controls (402,197.37) and (422.37,177) .. (447.5,177) .. controls (472.63,177) and (493,197.37) .. (493,222.5) .. controls (493,247.63) and (472.63,268) .. (447.5,268) .. controls (422.37,268) and (402,247.63) .. (402,222.5) -- cycle ;
\draw[line width=1pt,color=sqsqsq]   (268,78.5) .. controls (268,49.51) and (291.51,26) .. (320.5,26) .. controls (349.49,26) and (373,49.51) .. (373,78.5) .. controls (373,107.49) and (349.49,131) .. (320.5,131) .. controls (291.51,131) and (268,107.49) .. (268,78.5) -- cycle ;
\draw[line width=1pt,color=sqsqsq]   (242,152) .. controls (242,138.19) and (253.19,127) .. (267,127) .. controls (280.81,127) and (292,138.19) .. (292,152) .. controls (292,165.81) and (280.81,177) .. (267,177) .. controls (253.19,177) and (242,165.81) .. (242,152) -- cycle ;
\draw[line width=1pt,color=sqsqsq]   (343,151) .. controls (343,137.19) and (354.19,126) .. (368,126) .. controls (381.81,126) and (393,137.19) .. (393,151) .. controls (393,164.81) and (381.81,176) .. (368,176) .. controls (354.19,176) and (343,164.81) .. (343,151) -- cycle ;
\draw[line width=1pt,color=sqsqsq]   (303,50) .. controls (303,41.16) and (310.16,34) .. (319,34) .. controls (327.84,34) and (335,41.16) .. (335,50) .. controls (335,58.84) and (327.84,66) .. (319,66) .. controls (310.16,66) and (303,58.84) .. (303,50) -- cycle ;
\draw[line width=1pt,color=sqsqsq]   (273,80) .. controls (273,71.16) and (280.16,64) .. (289,64) .. controls (297.84,64) and (305,71.16) .. (305,80) .. controls (305,88.84) and (297.84,96) .. (289,96) .. controls (280.16,96) and (273,88.84) .. (273,80) -- cycle ;
\draw[line width=1pt,color=sqsqsq]   (333,78) .. controls (333,69.16) and (340.16,62) .. (349,62) .. controls (357.84,62) and (365,69.16) .. (365,78) .. controls (365,86.84) and (357.84,94) .. (349,94) .. controls (340.16,94) and (333,86.84) .. (333,78) -- cycle ;
\draw[line width=1pt,color=sqsqsq]   (304,110) .. controls (304,101.16) and (311.16,94) .. (320,94) .. controls (328.84,94) and (336,101.16) .. (336,110) .. controls (336,118.84) and (328.84,126) .. (320,126) .. controls (311.16,126) and (304,118.84) .. (304,110) -- cycle ;


\draw (310,212) node [anchor=north west][inner sep=0.75pt]   [align=left] {$V_1$};
\draw (184,217) node [anchor=north west][inner sep=0.75pt]   [align=left] {$V_2$};
\draw (435,220) node [anchor=north west][inner sep=0.75pt]   [align=left] {$V_3$};
\draw (255,147) node [anchor=north west][inner sep=0.75pt]   [align=left] {$V_{1,2}$};
\draw (352,145) node [anchor=north west][inner sep=0.75pt]   [align=left] {$V_{1,3}$};
\draw (309,136) node [anchor=north west][inner sep=0.75pt]   [align=left] {$V_{1,1}$};


\end{tikzpicture}
\caption{The partition structure of a  $\{P_1, P_2\}$-mixing construction $G$ with exactly three levels: the base for level-1 is $P_1$, while
the bases for the unique recursive parts at level 2 and 3 are $P_1$ and $P_2$ respectively.}
\label{fig:mixing-construction-2}
\end{figure}
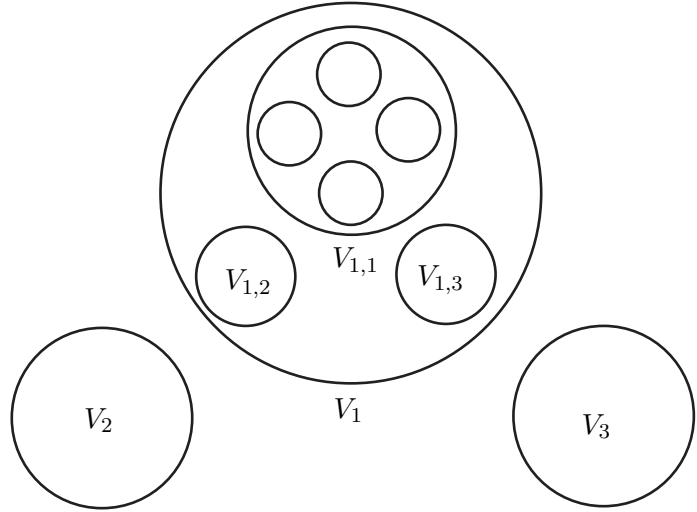

\begin{figure}[htbp]
\centering
\tikzset{every picture/.style={line width=0.75pt}} 
\begin{tikzpicture}[x=0.75pt,y=0.75pt,yscale=-1,xscale=1]

\draw[line width=1pt,color=sqsqsq]    (369,300) -- (468,395) ;
\draw[line width=1pt,color=sqsqsq]    (468,300) -- (468,395) ;
\draw[line width=1pt,color=sqsqsq]    (565,300) -- (468,395) ;
\draw[line width=1pt,color=sqsqsq]    (269,220) -- (369,300) ;
\draw[line width=1pt,color=sqsqsq]    (369,220) -- (369,300) ;
\draw[line width=1pt,color=sqsqsq]    (461,220) -- (369,300) ;
\draw[line width=1pt,color=sqsqsq]    (158,165) -- (269,220) ;
\draw[line width=1pt,color=sqsqsq]    (230,165) -- (269,220) ;
\draw[line width=1pt,color=sqsqsq]    (295,165) -- (269,220) ;
\draw[line width=1pt,color=sqsqsq]    (360,165) -- (269,220) ;

\draw [fill=uuuuuu] (468,395) circle (1.5pt);
\draw (468,395+12) node {\footnotesize{$((),\emptyset)$}};
\draw [fill=uuuuuu] (468,300) circle (1.5pt);
\draw (468,300-12) node {\footnotesize{$((2),1)$}};
\draw [fill=uuuuuu] (565,300) circle (1.5pt);
\draw (565,300-12) node {{\footnotesize{$((3),1)$}}};
\draw [fill=uuuuuu] (369,300) circle (1.5pt);
\draw (369-20,300+12) node {{\footnotesize{$((1),1)$}}};
\draw [fill=uuuuuu] (369,220) circle (1.5pt);
\draw (369,220-12) node {{\footnotesize{$((1,2),1)$}}};
\draw [fill=uuuuuu] (461,220) circle (1.5pt);
\draw (461,220-12) node {{\footnotesize{$((1,3),1)$}}};
\draw [fill=uuuuuu] (269,220) circle (1.5pt);
\draw (269-20,220+12) node {{\footnotesize{$((1,1),1)$}}};
\draw [fill=uuuuuu] (158,165) circle (1.5pt);
\draw (158,165-12) node {{\footnotesize{$((1,1,1),2)$}}};
\draw [fill=uuuuuu] (230,165) circle (1.5pt);
\draw (230,165-12) node {{\footnotesize{$((1,1,2),2)$}}};
\draw [fill=uuuuuu] (295,165) circle (1.5pt);
\draw (295,165-12) node {{\footnotesize{$((1,1,3),2)$}}};
\draw [fill=uuuuuu] (360,165) circle (1.5pt);
\draw (360,165-12) node {{\footnotesize{$((1,1,4),2)$}}};
\draw (16,395) node  {Level-$0\colon$};
\draw (16,300) node {Level-$1\colon$};
\draw (16,220)  node  {Level-$2\colon$};
\draw (16,165) node {Level-$3\colon$};
\end{tikzpicture}
\caption{The tree $\B T_G$ of $G$.}
\label{fig:ordered-tree-2}
\end{figure}

\section{Proof of Theorem~\ref{THM:mixed-pattern}}\label{SEC:Proof-Turan-number}


The main idea of the proof of Theorem~\ref{THM:mixed-pattern}\ref{it:a} is similar to the proof of  Theorem~3 in~\cite{PI14}.
The starting point is the easy observation
(Lemma~\ref{lm:FInfty}) that by forbidding $\C F_\infty$
we restrict ourselves to subgraphs of $P_{I}$-mixing constructions;
thus Part~\ref{it:a} of Theorem~\ref{THM:mixed-pattern} would trivially hold if infinite forbidden families were allowed. Our task is to show that, for some large $M$, the finite subfamily $\C F_M$ of $\C F_\infty$  still has the above properties.
The  Strong Removal Lemma of R\"odl and
Schacht~\cite{RS09} (stated as Lemma~\ref{lm:RS} here)
implies that for every $\e>0$
there is $M$ such that every $\C F_M$-free $r$-graph with
$n\ge M$ vertices can be made $\C F_\infty$-free by removing at most
$\frac{c_0}2 \binom{n}{r}$ edges. It follows that
every maximum $\C F_M$-free $r$-graph $G$ on $[n]$ is $c_0 \binom{n}{r}$-close
in the edit distance to a $P_{I}$-mixing construction (see Lemma~\ref{lm:edit}), where $c_0>0$ can be made arbitrarily small by choosing $M$ large.
Then our key Lemma~\ref{lm:Max2} (which heavily relies on another important result,
the existence of a ``rigid'' $F\in\Sigma P_I$ as proved in Lemma~\ref{lm:MaxRigid}) shows via stability-type arguments that some small constant
$c_0>0$ (independent of $n$) suffices to ensure that there is a partition
$V(G)=V_1\cup \dots\cup V_{m_i}$ for some $i\in I$ such that
$G\setminus (\bigcup_{j\in R_i} G[V_j])=\blow{E}{V_1,\dots,V_{m_i}}$, that is,
$G$ follows exactly the bottom level of some $P_i$-construction
(but nothing is stipulated about what happens inside the recursive parts~$V_j$).
The maximality of $G$ implies that each $G[V_j]$ with $j\in R_i$ is maximum
$\C F_M$-free (see Lemma~\ref{lm:FmFree}), allowing us to apply induction.

Part~\ref{it:b} of Theorem~\ref{THM:mixed-pattern} (which has no direct analogue in~\cite{PI14}) is needed in those applications where we have to analyse almost extremal constructions. It
does not directly follow from Lemma~\ref{lm:RS} (i.e.\ from the Removal Lemma), since the same constant $M$ in Theorem~\ref{THM:mixed-pattern}\ref{it:b} has to work for \textbf{every} $\varepsilon>0$. Similarly to Part~\ref{it:a}, the key idea here that, once we forced our $\C F_M$-free graph $G$ on $[n]$ to be $c_0\binom{n}{r}$-close to a $P_I$-mixing construction for some sufficiently small $c_0>0$ (but independent of $\e$)  then we can further bootstrap this to the required $\e \binom{n}{r}$-closeness by stability-type arguments.

Many simple lemmas that are needed for our proof can be borrowed from~\cite{PI14}, verbatim or with some minor modifications. However, new challenges arise to  accommodate our situation $|I| \ge 2$ and some new ideas are required here.

\subsection{Some additional assumptions and definitions related to $P_I$}\label{se:assumptions}

Recall that $I$ is finite, $P_I=\{P_i \colon i\in I\}$ with each pattern $P_i = (m_i, E_i, R_i)$ for $i\in I$ being minimal. Let us define
\begin{equation}\label{eq:Lambda}
 \lambda\coloneqq \max\{\lambda_{P_i}\colon i\in I\}\qquad\mbox{and}\qquad I'\coloneqq \{i\in I\colon \lambda_{P_i}=\lambda\}.
\end{equation}

We can assume that $\lambda>0$: if $\lambda=0$ then every $P_I$-mixing construction is edgeless and we can satisfy Theorem~\ref{THM:mixed-pattern} by letting $M = r$, with $\C F_{M} = \{K^r_r\}$  consisting of a single edge.

Furthermore, we can assume that $\lambda_{P_i}>0$ for every $i\in I$ by removing all patterns with zero Lagrangian. (This does not change the family of all $\Sigma P_i$-mixing constructions, since $E_i=\emptyset$ for every $i\in I$ with $\lambda_{P_i}=0$ and we always have the option to put the empty $r$-graph into a part.)

Note that if $\lambda = 1$, i.e.\ $\lambda_{P_t} = 1$ for some $t\in I$, then every complete $r$-graph is a $P_t$-construction by Lemma~\ref{lm:density=1}. Indeed, for example, if the second statement of Lemma~\ref{lm:density=1} holds, that is, $\multiset{\rep{i}{r-1},j}\in E_t$ for some $i\in R_t$, then this is attained by making all partitions to consist of only two non-empty parts, the $j$-th part consisting of a single vertex and the $i$-th (recursive) part being the rest. So, if $\lambda = 1$ then $\Lambda_{P_I}(n)=\binom{n}{r}$ and we can satisfy Theorem~\ref{THM:mixed-pattern} by letting $M = 0$  with $\mathcal{F}_{M} = \emptyset$ being the empty forbidden family. 
Therefore, let us assume that every $P_i$, in addition to being minimal, is also proper (that is, that $0<\lambda_{P_i}<1$ for each $i\in I$).

We can additionally assume that, for all distinct $i,j\in [m]$, the patterns $P_i$ and $P_j$ are \textbf{non-isomorphic} (that is, there is no bijection $f:[m_i]\to [m_j]$ with $f(R_i)=R_j$ and $f(E_i)=E_j$), by just keeping one representative from each isomorphism class.

Also, let us state here some definitions that apply in Section~\ref {SEC:Proof-Turan-number} (that is, for the rest of the proof of Theorem~\ref{THM:mixed-pattern}).
Since the set $I$ is finite, we can fix 
a constant $\beta>0$ that satisfies Part~\ref{it:separated}
of Lemma~\ref{lm:Omega1} for each pattern $P_i$ with $i\in I$.
Also, for $i\in I$, let us define
\begin{equation}\label{eq:fi}
f_i(\B x)\coloneqq  \lambda_{E_i}(\B x) + \lambda_{P_i} \sum_{j\in R_i} x_j^r,\quad \B x\in\IS_{m_i},
\end{equation}
 and let $\C X_i\subseteq \IS_{m_i}^*$ be the set of $P_i$-optimal vectors.
Similar to Fact~\ref{FACT:Lagrangian}\ref{FACT:Lagrangian-1}, for every $P_{I}$-mixing construction $G$ on $[n]$ with base $P_i$ and bottom partition $[n] = V_{1} \cup \cdots \cup V_{m_i}$, we have 
\begin{align}\label{equ:Lagrangian}
    \lambda
    \ge f_i\left(|V_1|/n, \ldots, |V_{m_i}|/n\right)
    = \rho(G) + o(1). 
\end{align}

\subsection{Basic properties of $P_I$-mixing constructions}\label{basic}
In this section we present some simple properties of $P_I$-mixing constructions.

The following two lemmas follow easily from the definitions.
We refer the reader to~{\cite[Lemmas~6 and 7]{PI14}}, where the proofs for Lemmas~\ref{lm:constr} and~\ref{lm:FInfty} are provided for $P$-constructions. 

\begin{lemma}\label{lm:constr}
Take any $G \in \Sigma P_{I}$.
Then
\begin{enumerate}[label=(\alph*)]
\item\label{it:constr1} every induced subgraph of $G$ is contained in $\Sigma P_{I}$, and
\item\label{it:constr2} every blowup of $G$ is a $P_I$-mixing subconstruction (that is, a subgraph in some element of $\Sigma P_{I}$).
\end{enumerate}
\end{lemma}

\begin{lemma}\label{lm:FInfty}
The following statements are equivalent for an arbitrary $r$-graph $G$ on at most $n$ vertices:
\begin{enumerate}[label=(\alph*)]
\item\label{it:FInfty1} $G$ is $\C F_n$-free;
\item\label{it:FInfty2} $G$ is $\C F_\infty$-free;
\item\label{it:FInfty3} $G$ is a $P_{I}$-mixing subconstruction.
\end{enumerate}
\end{lemma}

It follows from Lemma~\ref{lm:FInfty} that $\ex(n,\C F_n)=\ex(n,\C F_\infty)=\Lambda_{P_I}(n)$.

\begin{lemma}\label{lm:monotone}
 For all $t\ge s\ge r$, it holds that $\Lambda_{P_I}(s)/\binom{s}{r} \ge \Lambda_{P_I}(t)/\binom{t}{r}$.
\end{lemma}

\bpf Take a maximum $P_I$-mixing construction $G$ on $[t]$. By Lemma~\ref{lm:constr}\ref{it:constr1}, each $s$-set spans at most $\Lambda_{P_I}(s)$ edges. On the other hand, the average number of the edges spanned by a uniformly random $s$-subset of $[t]$ is $\Lambda_{P_I}(t) \binom{t-r}{s-r}/\binom{t}{s}$, giving the required.\epf

The following result states that the maximum asymptotic edge density of a $P_I$-mixing construction is the same as the largest one attained by a single pattern~$P_i$.

\begin{lemma}\label{LEMMA:lambda-sigma}
We have $\lim_{n\to \infty}\Lambda_{P_I}(n)/\binom{n}{r} = \lambda$.
\end{lemma}
\bpf By fixing some $i$ in $I'\not=\emptyset$, we trivially have that $\Lambda_{P_I}(n)\ge \Lambda_{P_i}(n)$ for each $n$. This implies that
\begin{align*}
\lim_{n\to \infty}\Lambda_{P_I}(n)/\binom{n}{r}\ge \lim_{n\to\infty} \Lambda_{P_i}(n)/\binom{n}{r}=\lambda.
\end{align*}

Let us show the converse inequality.
Let $\Sigma_hP_I$ consist of all $P_I$-mixing constructions whose partition structure has height at most~$h$. For example, $\Sigma_1P_I$ consists of all possible blowups of $E_i$, for $i\in I$, without putting any edges into their recursive parts.

Let us show by induction on $h\ge 1$ that $\lambda_h\le \lambda$, where we let
 \begin{align*}
\lambda_{h}\coloneqq \lim_{n\to\infty}\frac{\max\left\{ |G| \colon v(G) = n \text{ and } G\in \Sigma_h P_{I} \right\}}{\binom{n}{r}}.
\end{align*}
(The limit exists since the ratios are non-increasing in $n$ by the same argument as in Lemma~\ref{lm:monotone}.)

If $h=1$, then every $r$-graph in $\Sigma_hP_I$ is a $P_i$-construction for some $i\in I$ and has asymptotic density at most $\max\{\lambda_{P_i}\colon i\in I\}=\lambda$, as desired. Let $h\ge 2$. Take any maximum $G_n\in\Sigma_hP_I$ of order $n\to\infty$. By passing to a subsequence assume that there is $i\in I$ such that each $G_n$ has base pattern~$P_i$. Let  $V(G_n)=V_{n,1}\cup\cdots\cup V_{n,m_i}$ be the base partition of~$G_n$.
Let $x_{n,j} \coloneqq  |V_{n,j}|/n$ for $j\in [m_i]$. By passing to a subsequence again, assume that $x_{n,j}$ tends to some limit $x_j$ for each $j\in[m_i]$.
Then it follows from the definition of $\Sigma_h P_{I}$, continuity of $\lambda_{E_i}$, and induction that
\begin{align}
\rho(G_n)
& \le \lambda_{E_i}(x_{n,1}, \ldots, x_{n,m_i}) + \lambda_{h-1}\sum_{j\in R_i} x_{n,j}^{r} + o(1) \notag\\
& \le \lambda_{E_i}(x_1, \ldots, x_{m_i}) + \lambda\sum_{j\in R_i} x_j^{r} + o(1).  \notag
\end{align}
Furthermore, by Part~\ref{it:f} of Lemma~\ref{lm:Omega1} and by $\lambda_{P_i}\le \lambda$,
we have 
\begin{align*}
 \lambda_{E_i}(x_1, \ldots, x_{m_i}) \le \lambda_{P_i}\left(1-\sum_{j\in R_i} x_j^{r}\right)\le \lambda \left(1-\sum_{j\in R_i} x_j^{r}\right).
 \end{align*}
 By putting these together, we conclude that $\rho(G_n) \le \lambda+o(1)$, giving the required. This finishes the proof of the lemma.
\epf

Given a $P_I$-mixing construction $G$ and a vertex $v$ in $G$, a newly added vertex $v'$ is called a \textbf{clone} of $v$ if the link of $v'$ in the new $r$-graph is identical to the link of $v$ in $G$. 
Note that adding a clone of a vertex to a $P_I$-mixing construction results a $P_I$-mixing subconstruction. 

\begin{lemma}\label{lm:BlowInv}
For every $s\in\I N\cup\{\infty\}$, if an $r$-graph $G$ is $\mathcal{F}_{s}$-free then every blowup of $G$ is $\mathcal{F}_{s}$-free.
\end{lemma}
\bpf
By the definition of blowup, it suffices to show that cloning a vertex of $G$ will not create a copy of any member in $\mathcal{F}_{s}$.
So, let $G'$ be obtained from $G$ by adding a clone $x'$ of some vertex
$x$ of $G$. Take any $U\subseteq V(G')$ with $|U|\le s$. If at least one
of $x$ and $x'$ is not in $U$, then $G'[U]$
is isomorphic to a subgraph of $G$ and cannot be in $\C F_s$; so suppose
otherwise. Since $G$ is $\C F_s$-free, we have that
$G'[U\setminus\{x'\}] = G[U\setminus\{x'\}]$ is a $P_{I}$-mixing subconstruction. By Lemma~\ref{lm:constr}\ref{it:constr2},
$G'[U]$ is a $P_{I}$-mixing subconstruction.
So it follows that $G'$ is
$\C F_s$-free.\epf

\begin{lemma}\label{lm:FmFree}
Let $s\in\I N\cup\{\infty\}$ and $i\in I$. Let $G$ be an $r$-graph on
$V=V_1\cup\dots\cup V_{m_i}$ obtained by taking $\blow{E}{V_1,\dots,V_{m_i}}$
and putting arbitrary $\C F_s$-free $r$-graphs into parts $V_j$ for each $j\in R_i$.
Then $G$ is $\C F_s$-free.
\end{lemma}
\bpf Take an arbitrary $U\subseteq V(G)$ with $|U|\le s$.
Let $U_k\coloneqq V_k\cap U$ for all $k\in [m_i]$.
Notice that $G[U_k]$ has no edges for $k \in [m_i]\setminus R_i$.
On the other hand, for every $k\in R_i$ we have that $G[U_k]$ is a $P_{I}$-mixing subconstruction,
since $|U_k|\le s$ and $G[U_k]\subseteq G[V_k]$ is $\C F_s$-free.
By combining the partition structure of each $G[U_k]$ together with
the level-$1$ decomposition $U=U_1\cup\dots\cup U_{m_i}$, we see
that $G[U]$ is a $P_{I}$-mixing subconstruction.
Therefore, $G$ is $\mathcal{F}_{s}$-free.\epf

The following lemma says that the minimum degree of a maximum $\C F_s$-free $r$-graph is close to its average degree.

\begin{lemma}
\label{lm:regular}
For every $\e>0$ there is $n_0$
such that for every $s\in\I N\cup\{\infty\}$, every maximum $\C F_s$-free
$r$-graph $G$ with $n\ge n_0$ vertices has minimum degree at
least $(\lambda-\e)\binom{n-1}{r-1}$.
\end{lemma}

\bpf Fix $s\in\I N\cup\{\infty\}$. 
The difference between any two vertex degrees in a maximum $\C F_s$-free
$r$-graph $G$ on $n\to\infty$ vertices is at most $\binom{n-2}{r-2}$, as otherwise by deleting one vertex and cloning
the other we can strictly increase the number of edges, while the $r$-graph remains $\C F_s$-free by Lemma~\ref{lm:BlowInv}, a contradiction. Thus each degree is within $O(n^{r-2})$ of the average degree which in turn is at least $(\lambda+o(1))\binom{n-1}{r-1}$ by $\rho(G)+o(1)= \pi(\C F_s)\ge \pi(\C F_\infty)=\lambda$ (the last equality follows from Lemma~\ref{lm:FInfty}).\epf

\subsection{Properties of proper patterns}

Recall that all assumptions made in Section~\ref{se:assumptions} continue to apply to the fixed family $P_I$; in particular, we have $0<\lambda<1$. The following lemma shows that no bottom part of a $P_{I}$-mixing construction $G$ can contain almost all vertices if $G$ has large minimum degree.

\begin{lemma}\label{lm:MaxVi}
For every $c'>0$ there is $n_0$ such that for every $c>c'$ and every $r$-graph $G\in \Sigma P_{I}$ on $n\ge n_0$ vertices with $\delta(G)\ge c\binom{n-1}{r-1}$,
each bottom part $V_j$ of $G$ has at most $(1-c/r)n$ vertices.
\end{lemma}
\bpf Given $c'>0$, let $n\to\infty$ and take any $c$ and $G$ as in the lemma.
Let the base of $G$ be~$P_i$.
If  $j\in [m_i]\setminus R_i$, then it follows from $\delta(G)n/r\le |G|\le \binom{n}{r}-\binom{|V_j|}{r}$ that 
\begin{align*}
    \left(\frac{|V_j|-r}{n}\right)^r
    \le \frac{\binom{|V_j|}{r}}{\binom{n}{r}}
    \le 1-c. 
\end{align*}
Simplifying this inequality, we obtain 
\begin{align*}
    |V_j|
    \le (1-c)^{1/r} n +r
    \le \left(1-\frac{c}{r} - \frac{(r-1)c^2}{2r^2}\right)n + r 
    \le \left(1-\frac{c}{r}\right)n,  
\end{align*}
as desired. Here, we used the inequality $(1-x)^{1/r} \le 1-\frac{x}{r} - \frac{(r-1)x^2}{2r^2}$ for all $r \ge 1$ and $x \in [0,1]$.

Now suppose that $j\in R_i$. Since $V_j \neq V(G)$, pick any vertex $v \in V_k$ with $k\neq j$.
Since $\multiset{\rep{j}{r-1},k} \not\in E_i$ by Lemma~\ref{lm:density=1},
every edge through $v$ contains at least one other vertex outside of~$V_j$.
Thus $c\binom{n-1}{r-1} \le d_G(v)\le (n-|V_j|)\binom{n-2}{r-2}$, implying $|V_j| \le (1-c/(r-1)+o(1))n\le (1-c/r)n$.\epf


Informally speaking, the following lemma (which is a routine generalization of \cite[Lemma~15]{PI14}) implies, among other things,
that all part ratios
of bounded height in a $P_{I}$-mixing construction with large minimum degree
approximately follow some optimal vectors.
In particular,
 for each $i\in I'$, the set $\C X_i$ consists precisely of optimal limiting bottom ratios that
 lead to asymptotically maximum $P_{I}$-mixing constructions with base pattern~$P_i$.
Recall that $\beta>0$ is the constant that satisfies Part~\ref{it:separated} of Lemma~\ref{lm:Omega1} for every $i\in I$ while $\C X_i$ is the set of $P_i$-optimal vectors.

\begin{lemma}
\label{lm:OmegaN}
For every $\e>0$ and $\ell\in\I N$ there are constants
$\alpha_0<\e_0<\dots<\alpha_\ell<\e_\ell<\alpha_{\ell+1}$ in $(0,\e)$
such that the following holds for all sufficiently large~$n$. Let $G$ be
an arbitrary $P_{I}$-mixing construction on $n$ vertices with the partition
structure $\B V$ such that $\delta(G)\ge (\lambda-\alpha_0)\binom{n-1}{r-1}$.
Suppose that $\B i = (i_1, \ldots, i_{s})$ is a legal sequence of length $s$ with $0\le s\le \ell$,
and the induced subgraph $G[V_{\B i}]$ has base $P_j$ for some $j\in I$.
Let $\B v_{\B i}\coloneqq (|V_{\B i,1}|/|V_{\B i}|,\dots, |V_{\B i,m_j}|/|V_{\B i}|)$, where
$V_{\B i} = V_{\B i, 1} \cup \cdots \cup V_{\B i, m_j}$ is the bottom partition of $G[V_{\B i}]$. Then:
 \begin{enumerate}[label=(\alph*)]
  \item\label{it:OmegaI'} $j\in I'$ (that is, $\lambda_{P_j}=\lambda$);
 \item\label{it:OmegaNx}  $\|\B v_{\B i}-\B x\|_\infty\le \e_s$ for some $\B x\in\C X_j$;
 \item\label{it:OmegaNVi} $|V_{\B i,k}|\ge \left(\frac{\beta}{2}\right)^{s+1}n$ and $|V_{\B i,k}|\le \left(1-\frac{\lambda}{{{2r}}}\right)^{s+1}n$ for all $k\in[m_j]$;
\item\label{it:OmegaNDelta} $\delta(G[V_{\B i,k}])\ge (\lambda-\alpha_{s+1})\binom{|V_{\B i,k}|-1}{r-1}$ for all $k\in R_j$.
 \end{enumerate}
\end{lemma}
\bpf We choose positive constants in this order
\begin{align*}
 \alpha_{\ell+1}\gg \e_\ell\gg \alpha_\ell\gg \dots\gg \e_0\gg \alpha_0,
\end{align*}
each being sufficiently small depending on $ P_{I}$, $\e$, $\beta$,
and the previous constants. We use induction on $s=0,1,\dots,\ell$, with the induction step also working for the base case $s=0$, in which $\B i$ is the empty sequence. Take any $s\ge 0$ and suppose that the lemma holds for all smaller values of~$s$.

Let $n$ be large and let $G$ and $\B i$ be as in the lemma.
 Let $U_k\coloneqq V_{\B i,k}$ for $k\in
[m_j]$ and $U\coloneqq V_{\B i}$. Thus $U=U_1\cup\dots\cup U_{m_j}$, and $\B v_{\B
i}=(|U_1|/|U|,\dots,|U_{m_j}|/|U|)$.
Note that
 \begin{equation}\label{eq:deltaU}
 \delta(G[U])\ge (\lambda-\alpha_s)\binom{|U|-1}{r-1},
 \end{equation}
 which holds by Part~\ref{it:OmegaNDelta} of the inductive assumption applied to $G$  if $s>0$, and by the assumption of the lemma if $s=0$ (when $U=V_\emptyset=V(G)$). Thus, we have that
 \begin{equation}\label{eq:G[U]}
 |G[U]|\ge \delta(G[U])|U|/r \ge (\lambda-\alpha_s)\binom{|U|}{r}.
 \end{equation}

Note that $|U|$ can be made arbitrarily large by choosing $n$ large. Indeed, this follows from the  induction assumption for $s-1$ (namely Part~\ref{it:OmegaNVi}) applied to $G$ if $s\ge 1$, while $U=V(G)$ has $n$ elements if $s=0$.

We claim that $\lambda_{P_j}=\lambda$. Suppose on the contrary that the base of $G[U]$ is $j\in I\setminus I'$. Using the asymptotic notation with respect to $|U|\to\infty$, we have by Lemma~\ref{lm:Omega1}\ref{it:f} that
\begin{align*}
\rho(G[U])&\le \lambda_{E_j}(\B v_{\B i})+\lambda \sum_{i\in R_j}\B v_{\B i,i}^r+o(1)\\
&= f_j(\B v_{\B i})+(\lambda-\lambda_{P_j}) \sum_{i\in R_j}\B v_{\B i,i}^r+o(1)\\
 &\le  \lambda_{P_j} + (\lambda-\lambda_{P_j}) \sum_{i\in R_j}\B v_{\B i,i}^r+o(1).
\end{align*}
This is strictly smaller than $\lambda-\alpha_s$ since, by~\eqref{eq:deltaU} and Lemma~\ref{lm:MaxVi}, each $\B v_{\B i,i}$ is at most, say, $1-\lambda/2r$ and thus 
\begin{align*}
    \sum_{i\in R_j}\B v_{\B i,i}^r 
    \le \left(1-\frac{\lambda}{2r}\right) \cdot \sum_{i\in R_j}\B v_{\B i,i}^{r-1}
    \le \left(1-\frac{\lambda}{2r}\right) \cdot \left(\sum_{i\in R_j}\B v_{\B i,i}\right)^{r-1}
    \le 1-\frac{\lambda}{2r}
\end{align*} 
is bounded away from~1. This contradiction shows that $j\in I'$, proving Part~\ref{it:OmegaI'}.

Since $\alpha_s\ll \e_s$, Equation~\eqref{eq:G[U]}, Lemma~\ref{lm:MaxVi} and
Part~\ref{it:almost} of Lemma~\ref{lm:Omega1} (when applied to $P=P_j$ and $\B y=\B v_{\B
i}$)
give the desired $P_j$-optimal vector $\B x$, proving Part~\ref{it:OmegaNx}. Fix one such $\B x=(x_1,\ldots,x_{m_j})$.

For all $k\in[m_j]$, we have $|U_{k}|\ge (x_k-\e_s)|U| \ge (\beta/2)|U|$. This is at least
$(\beta/2)^{s+1}n$ by the inductive assumption on~$|U|$ if $s\ge 1$ (and is trivial if $s=0$). Likewise, $|U| \le \left(1 - \frac{\lambda}{{{2r}}}\right)^{s}n$.
Therefore, it follows from Lemma~\ref{lm:MaxVi} (with $c'\coloneqq \lambda-\alpha_s$) and~\eqref{eq:deltaU} that
\begin{align}
|U_k|
\le \left(1 - \frac{\lambda}{{{2r}}}\right)|U|
\le \left(1 - \frac{\lambda}{{{2r}}}\right)^{s+1}n. \notag
\end{align}
This proves Part~\ref{it:OmegaNVi}.

Finally, take arbitrary $k\in R_j$ and
$y\in U_k$. 
By the definition of Lagrange polynomials and Fact~\ref{FACT:Lagrangian}\ref{FACT:Lagrangian-2}, the degree of $y$ in $\blow{E_j}{U_1,\dots,U_{m_j}}$ is
\begin{align*}
d_{\blow{E_j}{U_1,\dots,U_{m_j}}}(y)
= \left(\lambda_{L_{E_j}(k)}(\B v_{\B
i}) +o(1)\right)\binom{|U|-1}{r-1}
 =\left(\frac1r\cdot \frac{\partial\lambda_{E_j}}{\partial_k}(\B v_{\B
i})+o(1)\right) \binom{|U|-1}{r-1}.
\end{align*}
Since $\|\B v_{\B
i}-\B x\|_\infty\le\e_s\ll \alpha_{s+1}$, we have by Part~\ref{it:f'}
of Lemma~\ref{lm:Omega1} and Lemma~\ref{LEMMA:Taylor-inequ} that, for example,
\begin{align}
\frac{\partial \lambda_{E_j}}{\partial_k}(\B v_{\B
i})-\alpha_{s+1}^2
\le \frac{\partial \lambda_{E_j}}{\partial_k}(\B x)
=\frac{\partial f}{\partial_k} (\B x)- r\lambda x_k^{r-1}
= r\lambda-r\lambda x_k^{r-1}. \notag
\end{align}
Combining this with the previous inequality, we obtain 
\begin{align*}
    d_{\blow{E_j}{U_1,\dots,U_{m_j}}}(y)
    \le \frac{1}{r}\left(r\lambda - r\lambda x_{k}^{r-1} + \alpha_{s+1}^2 +o(1)\right) \binom{|U|-1}{r-1}
    \le \left(\lambda - \lambda x_{k}^{r-1}  +\frac{2 \alpha_{s+1}^2}{r}\right) \binom{|U|-1}{r-1}.
\end{align*}
Thus, by~\eqref{eq:deltaU} and the fact $\|\B v_{\B
i}-\B x\|_\infty\le\e_s$, we have
\begin{align*}
 d_{G[U_k]}(y)
 &=  d_{G[U]}(y)-d_{\blow{E_j}{U_1,\dots,U_{m_j}}}(y)\\
 &\ge \left((\lambda-\alpha_s)- \left(\lambda-\lambda x_k^{r-1}+\frac{2 \alpha_{s+1}^2}{r}\right)\right)\binom{|U|-1}{r-1}\\
 &= \left(\lambda x_k^{r-1}-\alpha_s - \frac{2 \alpha_{s+1}^2}{r}\right)\binom{|U|-1}{r-1}\\
 &\ge  \left(\lambda (\B v_{\B i, k} - \varepsilon_{s})^{r-1}-\alpha_s - \frac{2 \alpha_{s+1}^2}{r}\right)\binom{|U|-1}{r-1}\\
 &\ge  \left(\lambda \B v_{\B i, k}^{r-1} - \lambda (r-1) \B v_{\B i, k}^{r-2} \varepsilon_{s}-\alpha_s - \frac{2 \alpha_{s+1}^2}{r}\right)\binom{|U|-1}{r-1}.
\end{align*}
Since $\B v_{\B i, k} \gg\alpha_{s+1}\gg \e_s\gg \alpha_s$, we have $\lambda (r-1) \B v_{\B i, k}^{r-2} \varepsilon_{s} + \alpha_s + \frac{2 \alpha_{s+1}^2}{r} \le \alpha_{s+1} \B v_{\B i, k}^{r-1}$.
Therefore, the inequality above on $d_{G[U_k]}(y)$ continues as 
\begin{align*}
    d_{G[U_k]}(y)
    \ge \left(\lambda \B v_{\B i, k}^{r-1} - \alpha_{s+1} \B v_{\B i, k}^{r-1}\right)\binom{|U|-1}{r-1}
    \ge \left(\lambda  - \alpha_{s+1} \right)\binom{\B v_{\B i, k} |U|-1}{r-1}
    = \left(\lambda  - \alpha_{s+1} \right)\binom{|U_k|-1}{r-1}.
\end{align*}
This finishes the proof of Part~\ref{it:OmegaNDelta}.\epf

For every $P_{I}$-mixing construction $G$ with base $P_i$,
the following lemma gives a bound for the degree of a vertex in $G$
in terms of the distance between the vector of the part ratios and the set~$\C X_i$.

\begin{lemma}\label{LEMMA:vtx-degree-in-pattern}
Let $G \in \Sigma P_{I}$ be an $r$-graph on $n$ vertices with base $P_{i}$ and
bottom partition $V\coloneqq  V(G) = V_1\cup \cdots \cup V_{m_i}$ for some $i\in I'$. 
Let $\B u \coloneqq  \left(|V_1|/|V|, \ldots, |V_{m_i}|/|V|\right)$ and $\B x\in \C X_i$.
Then for every $j\in [m_i]$ and for every $v\in V_j$ we have
\begin{align}
d_{G}(v)
\le
\begin{cases}
\left(\lambda -  \lambda x_j^{r-1} +rm_i \cdot \|\B u - \B x\|_\infty+o(1)\right) \binom{n-1}{r-1} + d_{G[V_{j}]}(v), & \text{if } j\in R_i, \\
\left(\lambda +rm_i \cdot \|\B u - \B x\|_\infty+o(1)\right) \binom{n-1}{r-1}, & \text{if } j\in [m_i]\setminus R_i.
\end{cases}
 \notag
\end{align}
\end{lemma}
\bpf
First, it is easy to see from the definition of $\Sigma P_{I}$ that
\begin{align}
d_{G}(v)
 =  d_{\blow{E_{i}}{V_{1}, \cdots, V_{m_{i}}}}(v) + d_{G[V_{j}]}(v)
  =  \left(\frac1r\cdot \frac{\partial\lambda_{E_{i}}}{\partial_{j}}(\B u)+o(1)\right) \binom{n-1}{r-1} + d_{G[V_{j}]}(v).  \notag
\end{align}
So it follows from Lemma~\ref{LEMMA:Taylor-inequ} that
\begin{align}
d_{G}(v)
\le \left(\frac1r\cdot \frac{\partial\lambda_{E_{i}}}{\partial_{j}}(\B x) +rm_i \cdot \|\B u - \B x\|_\infty +o(1)\right) \binom{n-1}{r-1} + d_{G[V_{j}]}(v). \notag
\end{align}
If $j\in R_i$, then it follows from Lemma~\ref{lm:Omega1}\ref{it:f'} that
\begin{align}
\frac{\partial\lambda_{E_{i}}}{\partial_{j}}(\B x)
 = \frac{\partial\left(\lambda_{E_{i}}+\lambda \sum_{k\in R_i}x_k^{r}\right)}{\partial_{j}}(\B x) - \frac{\partial\left(\lambda \sum_{k\in R_i}x_k^{r}\right)}{\partial_{j}}(\B x)
 = r \lambda - r \lambda x_j^{r-1}, \notag
\end{align}
and hence,
\begin{align}
d_{G}(v)
& \le
\left( \lambda -  \lambda x_j^{r-1} +rm_i \cdot \|\B u - \B x\|_\infty+o(1)\right) \binom{n-1}{r-1} + d_{G[V_{j}]}(v). \notag
\end{align}
If $j\in [m_i]\setminus R_i$, then $d_{G[V_{j}]}(u) = 0$ and we have by Lemma~\ref{lm:Omega1}\ref{it:f'} that
$\partial\lambda_{E_{i}}(\B x)/\partial_{j}
 = r \lambda$, and hence
\begin{align}
d_{G}(v)
 \le \left(\frac1r\cdot r \lambda +rm_i \cdot \|\B u - \B x\|_\infty+o(1)\right) \binom{n-1}{r-1}
 = \left( \lambda +rm_i \cdot \|\B u - \B x\|_\infty+o(1)\right) \binom{n-1}{r-1}. \notag
\end{align}
This proves Lemma~\ref{LEMMA:vtx-degree-in-pattern}.
\epf

The next lemma shows that every $ P_{I}$-mixing subconstruction with large minimum degree is almost regular.

\begin{lemma}\label{LEMMA:pattern-max-degree}
For every $\varepsilon > 0$ there exist $\delta > 0$ and $n_0$ such that every $ P_{I}$-mixing subconstruction $G$ on $n\ge n_0$ vertices with $\delta(G) \ge (\lambda-\delta)\binom{n-1}{r-1}$
satisfies $\Delta(G) \le (\lambda+\varepsilon)\binom{n-1}{r-1}$.
\end{lemma}
\bpf
Fix any $\varepsilon >0$. Choose a large constant $\ell\in \mathbb{N}$.
Let $\alpha_0<\e_0<\dots<\alpha_\ell<\e_\ell<\alpha_{\ell+1}$ in $(0,\e)$ be the constants given by Lemma~\ref{lm:OmegaN}, where $\ell$ and $\varepsilon$ corresponds to the same constants in both lemmas.
Next, choose sufficiently small constants $\delta \gg 1/n_0$. Let us show that they satisfy the lemma. So pick any $G$ as in the lemma.

By adding edges to given $G$, we see that it is enough to prove the lemma when $G\in \Sigma P_I$.
So suppose that $G$ is a $P_I$-mixing construction on $n\ge n_0$ vertices with partition structure $\B V$ such that $\delta(G) \ge (\lambda-\delta)\binom{n-1}{r-1}$.
Let $u\in V$ and let $\B i \coloneqq (i_1, \ldots, i_s) = \mathrm{br}_{\B V}(u)$ denote the branch of $u$, where $s$ is an integer depending on $u$. 
Let $b_j \coloneqq  b(G[V_{i_1, \ldots, i_{j}}])$ and $\B i_j \coloneqq  (i_1, \ldots, i_j)$ for $0 \le j \le s$.
Notice that $i_{j+1} \in R_{b_j}$ for $0 \le j \le s-1$ and $i_{s} \not\in R_{b_{s-1}}$ (in particular, $d_{G[V_{\B i_{s}}]}(u) = 0$).
Additionally, note that $\B i_0 = ()$ is the empty sequence,~$\B i_s = \B i$, and $V_{\B i_0} = V(G)$. 
Let $\B u_{j} \coloneqq  \left(|V_{\B i_{j-1}, 1}|/|V_{\B i_{j-1}}|, \ldots, |V_{\B i_{j-1}, m_{b_{j-1}}}|/|V_{\B i_{j-1}}|\right)$ for $1\le j \le s$.
Notice from Lemma~\ref{lm:OmegaN}\ref{it:OmegaI'} that $b_{j} \in I'$ for every $1 \le j \le \ell$.
Additionally, by Lemma~\ref{lm:OmegaN}\ref{it:OmegaNx} and our choice of constants, for each $j\in [\ell]$ the vector $\B u_{j}$ is within distance $\e_j$ from some $P_{b_{j-1}}$-optimal vector $\B x_{j}$.
Also, let $m\coloneqq r\cdot \max\{ m_i\colon i\in I\}$. 

By Lemma~\ref{LEMMA:vtx-degree-in-pattern}, we obtain
\begin{align*}
    d_{G}(u) 
    & \le \left(\lambda - \lambda \B x_{1,i_1}^{r-1} + rm \cdot  \|\B u_{1} - \B x_{1}\|_\infty\right) \binom{|V_{\B i_0}|-1}{r-1} + d_{G[V_{\B i_1}]}(u) \\
    & \le \left(\lambda - \lambda \left(\B u_{1,i_1} - \varepsilon_1 \right)^{r-1} + rm  \varepsilon_1\right) \binom{|V_{\B i_0}|-1}{r-1} + d_{G[V_{\B i_1}]}(u) \\
    & \le \left(\lambda - \lambda \B u_{1,i_1}^{r-1} + (r-1) \B u_{1,i_1}^{r-2} \varepsilon_{1} + rm  \varepsilon_1\right) \binom{|V_{\B i_0}|-1}{r-1} + d_{G[V_{\B i_1}]}(u) \\
    & = \lambda \binom{|V_{\B i_0}|-1}{r-1} - \lambda \B u_{1,i_1}^{r-1} \binom{|V_{\B i_0}|-1}{r-1} + \left((r-1) \B u_{1,i_1}^{r-2} \varepsilon_{1} + rm  \varepsilon_1\right) \binom{|V_{\B i_0}|-1}{r-1} + d_{G[V_{\B i_1}]}(u) \\
    & \le \lambda \binom{|V_{\B i_0}|-1}{r-1} - \lambda  \binom{\B u_{1,i_1} |V_{\B i_0}|-1}{r-1} + 2 rm \varepsilon_1 \binom{n-1}{r-1} + d_{G[V_{\B i_1}]}(u)  \\
    & = \lambda \binom{|V_{\B i_0}|-1}{r-1} - \lambda  \binom{|V_{\B i_1}|-1}{r-1} + 2 rm \varepsilon_1 \binom{n-1}{r-1} + d_{G[V_{\B i_1}]}(u).
\end{align*}
If $s < \ell$, then by repeating the argument above for $s$ times, we obtain 
\begin{align*}
    d_{G}(u) 
    & \le \sum_{j=0}^{s-1} \left(\lambda \binom{|V_{\B i_j}|-1}{r-1} - \lambda  \binom{|V_{\B i_{j+1}}|-1}{r-1} + 2 rm \varepsilon_{j+1} \binom{n-1}{r-1} \right) + d_{G[V_{\B i_s}]}(u) \\
    & = \lambda \binom{|V_{\B i_0}|-1}{r-1} - \lambda  \binom{|V_{\B i_{s}}|-1}{r-1} +   \sum_{j=0}^{s-1} \varepsilon_{j+1} \cdot 2 rm \binom{n-1}{r-1} \\
    & \le \lambda \binom{n-1}{r-1} + \sum_{j=0}^{s-1} \varepsilon_{j+1} \cdot 2 rm \binom{n-1}{r-1} 
    \le (\lambda + \varepsilon) \binom{n-1}{r-1}.
\end{align*}
%
If $s \ge \ell$, then by repeating the argument above for $\ell$ times and applying Lemma~\ref{lm:OmegaN}\ref{it:OmegaNVi} to $V_{\B i_{\ell}}$, we obtain 
\begin{align*}
    d_{G}(u) 
    & \le \sum_{j=0}^{\ell-1} \left(\lambda \binom{|V_{\B i_j}|-1}{r-1} - \lambda  \binom{|V_{\B i_{j+1}}|-1}{r-1} + 2 rm \varepsilon_{j+1} \binom{n-1}{r-1} \right) + d_{G[V_{\B i_{\ell}}]}(u) \\
    & = \lambda \binom{|V_{\B i_0}|-1}{r-1} - \lambda  \binom{|V_{\B i_{\ell}}|-1}{r-1}  + \sum_{j=0}^{\ell-1} \varepsilon_{j+1} \cdot 2 rm \binom{n-1}{r-1} + \binom{|V_{\B i_{\ell}}|-1}{r-1} \\
    & \le \lambda \binom{|V_{\B i_0}|-1}{r-1} + \sum_{j=0}^{\ell-1} \varepsilon_{j+1} \cdot 2 rm \binom{n-1}{r-1} + \binom{\left(1-\frac{\lambda}{2r}\right)^{\ell} n - 1}{r-1} \\
    & \le \left(\lambda + \sum_{j=0}^{\ell-1} \varepsilon_{j+1} \cdot 2 rm + \left(1-\frac{\lambda}{2r}\right)^{\ell}\right) \binom{n-1}{r-1} 
    \le (\lambda + \varepsilon) \binom{n-1}{r-1}. 
\end{align*}
This proves Lemma~\ref{LEMMA:pattern-max-degree}.
\epf


Given a collection $E$ of $r$-multisets on a set $V$, a subset $U\subseteq V$ is called \textbf{externally $E$-homogeneous} if
every permutation $\sigma$ of $V$ that fixes every vertex outside of $U$ is a
symmetry of the set of multisets in $E$ that intersect the complement of $U$, that
is, $\sigma\left(E\setminus \multisets{U}{r}\right)=E\setminus \multisets{U}{r}$.
Equivalently, every permutation of $V$ that moves only the elements of $U$ preserves the set of multisets from $E$ that intersect both $U$ and $V\setminus U$.

It is clear from the definition that if $|U| = 1$, then $U$ is always externally $E$-homogeneous.
In addition, we have the following simple fact. 
\begin{fact}\label{FACT:ext-homo}
    If $E$ is an $r$-graph (that is, all multisets in $E$ are simple sets), then a set $U \subseteq  V(E)$ is externally $E$-homogeneous if and only if for every $s \in [r-1]$ and for every $S \in \binom{V\setminus U}{s}$
    \begin{align*}
        \mathrm{either}\quad
        L_{E}(S) \cap \binom{U}{r-s} = \emptyset
        \quad\mathrm{or}\quad
        L_{E}(S) \cap \binom{U}{r-s} = \binom{U}{r-s}. 
    \end{align*}
\end{fact}

Given a pattern $P = (m, E, R)$,  a map $h:[m]\to[m]$ is
an \textbf{automorphism} of the pattern $P$ if $h$ is bijective, $h(R)=R$, and $h$ is
an automorphism of $E$ (that is, $h(E)=E$).
Let us call a $P_{I}$-mixing construction $G$
with the base $P_i$ and the bottom partition $V_1\cup\dots\cup V_{m_i}$ for some $i\in I'$ \textbf{rigid} if,
for every embedding $f$ of $G$ into a $P_{I}$-mixing construction
$H$ with some base $P_j$ and bottom partition $U_1\cup\dots\cup U_{m_j}$ such that $f(V(G))$ intersects at least
two different parts $U_{k}$, we have $j = i$ and there is an automorphism $h$ of $P_i$
such that $f(V_{k})\subseteq U_{h(k)}$
for every $k\in[m_i]$.

%
%

The following key lemma  generalizes \cite[Lemma 17]{PI14} by allowing more than one pattern.
Its proof requires some new ideas. For example, a trick that was used a number of times in~\cite{PI14}, in particular when proving Lemma~17 there, is that any embedding of a maximum
$P$-construction $G$ into another $P$-construction is induced (that is, non-edges are mapped to non-edges). However, a maximum  $P_I$-mixing construction whose base has to be $P_i$ for given $i\in I'$ need not be maximum (nor even maximal) among all $P_I$-mixing constructions and some different arguments are required.

\begin{lemma}\label{lm:MaxRigid}
There exist constants $\varepsilon_0>0$ and $n_0$ such that every  $P_{I}$-mixing construction $G$ on $n\ge n_0$ vertices with $\delta(G) \ge (\lambda-\varepsilon_0)\binom{n-1}{r-1}$ is rigid.
\end{lemma}

%
\bpf
For every $i\in I$ define
\begin{align}
\eta_i \coloneqq  \min\left\{\lambda-\lambda_{P_i-j}\colon j\in [m_i]\right\}. \notag
\end{align}
Since $P_i$ is minimal, we have $\eta_i > 0$.
Let $\eta \coloneqq  \min\left\{\eta_i \colon i\in I\right\}$. Clearly, $\eta>0$.

Recall that $\beta>0$ is a constant satisfying Lemma~\ref{lm:Omega1}\ref{it:separated} for each pattern~$P_i$.
Choose sufficiently small positive constants in this order $\varepsilon_2\gg \varepsilon_1\gg \varepsilon_0$.
Let $n$ be a sufficiently large integer
and $G$ be a $P_{I}$-mixing construction on $[n]$ that satisfies the assumptions in Lemma~\ref{lm:MaxRigid}.
Let $\B V$ denote the partition structure of $G$. Assume that the bottom partition is $V(G) = V_1 \cup \cdots \cup V_{m_i}$ for some~$i\in I$.

By our choice of the constants (and by Lemmas~\ref{lm:OmegaN} and~\ref{LEMMA:pattern-max-degree}) we have that
\begin{enumerate}[label=(\Alph*)]
\item\label{it:C} $|V_k| \ge \beta n/2$ for all $k\in [m_i]$,
\item\label{it:B} $\delta(G[V_k]) \ge (\lambda-\varepsilon_1)\binom{|V_k|-1}{r-1}$ for all $k\in R_i$, and
\item\label{it:A} every $P_{I}$-mixing subconstruction $G'$ on $n$ vertices with $\delta(G') \ge (\lambda-\varepsilon_0)\binom{n-1}{r-1}$ satisfies
$\Delta(G') \le (\lambda+\varepsilon_1)\binom{n-1}{r-1}$.
\end{enumerate}
Take any embedding $f$ of $G$ into some $P_{I}$-mixing construction
$H$ with the base $P_j$ and the bottom partition $V(H)=U_1\cup\dots\cup U_{m_j}$ for some $j\in I$ such that $f(V(G))$ intersects at least two different parts $U_{k}$.

Since $|V_k| \ge \beta n/2 > \varepsilon_2 m_j n$ for all $k\in [m_i]$,
by the Pigeonhole Principle there is a function $h\colon [m_i]\to [m_j]$ such that
\begin{align}\label{eq:h(i)}
|f(V_{k})\cap U_{h(k)}|\ge \varepsilon_2 n,\quad \mbox{for all $k\in [m_i]$}.
\end{align}

\bcl{h}{We can choose $h$ in~\req{h(i)} so that,
additionally, $h(R_i)\subseteq R_j$ and $h$
assumes at least two different values.}
\bcpf
First we prove that we can choose $h$ so that $h(R_i)\subseteq R_j$. This is trivially true if $R_i$ is empty, so suppose otherwise.

We claim that $R_j$ is also non-empty.
Suppose to the contrary that $R_j = \emptyset$.
Let $k\in R_i$.
If there exists $j_0\in [m_j]$ such that $|f(V_k)\cap U_{j_0}|\le \varepsilon_2 n$,
then there exists a subset $V_k'\subseteq V_k$ of size at least $|V_k|-\varepsilon_2 n$
such that $f(V_{k}') \subseteq U\setminus U_{j_0}$.
Since, by~\ref{it:C}, 
\begin{align}\label{eq:RhoGVk}
\rho(G[V_k'])
= |G[V_k']|/\binom{|V_k'|}{r}
& \ge {\left(|G[V_{k}]|-\varepsilon_2 n \binom{|V_k|-1}{r-1}\right)}/{\binom{|V_k|}{r}} \notag \\
& = \rho(G[V_{k}]) - \varepsilon_2 n \cdot \frac{r}{|V_k|}
\ge \rho(G[V_{k}]) - \frac{r \varepsilon_2}{\beta/2},
\end{align}
 which is  at least $\lambda -\frac{\eta}{2}$ by \ref{it:B}, the induced subgraph of $H$ on $f(V_{k}') \subseteq U\setminus U_{j_0}$ has edge density at least $\lambda-\eta/2$.
Since $n$ is large and $R_{j} = \emptyset$, we have $\lambda_{P_j - \{j_0\}} \ge \rho(H[f(V_{k}')]) \ge \lambda-\eta/2 > \lambda -\eta_j$, contradicting our definition of $\eta_j$ (that is, the minimality of~$P_j$).
Therefore, we have that $|f(V_k)\cap U_{s}|\ge \varepsilon_2 n$ for all $s\in [m_j]$.

For every $s\in [m_j]$, let $W_s \coloneqq  U_s\cap f(V_k)$.
Pick some $k'\in [m_i]\setminus\{k\}$ and a vertex $u\in V_{k'}$; note that $V_{k'}\not=\emptyset$ by~\ref{it:C}.
Since $|W_s|\ge \varepsilon_2 n\ge r$ for all $s\in [m_j]$ and the $E_j$-link of every element of $[m_j]$ is non-empty by~\eqref{eq:MinLinks},
there exists an edge $D\in H$ containing $f(u)$ such that $D\setminus \{f(u)\}$ is contained in $\bigcup_{s\in [m_j]}W_s$.
Let $X$ denote the profile of $D$ and $X'$ be the profile of $D\setminus \{f(u)\}$, both with respect to $U_{1}\cup \cdots \cup U_{m_j}$.
Since $H$ is a $P_{I}$-mixing construction,
all $r$-subsets of $V(H)$ with profile $X$ are contained in $H$.
Let $\mathcal{D}$ be the collection of all $r$-subsets $D'\ni f(u)$ such that the profile of $D'\setminus f(u)$ with respect to $W_1,\ldots,W_{m_j}$ is~$X'$.
No element of the set $f^{-1}(\mathcal{D})$ is contained in $G$
because every member in $f^{-1}(\mathcal{D})$ has profile $\multiset{k^{r-1},k'}$ which cannot belong to
$E_i$ by Lemma~\ref{lm:density=1}.
On the other hand, if we add $f^{-1}(\mathcal{D})$ into the edge set of $G$, the new $r$-graph $G'$ is still embedded into $H$ by the same function~$f$.
In other words, $G'$ is a $P_{I}$-mixing subconstruction.
However, the degree of $u$ in $G'$ satisfies
\begin{align}
d_{G'}(u)
\ge d_{G}(u)+ |f^{-1}(\mathcal{D})|
& \ge (\lambda-\varepsilon_0)\binom{n-1}{r-1}+\prod_{s\in [m_j]}\frac{(\varepsilon_2 n)^{{X'(s)}}}{X'(s)!} \notag\\
& \ge (\lambda-\varepsilon_0)\binom{n-1}{r-1} + \frac{\varepsilon_2^{r-1}n^{r-1}}{(r-1)!}
 >  (\lambda+\varepsilon_1)\binom{n-1}{r-1}, \notag
\end{align}
which contradicts \ref{it:A} above.
This proves that $R_j \neq \emptyset$.

Suppose to the contrary that we cannot satisfy the first part of the
claim for some $k \in R_i$, that is, for each $\ell \in R_j$ we have $|f(V_{k})\cap U_{\ell}|<\varepsilon_2 n$.
Since $|f(V_{k})\cap U_{\ell}|<\varepsilon_2 n$ for all $\ell \in R_j$,
there exists a subset $V_k' \subseteq V_k$ of size at least $|V_k| - \varepsilon_2 n |R_j|$ such that $f(V_k') \subseteq U_{[m_j]\setminus R_j} \coloneqq \bigcup_{\ell\in [m_j]\setminus R_j} U_{\ell}$,
that is, the induced subgraph $G[V_{k}']$ is embedded into $H[U_{[m_j]\setminus R_j}]$.
Since $|V_{k}| \ge \beta n/2 \gg \varepsilon_2 n |R_j|$ and $\rho(G[V_{k}]) \ge \lambda - \varepsilon_1$,
we have
\begin{align*}
\rho(G[V_k']) 
= \frac{|G[V_{k}']|}{\binom{|V_{k}'|}{r}}
\ge \frac{|G[V_{k}]| - \varepsilon_2 n |R_j| \binom{|V_k|-1}{r-1}}{\binom{|V_k|}{r}}
& = \rho(G[V_{k}]) - \varepsilon_2 n |R_j| \cdot \frac{r}{|V_{k}|} \\  
& > \lambda - \varepsilon_1 - \frac{r\, |R_j| \varepsilon_2}{\beta/2} > \lambda - \frac{\eta}{2}. 
\end{align*}
This means that there are arbitrarily large $(P_j-R_j)$-constructions of edge density at least
$\lambda - \eta/2$, that is, $\lambda_{P_j-R_j}\ge \lambda - \eta/2 > \lambda - \eta$. This contradicts the minimality of $P_j$ since $R_j\not=\emptyset$.

Let us restrict ourselves to those $h$ with $h(R_i)\subseteq R_j$.
Suppose that we cannot fulfil the second part of the claim.
Then there is $s \in [m_j]$ such that
$|f(V_{k})\cap U_{s}|\ge \varepsilon_2 n$ for every $k \in[m_i]$.
Since $E_i \neq \emptyset$, the induced subgraph $G[f^{-1}(U_{s})]$ is non-empty
(it has at least $\varepsilon_2 n$ vertices from each bottom part of $G$) and is mapped entirely into $U_{s}$.
Thus $s\in R_j$.
Since $f(V(G))$ intersects at least two different bottom parts of $H$,
we can pick some $v\in V_{t}$ for some $t\in [m_i]$ such that $f(v)\in U_{s'}$ and $s' \in [m_j]\setminus \{s\}$.
Fix some $(r-1)$-multiset $X$ in $L_{E_i}(t)$; note that $L_{E_i}(t)\neq \emptyset$ by \req{MinLinks}, that is, by the minimality of~$P_i$.
Take an edge $D\in G$ containing $v$ so that $D\setminus\{v\}$ is a subset of $f^{-1}(U_{s})$ and has profile $X$;
it exists because each bottom part of $G$ contains at least $\varepsilon_2 n\ge r$ vertices of $f^{-1}(U_{s})$.
The $r$-set $f(D)$ is an edge of $H$ as $f$ is an embedding.
However, it has $r-1$ vertices in $U_{s}$ and one vertex in~$U_{s'}$.
Thus the $r$-multiset $\multiset{\rep{s}{r-1},s'}$ belongs to $E_j$.
Since $s\in R_j$, this contradicts Lemma~\ref{lm:density=1}. The claim is proved.\ecpf

\bcl{bijection}{Each $h$ satisfying Claim~\ref{cl:h} is a bijection. Moreover, $j=i$ and $h$ is an automorphism of $P_i$.}
\bcpf
Fix a map $h$ that satisfies Claim~\ref{cl:h}.
First we prove that $h$ is injective.
For every $s\in [m_j]$ let $U_{s}'\coloneqq \bigcup_{t\in h^{-1}(s)} f(V_t)\subseteq V(H)$.
Note that $U_s'=\emptyset$ for $s$ not in the image of $h$.
Let $H'$ be the $P_{I}$-mixing construction
on $f(V(G))$ such that $U_1'\cup\dots\cup U_{m_j}'$ is the bottom partition of $H'$ with base $P_j$ (thus $\blow{E_j}{U_1',\dots,U_{m_j}'} \subseteq H'$),
and,  for each $s\in R_j$, $H'[U_s']$ is the image of the $P_{I}$-mixing construction
$G[f^{-1}(U_s')]$ under the bijection $f$.

We have just defined a new $P_{I}$-mixing construction $H'$
so that each part $V_s$ of $G$ is entirely mapped by $f$ into the
$h(s)$-th part of $H'$, that is,
 \begin{align}\label{eq:UncutUnderH}
 f(V_s)\subseteq U_{h(s)}',\quad\mbox{for every $s\in [m_i]$.}
\end{align}
This $H'$ will be used only for proving Claim~\ref{cl:bijection}.

Let us show first that the same map $f$
is an embedding of $G$ into $H'$.  Take any edge $D\in G$. First, suppose that $f(D)$ intersects at least two different
parts $U_t'$. By~\eqref{eq:UncutUnderH}, $D$ has to be a bottom edge of~$G$. Let $D'\in G$ have the same profile as $D$ with respect to $V_1,\ldots,V_{m_i}$ and satisfy
 \begin{align}\label{eq:D'}
 D'\subseteq \bigcup_{s\in[m_i]} \left( V_s\cap f^{-1}(U_{h(s)}) \right),
 \end{align}
 which is possible because there are at least $\varepsilon_2 n$ vertices available in
each part $V_s\cap f^{-1}(U_{h(s)})$. Since $f(D')\cap U_t= f(D')\cap U_t'$ for all $t\in[m_j]$,
the bottom edge $f(D')$ of $H$ has the same profile $X$ with
respect to the partitions $U_1\cup\dots\cup U_{m_j}$ and
$U_1'\cup\dots\cup U_{m_j}'$. Thus $X\in E_j$.
Next, as each $f(V_s)$ lies entirely inside $U_{h(s)}'$, the
sets $f(D)$ and $f(D')$ have the same profiles with
respect to the partition $U_1'\cup\dots\cup U_{m_j}'$.
Thus $f(D)$ is an edge of $\blow{E_j}{U_1',\dots,U_{m_j}'}$,
as required. It remains to consider the case that $f(D)$ lies inside some part~$U_t'$.
Let $G'\coloneqq G[f^{-1}(U_t')]$. Assume that
$t\in [m_j]\setminus R_j$ for otherwise $f(D)\in f(G')$, which is a subset of $H'$ by the
definition of $H'$. We claim that $G'$ has no edges in this case (obtaining a contradiction to $D\in G'$). Since
$h(R_i)\subseteq R_j$, we have
$h^{-1}(t)\cap R_i=\emptyset$ and $D$ must be a bottom edge of $G$. As before,
we can find an edge $D'\in G$ that satisfies~\req{D'} and has the same
profile as $D$ with respect to $V_1,\dots,V_{m_i}$.
However, $f$ maps this $D'$ inside a non-recursive part $U_t$ of $H$,
a contradiction. Thus $f$ is an embedding of $G$ into $H'$.

Thus, by considering $H'$ instead of $H$ (and without changing
$h$) we have that $f(V_s)\subseteq U_{h(s)}'$ for all $s\in[m_i]$.

Suppose on the contrary to the claim that $|h^{-1}(s)|\ge 2$ for some $s\in [m_j]$.  Let $A\coloneqq h^{-1}(s)$
and $B\coloneqq [m_i]\setminus A$. Since $h$ assumes at least two different values,
the set $B$ is non-empty.

Trivially, the set $U_s'$ is externally $H'$-homogeneous.
We claim that $f^{-1}(U_s')=V_A$ is externally $G$-homogeneous (recall that we denote $V_A\coloneqq \bigcup_{\ell\in A} V_{\ell}$).
Indeed, suppose that $V_A$ is not externally $G$-homogeneous. Then there is $D\in G$ that intersects both $V_A$ and its complement such that for some bijection $\sigma$ of $V(G)$ that moves only points of $V_A$ the $r$-set $D'\coloneqq \sigma(D)$ is not in~$G$. The profile of $D'$ with respect to $V_1,\dots,V_{m_i}$ contains at least two different elements of $[m_i]$, one in $A$ and another in $[m_i]\setminus A$. Let $\C D'$ consist of all $r$-sets that have the same profile with respect to $V_1,\dots,V_{m_i}$ as~$D'$. Since each bottom part $V_k$ has at least $\beta n/2$ vertices, it holds that $|\C D'|\ge (\beta n/2)^r/r!$. Since $D'\not\in G$, no $r$-set in $\C D'$ is an edge of~$G$. With respect to $U_1',\dots,U_{m_j}'$, all $r$-sets in $f(\C D')$ have by~\eqref{eq:UncutUnderH} the same profile as $f(D')$, which in turn is the same as the profile as that of $f(D)$ (because the bijection $f\circ \sigma\circ f^{-1}$ of $V(H')$ moves only elements inside $f(V_A)\subseteq U_s'$). Since $f(D)\in H'$, we must have $f(\C D')\subseteq H'$. Thus we can add $\C D'$ to $G$, keeping it a $P_I$-mixing subconstruction. However, the new $r$-graph has edge density at least $\lambda-\varepsilon_0+(\beta/2)^r$, which is a contradiction to~\ref{it:A} since $\beta \gg \varepsilon_{1} \gg \varepsilon_0\gg 1/n$.

Thus $V_A$ is externally $G$-homogeneous. It follows that $A$ is externally $E_i$-homogeneous (since each $V_t$ has at least $\beta n/2\ge r$ elements).

Next, let us show that $A\cap R_i = \emptyset$. Suppose that this is false. Then fix $i_{\ast} \in A\cap R_i$. 
Let $\hat{G}$ be the $r$-graph obtained from $G$ by replacing $G[V_A]$ with a maximum $P_I$-mixing construction. 
By Fact~\ref{FACT:ext-homo}, the $r$-graph $\hat{G}$ remains a $P_I$-mixing construction, with $A$ playing the role of $i_{\ast}$. 
This implies that $\rho(G[V_A]) \ge \lambda-3\varepsilon_0/\beta^r$, for otherwise we would have 
\begin{align*}
    |\hat{G}|
    = |G| - |G[V_{A}]| + |\hat{G}[V_{A}]| 
    & \ge (\lambda - \varepsilon_{0}) \binom{n}{r} - \left(\lambda- \frac{3\varepsilon_0}{\beta^r}\right) \binom{|V_{A}|}{r} + (\lambda+o(1)) \binom{|V_{A}|}{r} \\
    & \ge (\lambda - \varepsilon_{0}) \binom{n}{r} + \frac{2\varepsilon_0}{\beta^r} \binom{|V_{A}|}{r} \\
    & \ge (\lambda - \varepsilon_{0}) \binom{n}{r} + 2\varepsilon_0 \binom{|V_{A}|/\beta}{r} 
     \ge (\lambda - \varepsilon_{0}) \binom{n}{r} + 2\varepsilon_0 \binom{n}{r}, 
\end{align*}
a contradiction to~\eqref{equ:Lagrangian}. 
Here, the last inequality follows from~\ref{it:C} and the assumption that $|A| \ge 2$.

Recall that $B=[m_j]\setminus A$. Consider the pattern $Q\coloneqq  P_i-B$ obtained from $P_i$ by removing~$B$.
Let $E' \coloneqq  E_i[A]$ and $R' \coloneqq  A\cap R_i$.
Without loss of generality we may assume that $A = [a]$ for some $a\in [m_i-1]$.
Let $\B x\coloneqq (x_1,\dots,x_a)$, where  $x_{k} \coloneqq  |V_k|/|V_A|$ for $k\in [a]$.
Then it follows from $\rho(G[V_A]) \ge \lambda-3\varepsilon_0/\beta^r$ that
the obtained vector $\B x\in \mathbb{S}_{a}$ satisfies that
\begin{align}\label{equ:lambda}
\lambda_{E'}(\B x) + \lambda \sum_{k\in R'}x_{k}^r
= \rho(G[V_A]) + o(1)
\ge \lambda - 4\varepsilon_0/\beta^r.
\end{align}
This inequality does not contradict the minimality of $P_i$ yet, since $G[V_A]$ is not necessarily a $Q$-construction (for $k\in R'$, the $P_I$-mixing construction $G[V_k]$ can use parts indexed by~$B$).
However, if we replace $G[V_k]$ by a maximum $Q$-construction for every $k\in R'$,
then the resulting $r$-graph has edge density $\lambda_{E'}(\B x) + \lambda_{Q} \sum_{k\in R'}x_{k}^r + o(1)$.
By the definition of Lagrangian, we have
\begin{align}\label{equ:lambda-Q}
\lambda_{Q} \ge \lambda_{E'}(\B x) + \lambda_{Q} \sum_{k\in R'}x_{k}^r.
\end{align}
By $|A|\ge2$ and \ref{it:C}, we have that  $x_k\le 1 - \beta/2$ for every $k\in [a]$.
So it follows from $\sum_{k=1}^{a}x_k = 1$ that, say, $\sum_{k=1}^{a}x_k^r \le  1 - \varepsilon_2$.
In particular, $\sum_{k\in R'}x_{k}^r \le  1 - \varepsilon_2$.
Therefore, by~\eqref{equ:lambda} and~\eqref{equ:lambda-Q}, we obtain
\begin{align}
\lambda_Q\left(1-\sum_{k\in R'}x_{k}^r\right) \ge \lambda_{E'}(\B x)  \ge  \lambda \left(1-\sum_{k\in R'}x_{k}^r\right) - 4\varepsilon_0/\beta^r, \notag
\end{align}
which implies that $\lambda_{Q} \ge \lambda-\eta/2$, contradicting the minimality of $P_i$.
Therefore, $A\cap R_i= \emptyset$.

Since $|A| \ge 2$, we can choose $t_1, t_2 \in A$ such that $t_1 \neq t_2$.
Since $A$ is externally $E_i$-homogeneous and, by Lemma~\ref{lm:DNs},  we have $L_{E_i}(t_1) \not= L_{E_i}(t_2)$, there exists a multiset in $E_i$ that is completely contained inside~$A$.

So, $G[V_A] \neq \emptyset$.
Since $U_s' = f(V_A)$, we have $H'[U_s'] \neq \emptyset$, and hence, $s\in R_j$.
Notice that $H'$ is a $P_{I}$-mixing construction on $n$ vertices with $\delta(H') \ge \delta(f(G))\ge (\lambda-\varepsilon_0)\binom{n-1}{r-1}$.
So, similarly to \ref{it:B}, we have $\rho(H'[U_k']) \ge \lambda - \varepsilon_1$ for all $k\in R_j$.
In particular, we have $\rho(H'[U_s']) \ge \lambda - \varepsilon_1$.
Therefore, $|H'\setminus H'[U_s']| \le (\lambda+\varepsilon_1)\binom{n}{r} - (\lambda - \varepsilon_1)\binom{|U_s'|}{r}$,
which implies that 
\begin{align}
|G\setminus G[V_A]|
\le |H'\setminus H'[U_s']|
& \le (\lambda+\varepsilon_1)\binom{n}{r} - (\lambda - \varepsilon_1)\binom{|U_s'|}{r}. \notag
\end{align}
Note that $|U_s'|=|V_A|$ since $f$ gives a bijection between these two sets. Therefore,
\begin{align}
|G[V_A]|
= |G| - |G\setminus G[V_A]|
& \ge (\lambda-\varepsilon_0)\binom{n}{r}-\left((\lambda+\varepsilon_1)\binom{n}{r} - (\lambda - \varepsilon_1)\binom{|V_A|}{r}\right) \notag\\
& \ge \lambda \binom{|V_A|}{r} - 3\varepsilon_1 \binom{n}{r}. \notag
\end{align}
Thus,
\begin{align}
\rho(G[V_A])
= \frac{|G[V_{A}]|}{\binom{|V_{A}|}{r}}
\ge \frac{\lambda \binom{|V_A|}{r} - 3\varepsilon_1 \binom{n}{r}}{\binom{|V_A|}{r}} 
\ge \lambda - \frac{3\varepsilon_1 \binom{n}{r}}{\binom{\beta n}{r}} 
= \lambda - \frac{3\varepsilon_1}{\beta^r + o(1)}
\ge \lambda - \frac{\eta}{2}. \notag
\end{align}
This implies that $\lambda_{P_i-B}\ge \rho(G[V_A]) \ge \lambda - \eta/2 > \lambda -\eta$, which contradicts the minimality of~$P_i$.
Therefore, $h$ is injective.

Since each multiset in $E_i$ corresponds to a non-empty set of edges of $G$ by \ref{it:C} and $f$ is an embedding, we have $h(E_i)\subseteq E_j$.

Suppose to the contrary that $h$ is not surjective, that is $m_j \ge m_i + 1$.
Let $\B x = (x_1, \ldots, x_{m_i}) \in \C X_i$ be a $P_i$-optimal vector.
Let $\B y \coloneqq  (x_{h^{-1}(1)}, \ldots, x_{h^{-1}(m_j)})$, where $x_{\emptyset}$ means 0. Thus,
$\B y$ is defined by
rearranging the entries of $\B x$ according to $h$ and padding them with $m_j-m_i$ zeros. By
$ h(R_i) \subseteq R_j$ and $h(E_i) \subseteq E_j$,
we have
\begin{align*}\lambda_{E_j}(\B y) + \lambda \sum_{k\in R_j}y_k^{r} \ge \lambda_{E_i}(\B x) + \lambda \sum_{k\in R_i}x_k^r = \lambda.\end{align*}
On the other hand, by Lemma~\ref{lm:Omega1}\ref{it:f}, we have $\lambda_{E_j}(\B y) + \lambda \sum_{k\in R_j}y_k^{r} \le \lambda$.
Therefore, $\lambda_{E_j}(\B y) + \lambda \sum_{k\in R_j}y_k^{r} = \lambda$, which implies that $\B y \in \C X_j$ is a $P_j$-optimal vector.
However, $\B y$ has at least one entry 0, which contradicts Lemma~\ref{lm:Omega1}\ref{it:boundary}.
Therefore, $h$ is surjective, and hence, is bijective.

None of the inclusions $h(E_i)\subseteq E_j$ and $h(R_i)\subseteq R_j$ can be strict as otherwise, since every $P_j$-optimal vector has all coordinates positive by the minimality of $P_j$, we have $\lambda_{P_j}>\lambda_{P_i}=\lambda$, a contradiction. Thus $h$ gives an isomorphism between $P_i$ and $P_j$, and we conclude that $j=i$. This finishes the proof of the claim.\ecpf

By
relabelling the parts of $H$, we can assume
for notational convenience that $h$ is the identity mapping.
Now we are ready to prove the lemma,
namely that $f(V_s)\subseteq U_s$ for every $s\in [m_i]$.

Suppose on the contrary that $f(v)\in U_t$ for some $v\in V_s$ and $t\in [m_i]\setminus\{s\}$.
It follows that $L_{E_i}(s) \subseteq L_{E_i}(t)$.
By Lemma~\ref{lm:DNs} this inclusion is strict and $s\in R_i$.
Pick some $X$ from $L_{E_i}(t)\setminus L_{E_i}(s) \neq \emptyset$.
For every $D\in H$ containing $f(v)$ such that $D\setminus\{f(v)\}$ has the profile $X$ with respect to both $U_1\cup \dots\cup U_{m_i}$
and $f(V_1)\cup\dots\cup f(V_{m_i})$, we add $f^{-1}(D)$ into $G$. Denote by $\tilde{G}$ the new $r$-graph.
Observe that $f$ is also an embedding of $\tilde{G}$ into $H$. Thus,  
$\tilde{G}$ is a $P_I$-mixing subconstruction.
Notice that since $|f(V_j)\cap U_j| \ge \varepsilon_2 n$, there are at least $\prod_{j\in [m_i]}(\varepsilon_2 n)^{X(j)}/X(j)! \ge \varepsilon_2^{r-1} n^{r-1}/(r-1)!$ such edges~$D$.
Since all edges in $f^{-1}(D)$ contain $v$,
we have $d_{\tilde{G}}(v) \ge d_{G}(v) + \varepsilon_2^{r-1} n^{r-1}/(r-1)! > (\lambda+\varepsilon_1)\binom{n-1}{r-1}$ contradicting \ref{it:A} above.
This shows that $G$ is rigid.\epf

Recall that $I'$ consists of the elements $P_i\in I$ whose Lagragian $\lambda_{P_i}$ attains the maximum value $\lambda=\lambda_{P_I}$. Let us denote
\begin{align*}
P_{I'}\coloneqq \{P_i\mid i\in I'\}.
\end{align*}
Later, in the proof of Lemma~\ref{lm:Max2}, we will need, for a given feasible tree $\B T$ satisfying some technical conditions, the existence of a rigid $P_{I}$-mixing construction $F$ whose tree is $\B T$ and every part of $F$ is sufficiently large, specifically 
\begin{align}\label{eq:RigidLargeVi}
 |V_{\B i}|\ge (r-1)\max\left\{r, \max\{m_k \colon k\in I\}\right\},\quad \mbox{for every legal (with respect to $F$) sequence $\B i$.}
\end{align}

The following two lemmas provide such $F$ (in fact, each obtained $F$ will be a $P_{I'}$-mixing construction). The proofs are slightly different depending on whether $\B T$ is extendable or not. (Recall that a feasible tree $\B T$ is extendable if it is a subtree of some strictly larger feasible tree.)

Recall the definition of ``clone" from the paragraph above Lemma~\ref{lm:BlowInv}.
Note that if we add a clone $v'$ of a vertex $v$ of a $P_I$-mixing construction $G$, we generally obtain a subconstruction rather than a $P_I$-mixing construction. 
For example, if $V_j$ is the bottom part containing $v$ while some edge $D$ of the base multiset $E_i$ contains $j$ with multiplicity more than 1 then the blowup of $D$ in a $P_I$-mixing construction would additionally include some edges containing both $v$ and~$v'$. So let us define the operation of \textbf{doubling} $v$ in $G$ where we take a new vertex $v'$ (called the \textbf{double} of $v$), put it in the partition structure of $G$ so that $v'$ has the same branch as $v$, and add all edges through $v'$ as stipulated by the new partition structure. Of course, the degree of the new vertex $v'$ is at least the degree of $v$ in the old $r$-graph~$G$.

\begin{lemma}\label{LEMMA:non-extendable-partition}
For every non-extendable feasible tree $\B T$ with all indices in $I'$,
there exists a rigid $P_{I}$-mixing construction $F$ such that $\B T_F=\B T$ and~\eqref{eq:RigidLargeVi} holds. 
\end{lemma}
\bpf
Take any tree $\B T$ as in the lemma. Let $\varepsilon_0>0$ and $n_0$ be the constants given by Lemma~\ref{lm:MaxRigid}.
Let $n$ be a sufficiently large integer; in particular, so that $n \ge n_0$ and we can apply Lemma~\ref{lm:OmegaN} with $\ell$ equal to the height of $\B T$.
Let $F$ be a maximum $n$-vertex $P_{I'}$-mixing construction under the requirement that its tree $\B T_F$ is a subtree of~$\B T$.
By taking (near) optimal parts ratio for the bottom partition and for
 every recursive part in $\B T$, it is easy to see that $\rho(F) \ge \lambda-\varepsilon_0/2$. 
 
 We claim that 
 \begin{align}\label{eq:FMinDeg}
 \delta(F)\ge (\lambda-\varepsilon_0)\binom{n-1}{r-1}.
 \end{align}
 Indeed, suppose that a vertex $u$ violates~\eqref{eq:FMinDeg}. Remove $u$ from $F$ and double a vertex $v$ with degree at least $(\lambda-\varepsilon_0/2)\binom{n-1}{r-1}$ in~$F$. The new $r$-graph $F'$ has strictly larger number of edges:
\begin{align}\label{eq:F'-F}
|F'| - |F| \ge (\lambda-\varepsilon_0/2)\binom{n-1}{r-1} - (\lambda-\varepsilon_0)\binom{n-1}{r-1} - \binom{n-2}{r-2} >0,
\end{align}
However, the tree of $F'$ is still a subtree of $\B T$ (even if the part that contained $u$ became edgeless),
contradicting the maximality of~$F$.

It follows from Lemma~\ref{lm:MaxRigid} that $F$ is rigid. Also,~\eqref{eq:RigidLargeVi} holds by~\eqref{eq:FMinDeg} and Lemma~\ref{lm:OmegaN}\ref{it:OmegaNVi}. In particular, we have that $\B T_{F}=\B T$, finishing the proof.
\epf

Call a feasible tree $\B T$ of height $\ell$ \textbf{maximal} if every leaf of height less than $\ell$ is non-recursive (or, equivalently, $\B T$ cannot be extended to a larger feasible tree of the same height).

\begin{lemma}\label{lm:tight}
There exists constant $\ell_0 \in \mathbb{N}$ such that, for every feasible extendable tree $\B T$ with all indices in $I'$ which is maximal of height $\ell_0$,
there exists a rigid $ P_{I}$-mixing construction $F$ such that $\B T_F=\B T$ and~\eqref{eq:RigidLargeVi} holds.
\end{lemma}
\bpf Let $\varepsilon_0$ and $n_0$ be given by Lemma~\ref{lm:MaxRigid}. Fix a sufficiently large $\ell_0$. Let $\B T$ be a tree as in the lemma.
Choose $n \in \mathbb{N}$ to be sufficiently large such that, in particular, $n\ge n_0$, $(\beta/2)^{\ell_0}n \ge (r-1)\max\left\{r, \max\{m_k \colon k\in I\}\right\}$, and we can apply Lemma~\ref{lm:OmegaN} with $\ell$ equal to $\ell_0$.
Let $G$ be a maximum $P_{I'}$-mixing construction provided that $\B T_G\TruncatedLevel{\ell_0}$, the tree $\B T_G$ restricted to levels up to $\ell_0$, is a subtree of $\B T$.
Let $\B V$ denote the partition structure of $G$.

\bcl{mindegreeG}{We have $\delta(G) \ge (\lambda-\varepsilon_0/2)\binom{n-1}{r-1}$.}
\bcpf
If we take (near) optimal parts ratio for all partitions up to level $\ell_0$ and and put a maximum $P_{I'}$-mixing construction into each part corresponding to a recursive leaf of $\B T$, then the obtained $r$-graph $G'$ has edge density $\lambda+o(1)$ as $n\to\infty$.
Since $n$ is sufficiently large and $|G|\ge |G'|$ by the definition of $G$, we can assume that $\rho(G) \ge \lambda-\varepsilon_0/4$.
Now, there cannot be a vertex $v\in V(G)$ with $d_G(v) < (\lambda-\varepsilon_0/2)\binom{n-1}{r-1}$,
for otherwise we would remove $v$ from $G$ and double a vertex $u\in V(G)$ with maximum degree (which is at least $\rho(G)\binom{n}{r}r/n \ge (\lambda-\varepsilon_0/4)\binom{n-1}{r-1}$), getting a contradiction as in~\eqref{eq:F'-F}.
\ecpf

Let $F$ be obtained from $G$ by removing edges in $G[\B V_{\B i}]$ for all legal sequences $\B i$ in $G$ of length at least~$\ell_0+1$.

\bcl{mindegreeF}{We have $\delta(F) \ge (\lambda-\varepsilon_0)\binom{n-1}{r-1}$.}
\bcpf
By Lemma~\ref{lm:OmegaN}\ref{it:OmegaNVi} we have $|V_{\B i}| \le \left(1-\frac{\lambda}{{{2r}}}\right)^{\ell_0}n \ll \varepsilon_0 n$ for all legal sequences $\B i$ in $G$ of length at least $\ell_0$.
So it holds that $d_{F}(u) = d_{G}(u) \ge (\lambda-\varepsilon_0)\binom{n-1}{r-1}$ if the length of $\mathrm{br}_{\B V'}(u)$ is at most $\ell_0$,
and $d_{F}(u) \ge  d_{G}(u) - \left(1-\frac{\lambda}{{{2r}}}\right)^{\ell_0\cdot (r-1)}n^{r-1} \ge (\lambda-\varepsilon_0)\binom{n-1}{r-1}$ if the length of $\mathrm{br}_{\B V'}(u)$ is at least $\ell_0+1$.
\ecpf

Since $\B T$ is maximal of height $\ell_0$, the tree of $F$ is a subtree of~$\B T$.
Our choice of $n$ also makes sure that, for every legal $\B i$ in $\B T$ of length at most $\ell_0$, we have $|V_{\B i}| \ge (\beta/2)^{\ell_0}n \ge (r-1)\max\left\{r, \max\{m_k \colon k\in I\}\right\}$. In particular, we have $\B T_{F}=\B T$.
Since $F$ is a $ P_{I}$-mixing construction with $\delta(F) \ge (\lambda-\varepsilon_0)\binom{n-1}{r-1}$,
it follows from Lemma~\ref{lm:MaxRigid} that $F$ is rigid, finishing the proof of Lemma~\ref{lm:tight}.
\epf

\subsection{Key lemmas}
The following lemma, which is proved via  stability-type arguments, is the key for the proof of Theorem~\ref{THM:mixed-pattern}. Its conclusion, informally speaking, implies there is a way to replace a part of a ``reasonably good" $\C F_{M_0}$-free $r$-graph $G$ by a blowup $G'$ of some $E_i$ on $V(G)$ so that  
the new $r$-graph is still $\C F_{M_0}$-free and satisfies $|G'|\ge |G|+\Omega(|G\bigtriangleup G'|)$ . This is closely related to the so-called \textbf{local} (or \textbf{perfect}) stability, see e.g.~\cite{NorinYepremyan17,PikhurkoSliacanTyros19}.

\begin{lemma}\label{lm:Max2}
There are $c_0>0$ and $M_0\in\I N$ such that the following holds.
Suppose that $G$ is a  $\C F_{M_0}$-free $r$-graph on $n\ge r$ vertices that is
$c_0\binom{n}{r}$-close to some $P_{I}$-mixing construction and satisfies
 \begin{align}\label{eq:DeltaG}
 \delta(G)\ge (\lambda-c_0) \binom{n-1}{r-1}.
 \end{align}
Then there are $i\in I$ and a partition
$V(G)=V_1\cup\dots\cup V_{m_i}$ such that $|V_{j}| \in [\beta n/2, (1-\lambda/2r)n]$ for all $j\in [m_i]$ and
\begin{align}\label{eq:Bad-vs-Missing}
10\, |B|\le |A|,
\end{align}
where 
 \begin{align}
 A&\coloneqq  \blow{E_i}{V_1,\dots,V_{m_i}}\setminus G,\label{eq:A}\\
 B&\coloneqq  G\setminus\Bigg(  \blow{E_i}{V_1,\dots,V_{m_i}}\cup \bigcup_{j\in R_i} G[V_j]\Bigg).\label{eq:B}
\end{align}
\end{lemma}

\bpf Clearly, it is enough to establish the existence of $M_0$
such
that the conclusion of the lemma holds for every sufficiently large $n$. (Indeed,
it clearly holds for $n\le M_0$ by Lemma~\ref{lm:FInfty}, so we can simply increase $M_0$
at the end to take care of finitely many exceptions; alternatively, one
can decrease $c_0$.)

Let $\ell_0$ be the constant returned by Lemma~\ref{lm:tight}. Then let $M_0$ be sufficiently large. Next, choose some constants $c_i$ in this order $c_4\gg c_3\gg c_2\gg
c_1\gg c_0>0$, each being sufficiently small depending on the previous ones. Let
$n$ tend to infinity.

Let $G$ be a
$\C F_{M_0}$-free $r$-graph on $[n]$ that satisfies~\eqref{eq:DeltaG} and is $c_0\binom{n}{r}$-close to some $P_{I}$-mixing construction~$H$.
We can assume that the vertices of $H$ are already
labelled so that $|G\bigtriangleup H|\le c_0\binom{n}{r}$.
Let $\B V$ be the partition structure of $H$.
In particular, the bottom partition of $H$ is $V_1\cup\dots\cup V_{m_i}$ for some $i\in I$.

One of the technical difficulties that we are going to face is that
some part $V_j$ with $j\in R_i$ may in principle contain almost every vertex of $V(G)$
(so every other part $V_{k}$ has $o(n)$ vertices).
This means that the ``real'' approximation to $G$ starts only at some
higher level inside $V_j$. On the other hand, Lemma~\ref{lm:OmegaN}
gives us a way to rule out such cases: we have to ensure that
the minimal degree of $H$ is close to $\lambda \binom{n-1}{r-1}$. So, as our first step, we are going to modify the
$P_{I}$-mixing construction $H$
(perhaps at the expense of increasing
$|G\bigtriangleup H|$ slightly) so that its minimal degree is large. When we change $H$ here (or later), we update the partition structure $\B V$ of $H$ appropriately. (Note that the partition tree $\B T_{H}$ may change too, when some part $V_{\B i}$ stops spanning any edges.)

So, let $Z\coloneqq \{x\in [n]\mid d_H(x)< (\lambda-c_1) \binom{n-1}{r-1}\}$.
Due to~\req{DeltaG}, every vertex of $Z$ contributes at least
$(c_1 - c_0) \binom{n-1}{r-1} \ge c_1\binom{n-1}{r-1}/2 \ge  c_1\binom{n-1}{r-1}/r$ to $|G\bigtriangleup H|$. We conclude that
\begin{align*}
    |Z|
    \le \frac{|G\bigtriangleup H|}{c_1\binom{n-1}{r-1}/r} 
    \le \frac{c_0\binom{n}{r}}{c_1\binom{n-1}{r-1}/r} 
    = \frac{c_0 n}{c_1}.
\end{align*}
Fix an arbitrary $y\in [n]\setminus Z$.
Let us change $H$ by making all vertices in
$Z$ into doubles of $y$. Clearly, we have now
\begin{align}\label{eq:MinDegStage1}
    \delta(H)
     \ge (\lambda-c_1)\binom{n-1}{r-1}-|Z|\binom{n-2}{r-2}
     & \ge (\lambda-c_1)\binom{n-1}{r-1} - \frac{c_0 n}{c_1} \cdot \frac{r-1}{n-1} \binom{n-1}{r-1} \notag \\
     & \ge (\lambda-2c_1)\binom{n-1}{r-1}
\end{align}
 while $|G\bigtriangleup H|\le c_0\binom{n}{r}+ |Z|\binom{n-1}{r-1}\le c_0\binom{n}{r} + \frac{c_0 n}{c_1} \cdot \frac{r}{n} \binom{n}{r}  \le
c_1\binom{n}{r}$.

%


By Lemma~\ref{lm:OmegaN} we can conclude that, in the new $P_I$-mixing construction $H$ satisfying~\eqref{eq:MinDegStage1}, part ratios up to height $\ell_0$ are close to optimal ones and $|V_{\B i}|\ge 2c_4n$ for each legal sequence $\B i$ of length at most~$\ell_0$.

In order to satisfy the lemma, we may need to modify  the current partition $V_1,\dots,V_{m_i}$ of~$V(G)$ further. It will be convenient now to keep track of the sets $A$ and $B$ defined by~\eqref{eq:A} and~\eqref{eq:B} respectively, updating them when the partition changes. Recall that the set $A$  consists of edges that are in $\blow{E_i}{V_1,\dots,V_{m_i}}$ but not in~$G$. Let us call these edges \textbf{absent}. Call an $r$-multiset $D$ on $[m_i]$ \textbf{bad} if $D\not\in E_i$
and $D\not=\multiset{\rep{j}{r}}$ for some $j\in R_i$. Call an edge of $G$ \textbf{bad} if
its profile with respect to $V_1,\dots,V_{m_i}$ is bad. Thus $B$ is precisely the set of bad edges and our aim is to prove that there are at least 10 times more absent edges than bad edges. 

Our next modification is needed to ensure later that~\req{Bxi} holds. Roughly speaking, we want
a property that the number of bad edges cannot be decreased much if
we move one vertex between parts.
Unfortunately, we cannot just take a partition structure $\B V$ that minimises $|B|$ because then we
do not know how to guarantee that~\eqref{eq:MinDegStage1} holds (another property important in our proof).
Nonetheless, we can simultaneously satisfy both properties, although with weaker bounds.

Namely, we modify $H$
as follows (updating $A$, $B$, $\B V$, etc, as we proceed).
If there is a vertex $x\in [n]$ such that by moving it to
another part $V_j$ we decrease $|B|$ by at least
$c_2\binom{n-1}{r-1}$, then we pick $y\in V_j$ of maximum $H$-degree
and make $x$ a double of~$y$.
Clearly, we perform this operation at most
$c_1\binom{n}{r}/c_2\binom{n-1}{r-1}=c_1n/(c_2r)$ times because we initially
had $|B|\le |G\bigtriangleup H|\le c_1\binom{n}{r}$. Thus, we have
at all steps
of this process (which affects at most $c_1n/(c_2r)$ vertices of $H$) that,
trivially,
 \begin{align}
  |V_{\B j}| &\ge 2c_4n -\frac{c_1n}{c_2r}\ \ge\ c_4n,\qquad
\mbox{for all legal
$\B j$ with $|\B j|\le \ell_0$},\label{eq:PartSizes}\\
 |G\bigtriangleup H|&\le c_1\binom{n}{r} + \frac{c_1n}{c_2r} \binom{n-1}{r-1}
\ \le\ c_2 \binom{n}{r}.\label{eq:M0End}
 \end{align}
 It follows that at every step each part $V_j$ had a vertex
of degree at least $(\lambda-c_2/2)\binom{n-1}{r-1}$
for otherwise, by \req{DeltaG} and \req{PartSizes}, the edit distance between $H$ and $G$ at that moment would be at least
\begin{align*}
    \frac{1}{r} \cdot |V_j| \cdot \left(\left(\lambda-c_0\right) \binom{n-1}{r-1} - \left(\lambda-\frac{c_2}{2}\right)\binom{n-1}{r-1} \right)
    \ge \frac{1}{r} \cdot c_4n  \cdot \frac{c_2}{3} \binom{n-1}{r-1}
    = \frac{c_2 c_4}{3} \binom{n}{r}, 
\end{align*}
contradicting the first inequality in \req{M0End}. This implies
that every time we double a vertex it has a high degree. Thus  we have
by \req{MinDegStage1} that, additionally to \req{PartSizes} and~\req{M0End},
the following holds at the end of this process:
 \begin{equation}
 \delta(H) \ge \big(\lambda-\max\{3c_1,c_2/2\}\big)
\binom{n-1}{r-1} - \frac{c_1n}{c_2r} \binom{n-2}{r-2}\ \ge\
(\lambda-c_2) \binom{n-1}{r-1}.\label{eq:DeltaH1}
 \end{equation}
So by Lemma~\ref{lm:OmegaN}\ref{it:OmegaNVi}, $|V_j| \ge \beta n/2$ and $|V_j| \le (1-\lambda/2r)n$
(note that we may choose $c_2>0$ small enough and $n$ large enough in the beginning so that Lemma~\ref{lm:OmegaN}\ref{it:OmegaNVi} applies).


Suppose that $B\not=\emptyset$ for otherwise the lemma holds trivially.

Let $H'$ be obtained from $H$ by removing edges contained in $H[\B V_{\B j}]$ for all legal $\B j$ of length at least $\ell_0+1$. 
Let $\B T\coloneqq \B T_{H'}$. 
It follows from the definition that $\B T=\B T_{H}\TruncatedLevel{\ell_0}$,  that is, $\B T$ is obtained from $\B T_H$ by restricting it to levels up to~$\ell_0$. By Lemma~\ref{lm:OmegaN}\ref{it:OmegaI'}, all indices in $\B T$ belong to~$I'$.
If $\B T$ is non-extendable (and thus $\B T=\B T_H$), then we let $F$ be the rigid construction given by by Lemma~\ref{LEMMA:non-extendable-partition} whose tree $\B T_F$ is the same as~$\B T$.
If $\B T$ is extendable then, due to Lemma~\ref{lm:OmegaN}\ref{it:OmegaNDelta} and $\ell_0\ll 1/c_1\ll n$, every recursive part $V_{\B i}$ with $|\B i|< \ell_0$ spans at least one edge in $H'$ and thus the tree $\B T$ is maximal of height~$\ell_0$.
In this case, we let $F$ be the rigid construction given by Lemma~\ref{lm:tight} whose tree is~$\B T$.
In either case, let $\B W=(W_{\B i})$ be the partition structure of $F$. 
Since the number of possible trees~$\B T$ is bounded by a function of $\ell_0$ and we have $\ell_0\ll M_0$, we can assume that
 \begin{equation}\label{eq:M0}
 M_0\ge v(F)+r.
 \end{equation}

Let us show that the maximal degree of $B$ is small, namely that
\begin{align}\label{eq:DeltaB1}
\Delta(B)< c_3 \binom{n-1}{r-1}.
\end{align}
Suppose on the contrary that $d_B(x)\ge c_3 \binom{n-1}{r-1}$ for some $x\in [n]$.
For $j\in[m_i]$, let the $(r-1)$-graph
$B_{x,j}$ consist of those $D\in L_{G}(x)$ such that
if we add $j$ to the profile of $D$ then the obtained $r$-multiset
is bad. In other words, if we move $x$ to $V_j$, then $B_{x,j}$
will be the link of $x$ with respect to the updated bad $r$-graph $B$.
By the definition of $H$, we have
\begin{align}\label{eq:Bxi}
|B_{x,j}|\ge (c_3-c_2)\binom{n-1}{r-1} \quad \mbox{for every $j\in[m_i]$.}
\end{align}

For $\B D=(D_1,\dots,D_{m_i})\in \prod_{j=1}^{m_i} B_{x,j}$,
let $F_{\B D}$ be the $r$-graph that is constructed as follows.
Recall that $F$ is the rigid $P_I$-mixing construction
given by Lemma~\ref{LEMMA:non-extendable-partition} or \ref{lm:tight}, and $\B W$ is its
partition structure.
By relabelling vertices of $F$, we can assume that $x\not\in V(F)$ while
$D\coloneqq \bigcup_{j=1}^{m_i} D_j$ is a subset
of $V(F)$ so that for every $y\in D$ we have
$\branch{F}{y}=\branch{H'}{y}$, that is, $y$ has the
same branches in both $F$ and~$H'$. This is possible
because these $P_I$-mixing constructions have the same trees
of height at most $\ell_0$ while each part of $F$ of height at most $\ell_0$ has at least $m_i(r-1)\ge |D|$
vertices. Finally, add $x$ as a new vertex
and the sets $D_j\cup\{x\}$
for $j\in [m_i]$ as edges, obtaining the $r$-graph~$F_{\B D}$.

\bcl{FBD}{For every $\B D\in
\prod_{j=1}^{m_i} B_{x,j}$ we have $F_{\B D}\in \C F_{M_0}$.}

\bcpf Recall that $M_0$ was chosen to be sufficiently large depending on~$\ell_0$. When we applied Lemma~\ref{LEMMA:non-extendable-partition} or \ref{lm:tight}, the input tree $\B T$ had height $\ell_0$. By~\eqref{eq:M0}, we have $v(F_{\B D})=v(F)+1\le M_0$.

So, suppose on the contrary that we have an embedding $f$ of $F_{\B D}$
into some $P_I$-mixing construction $F'$ with the partition structure $\B U$. By
the rigidity of $F$, we can assume that the base of $F'$ is $P_i$ and that
$f(W_j)\subseteq U_j$ for every $j\in[m_i]$.
Let $j\in[m_i]$ satisfy $f(x)\in U_j$.
Then the edge $D_j\cup\{x\}\in F_{\B D}$ is mapped
into a non-edge because $f(D_j\cup\{x\})$ has bad profile
with respect $U_1,\dots,U_{m_i}$ by the
choice of $D_j\in B_{x,j}$, a contradiction.\ecpf

For every vector $\B D=(D_1,\dots,D_{m_i})\in \prod_{j=1}^{m_j} B_{x,j}$ and every
map $f\colon V(F_{\B D})\to V(G)$ such that $f$ is the identity on
$\{x\}\cup(\bigcup_{j=1}^{m_i} D_j)$ and $f$ preserves branches of height up to $\ell_0$ on all
other vertices, the image $f(F_{\B D})$ has to contain some $X\in
\OO G$ by Claim~\ref{cl:FBD}. (Recall that $G$ is $\C F_{M_0}$-free.) Also,
 \begin{align*}
 f\big(F_{\B D}\setminus\{D_1\cup\{x\},\dots,D_{m_i}\cup\{x\}\}\big)\subseteq H',
 \end{align*}
  that is, the underlying copy  of $F$ on which $F_{\B D}$ was built is embedded by $f$ into $H'$. On the other hand,
 each of the edges $D_1\cup\{x\},\dots,D_{m_i}\cup\{x\}$ of
$F_{\B D}$ that contain $x$
is mapped to an edge of $G$ (to itself). Thus $X\in
H'\setminus G$ and $X\not\ni x$. Any such $X$ can
appear, very roughly, for at most $\binom{w}{r-1}^{m_i}\, (w+1)!\, n^{w-r}$ choices
of $(\B D,f)$, where
$w\coloneqq v(F)=v(F_{\B D})-1$.  On the other hand, the total number of choices of
$(\B D, f)$ is at least $\prod_{j=1}^{m_i} |B_{x,j}| \times (c_4n/2)^{w-(r-1)m_i} \ge
\big((c_3-c_2)\binom{n-1}{r-1}\big)^{m_i}\times (c_4n/2)^{w-(r-1)m_i}$ (since every part of
$H'$ has at least $c_4n$ vertices by~\req{PartSizes}). We conclude that
\begin{align*}
|H\setminus G|\ge |H'\setminus G|
 \ge \frac{\left((c_3-c_2)\binom{n-1}{r-1}\right)^{m_i} \times(c_4n/2)^{w-(r-1)m_i}}{\binom{w}{r-1}^{m_i}\,(w+1)!\, n^{w-r}}>c_2\binom{n}{r}.
\end{align*}
 However, this contradicts~\req{M0End}. Thus~\req{DeltaB1} is proved.

Take any bad edge $D\in B$. We are going to show (in Claim~\ref{cl:M1Deg} below) that
$D$ must intersect $\Omega(c_3n^{r-1})$ absent edges. We need some preparation first.

For each $j\in R_i$ and $y\in D\cap V_j$ pick
some $D_y\in G[V_j]$ such that $D_y\cap D=\{y\}$; it exists
by Part~\ref{it:OmegaNDelta} of
Lemma~\ref{lm:OmegaN}, which gives that
\begin{align}\label{eq:cl:Level1Deg}
  d_{G[V_j]}(y)\ge c_4\binom{n-1}{r-1} \quad \mbox{for all $j\in R_i$ and $y\in
V_j$}.
 \end{align}
Let
$\B D\coloneqq (D,\{D_y\mid y\in D\cap V_{R_i}\})$. (Recall that we denote $V_{R_i}=\bigcup_{j\in R_i} V_j$.)
We define the $r$-graph
$F^{\B D}$
using the rigid $r$-graph $F$ as follows. By relabelling $V(F)$, we can assume
that $X\subseteq V(F)$,
where
 \begin{align}\label{eq:X}
 X\coloneqq \bigcup_{y\in D\cap V_{R_i}} D_y\setminus\{y\},
 \end{align}
so that for every $x\in X$ its branches in $F$ and $H'$
coincide. Again, there is enough space inside $F$ to accommodate all
$|X|\le r(r-1)$ vertices of $X$. Assume also that $D$ is disjoint
from $V(F)$. The vertex
set of $F^{\B D}$ is $V(F)\cup D$. The edge set of  $F^{\B D}$ is defined as follows.
Starting with the edge-set of $F$, add
$D$ and each $D_y$ with $y\in D\cap V_{R_i}$.
Finally, for every $y\in D\cap V_j$ with $j\in [m_i]\setminus R_i$ pick
some $z\in W_j$ and add $\{Z\cup \{y\}\mid Z\in F_z\}$ to the edge set,
obtaining the $r$-graph $F^{\B D}$.
The last step can be viewed as enlarging the part $W_j$ by $D\cap V_j$
and adding those edges that are stipulated by the pattern $P_i$ and
intersect $D$ in at most one vertex.

\bcl{FBD2}{For every $\B D$ as above, we have $F^{\B D}\in\C F_{M_0}$.}

\bcpf By~\eqref{eq:M0},  we have that $v(F^{\B D})=v(F)+r\le M_0$.
Suppose on the contrary that we have an embedding $f$ of $F^{\B D}$
into some $P_I$-mixing construction $F'$ with the partition structure $\B U$. We
can assume by the rigidity of $F$, that the base of $F'$ is $P_i$ and $f(W_j)\subseteq U_j$
for each $i$.

Let us show that for any $y\in D$ we have $f(y)\in U_j$, where the index $j\in [m_i]$ satisfies $y\in W_j$.
First, suppose that $j\in R_i$.
The $(r-1)$-set $f(D_y\setminus\{y\})$ lies entirely
inside $U_j$. We cannot have $f(y)\in U_{k}$ with $k \not= j$ because
otherwise the profile of the edge $f(D_y)$ is $\multiset{\rep{j}{r-1},k}$,
contradicting Lemma~\ref{lm:density=1}.  Thus $f(y)\in U_j$, as claimed.
Next, suppose that $j\in [m_i]\setminus R_i$. Pick some $z\in W_j$.
By the rigidity of $F$, if
we fix the restriction of $f$ to
$V(F)\setminus\{z\}$, then $U_j$ is the only part where $z$ can be mapped to.
By definition,
$y$ and $z$ have the same link $(r-1)$-graphs in  $F^{\B D}$ when restricted to
$V(F)\setminus\{y,z\}$. Hence, $f(y)$ necessarily belongs to $U_j$, as claimed.

Thus the edge $f(D)$ has the same profile as $D\in B$,
a contradiction.\ecpf

\bcl{M1Deg}{For every $D\in B$ there are at least $10rc_3\binom{n-1}{r-1}$
absent edges $Y\in A$ with $|D\cap Y|=1$.}

\bcpf Given $D\in B$, choose the sets $D_y$ for $y\in D\cap V_{R_i}$ as before
Claim~\ref{cl:FBD2}. The condition $D_y\cap D=\{y\}$ rules out
at most $r \binom{n-2}{r-2}$ edges for this $y$. Thus
by~\req{cl:Level1Deg} there are, for example, at least
$(c_4/2)\binom{n-1}{r-1}$ choices of each $D_y$. Form the $r$-graph
$F^{\B D}$ as above and consider potential injective embeddings $f$ of
$F^{\B D}$ into
$G$ that are the identity on $D\cup X$ and map every other
vertex of $F$ into a vertex of $H'$ with the same
branch, where $X$ is defined by \req{X}. For every vertex $x\not\in D\cup X$ we have at least $c_4n/2$
choices for $f(x)$ by \req{PartSizes}. By Claim~\ref{cl:FBD2}, $G$ does not contain $F^{\B D}$
as a subgraph so its image under $f$ contains some $Y\in\OO G$.
Since $f$ maps $D$ and each $D_y$ to an edge of $G$
(to itself) and
 \begin{align*}
 f\left(F^{\B D}\setminus(\{D\}\cup \{D_y\mid y\in D\cap V_{R_i}\})\right)\subseteq H',
 \end{align*}
 we have that $Y\in H'$.
The number of choices of $(\B D,f)$ is at least
 \begin{align*}
 \left((c_4/2)\binom{n-1}{r-1}\right)^{|D\cap V_{R_i}|}\times
(c_4n/2)^{w-(r-1)|D\cap V_{R_i}|}
\ge \left(\frac{c_4n}{4r}\right)^{w},
 \end{align*}
 where $w\coloneqq v(F)$. Assume that for at least half of the time the obtained set
$Y$ intersects $D$ for otherwise we get a contradiction to~\req{M0End}:
 \begin{align*}
 |H'\setminus G|
 \ge \frac12 \times \frac{(c_4n/4r)^{w}}{\binom{w}{r-1}^r\, (w+r)!\, n^{w-r}}>c_2\binom{n}{r}.
 \end{align*}

By the definitions of
$F^{\B D}$ and $f$, we have that $|Y\cap D|=1$ and
$Y\in A$.
Each such $Y\in A$ is counted for at most $\binom{w}{r-1}^r\,(w+r)!\, n^{w-r+1}$
choices of $f$ and $F^{\B D}$. Thus the number of such $Y$ is at least
$\frac1{2}(c_4n/4r)^{w}/(\binom{w}{r-1}^r\,(w+r)!\, n^{w-r+1})$, implying the claim.\ecpf

Let us count the number of pairs $(Y,D)$ where $Y\in A$, $D\in B$,
and $|Y\cap D|=1$. On one hand, each bad edge $D\in B$ creates at least
$10rc_3\binom{n-1}{r-1}$ such pairs by Claim~\ref{cl:M1Deg}.
On the other hand, we trivially
have at most $r|A|\cdot \Delta(B)$ such pairs.
Therefore, $|B|\cdot 10rc_3\binom{n-1}{r-1}\le r|A| \Delta(B)$,
which, by~\req{DeltaB1}, implies that $|A| \ge 10\,|B|$, as desired. This proves Lemma~\ref{lm:Max2}.\epf

Let us state a special case of a result of R\"odl and
Schacht~\cite[Theorem~6]{RS09} that we will need.

\begin{lemma}[Strong Removal Lemma~\cite{RS09}]\label{lm:RS} For every $r$-graph family $\C F$
and $\e>0$ there are $\delta>0$, $M_1$, and $n_0$ such
that the following holds. Let $G$ be an $r$-graph on $n\ge n_0$ vertices
such that for every $F\in \C F$ with $v(F)\le M_1$ the number of $F$-subgraphs in $G$
is at most $\delta n^{v(F)}$. Then $G$ can be made $\C F$-free
by removing at most $\e \binom{n}{r}$ edges.
\end{lemma}

\begin{lemma}\label{lm:edit}
For every $c_0>0$ there is $M_1$ such that
every maximum $\C F_{M_1}$-free $G$ with $n\ge M_1$ vertices is
$c_0\binom{n}{r}$-close
to a $P_{I}$-mixing construction.
\end{lemma}
\bpf Lemma~\ref{lm:RS} for $c_0/2$
gives $M_1$ such that any $\C F_{M_1}$-free $r$-graph $G$
on $n\ge M_1$ vertices can be made into an $\C F_\infty$-free
$r$-graph $G'$ by removing at most $c_0\binom{n}{r}/2$ edges. By
Lemma~\ref{lm:FInfty}, $G'$ embeds into some $P_{I}$-mixing construction $H$
with $v(H)=v(G')$.
Assume that $V(H)=V(G')$ and the identity map is an embedding of
$G'$ into~$H$.

Since $H$ is $\C F_{M_1}$-free,
the maximality of $G$ implies that $|G|\ge |H|$. Thus
$|H\setminus G'|\le c_0\binom{n}{r}/2$ and we can transform $G'$ into
$H$ by changing at most $c_0\binom{n}{r}/2$ further edges.\epf

\subsection{Proof of Theorem~\ref{THM:mixed-pattern}: Putting All Together}
We are ready to prove Part~\ref{it:a} of Theorem~\ref{THM:mixed-pattern}. Let all assumptions of Section~\ref{se:assumptions} apply.

\bpf[Proof of Theorem~\ref{THM:mixed-pattern}\ref{it:a}.]
Let Lemma~\ref{lm:Max2}
return $c_0$ and $M_0$. Then let Lemma~\ref{lm:edit} on input $c_0$ return
some~$M_1$. Finally, take sufficiently large $M$ depending on the previous constants.


Let us argue that this $M$ works in Theorem~\ref{THM:mixed-pattern}\ref{it:a}. As every graph in $\Sigma P_I$ is $\C F_M$-free, it is enough to show that every maximum $\C F_M$-free $r$-graph is a $P_I$-mixing construction. We use induction on the number of vertices~$n$. Let $G$ be any maximum
$\C F_M$-free $r$-graph on $[n]$. Suppose that $n> M$ for otherwise we
are done by Lemma~\ref{lm:FInfty}. Thus Lemma~\ref{lm:edit} applies
and shows that $G$ is $c_0\binom{n}{r}$-close to some
$P_{I}$-mixing construction. 
Lemma~\ref{lm:regular} shows additionally that the minimum degree of $G$ is at least $(\lambda-c_0)\binom{n-1}{r-1}$.
Thus Lemma~\ref{lm:Max2} applies and returns a partition
$[n]=V_1\cup\dots\cup V_{m_i}$ for some $i\in I$ such that \req{Bad-vs-Missing} holds, that is, $|A|\ge 10\, |B|$, where the sets $A$ of absent and $B$ of bad edges are defined by~\eqref{eq:A} and~\eqref{eq:B} respectively. Now, if we take the union of $\blow{E_i}{V_1,\dots,V_{m_i}}$ with $\bigcup_{j\in R_i} G[V_j]$, then the
obtained $r$-graph is still $\C F_{M}$-free by Lemma~\ref{lm:FmFree} and has
exactly $|A|-|B|+|G|$ edges.
The maximality of $G$ implies that $|B| \ge |A|$. By~\req{Bad-vs-Missing}, this is possible only if $A=B =\emptyset$. Thus, $G$ coincides with the blowup $\blow{E_i}{V_1,\dots,V_{m_i}}$, apart edges inside the recursive parts $V_j$, $j\in R_i$. 

Let $j\in R_i$ be arbitrary. By Lemma~\ref{lm:FmFree} if we replace
$G[V_j]$ by any $\C F_M$-free $r$-graph, then the new
$r$-graph on $V$ is still $\C F_M$-free.
By the maximality of $G$, we conclude
that $G[V_j]$ is a maximum $\C F_M$-free $r$-graph. By the induction
hypothesis (note that $|V_j|\le n-1$), the induced subgraph $G[V_j]$ is a
$P_{I}$-mixing construction.
It follows that $G$ is a $P_{I}$-mixing construction itself, which implies Theorem~\ref{THM:mixed-pattern}\ref{it:a}.
\epf

In order to prove Part~\ref{it:b} of Theorem~\ref{THM:mixed-pattern}, 
we prove the following partial result first, where stability is proved for the ``bottom'' edges only.

\begin{lemma}\label{LEMMA:stability-level-1}
There exists $M_2 \in \mathbb{N}$ so that
for every $\varepsilon>0$ there exist $\delta_0>0$ and $n_0$ such that the following holds for all $n\ge n_0$.
Suppose that $G$ is a $\mathcal{F}_{M_2}$-free $r$-graph on $n$ vertices with $|G| \ge (1-\delta_0)\mathrm{ex}(n,\mathcal{F}_{M_2})$.
Then there exist $i\in I$ and a partition $V(G) = V_1 \cup \cdots \cup V_{m_i}$ with $|V_{j}| \in [\beta n/4, (1-\lambda/4r)n]$ for all $j\in [m_i]$ 
such that $G'\coloneqq  G\setminus \left(\bigcup_{j\in R_i}G[V_j]\right)$ satisfies $|G'\triangle \blow{E_i}{V_1,\dots,V_{m_i}}| \le \varepsilon n^r$.
\end{lemma}
\bpf
Let $c_0$ and $M_0$ be the constants returned by Lemma~\ref{lm:Max2}. Then let $M_2$ be sufficiently large, in particular so that it satisfies Lemma~\ref{lm:RS} for $c_0/4$. Let us show that this $M_2$ works in the lemma. Given any $\varepsilon>0$, choose sufficiently small positive constants $\delta_1\gg \delta_0$.
Let $G$ be an $\mathcal{F}_{M_2}$-free $r$-graph on $n\to\infty$ vertices with $|G| \ge (1-\delta_0)\mathrm{ex}(n,\mathcal{F}_{M_2})$.
By our choice of $M_2$, the $r$-graph $G$ can be embedded into some $n$-vertex $P_{I}$-mixing construction $H$
by removing at most $c_0\binom{n}{r}/4$ edges. 
Since $|G| \ge (1-\delta_0)\mathrm{ex}(n,\mathcal{F}_{M_2}) \ge |H| - \delta_0 \binom{n}{r}$,
we have $|G\triangle H| \le 2 \cdot c_0\binom{n}{r}/4 + \delta_0 \binom{n}{r} \le 3c_0 \binom{n}{r}/4$.

Define
\begin{align*}
Z\coloneqq  \left\{v\in V(G) \colon d_G(v)\le  (\lambda-r\delta_1)\binom{n-1}{r-1}\right\}.
\end{align*}

\bcl{Size-ZG}{We have that $|Z|< \delta_1 n$.}

\bcpf
Suppose to the contrary that $|Z| \ge \delta_1 n$. Let $G'$ be obtained from $G$ by removing some $\delta_1 n$ vertices of $Z$. 
We have that
\begin{align*}
|G'|-\lambda \binom{n-\delta_1 n}{r}
& \ge |G| - \delta_1 n(\lambda-r \delta_1)\binom{n-1}{r-1}- \lambda (1-\delta_1)^r \binom{n}{r} +o(n^r)\\
& \ge \left((\lambda-\delta_0) - r\delta_1(\lambda- r \delta_1)-\lambda\left(1-r\delta_1+ \binom{r}{2}\delta_1^2\right)\right)\binom{n}{r} +o(n^r)\\
& \ge \left(- \delta_0  +r^2\delta_1^2 -\lambda \binom{r}{2}\delta_1^2 \right) \binom{n}{r}+o(n^r)\ > \Omega(n^r).
\end{align*}
 Thus the $\C F_{M_2}$-free $r$-graph $G'$ contradicts the consequence of Theorem~\ref{THM:mixed-pattern}\ref{it:a} 
 that $\pi(\C F_{M_2})=\lambda$.
\ecpf

Let  $n_1 \coloneqq  n-|Z|$ and $G_1 \coloneqq  G-Z$.
So $G_1$ is an $r$-graph on $n_1 \ge (1-\delta_1)n$ vertices with
\begin{align*}
\delta(G_1) & \ge (\lambda-r \delta_1)\binom{n-1}{r-1} - \delta_1 n\binom{n-2}{r-2} \ge (\lambda-2r\delta_1)\binom{n-1}{r-1}, \quad \text{and}\\
|G_1| & \ge (1-\delta_0)\mathrm{ex}(n,\mathcal{F}_{M_2}) - \delta_1 n\binom{n-1}{r-1} \ge \mathrm{ex}(n,\mathcal{F}_{M_2}) -2r\delta_1 \binom{n}{r}.
\end{align*}
Let $H_1$ be the induced subgraph of $H$ on $V(G)\setminus Z$ and note that $H_1$ is also a $P_I$-mixing construction.
Since $|G_1\triangle H_1| \le |G\triangle H| \le 3c_0\binom{n}{r}/4 \le c_0 \binom{n_1}{r}$,
by Lemma~\ref{lm:Max2} there are $i\in I'$ and a partition
$V(G)\setminus Z=V_1'\cup\dots\cup V_{m_i}'$ such that $|V_j'| \in [\beta n_1/2, (1-\lambda/2r)n_1]$ for all $j\in [m_i]$ and
it holds that $|A_1|\ge 10\, |B_1|$, where $A_1\coloneqq \blow{E_i}{V_1',\dots,V_{m_i}'}\setminus G$ and $B_1 \coloneqq  G_1\setminus \big(\blow{E_i}{V_1',\dots,V_{m_i}'}\cup \bigcup_{j\in R_i} G_1[V_j]\big)$.
If we take the union of $\blow{E_i}{V_1',\dots,V_{m_i}'}$ with $\bigcup_{j\in R_i} G_1[V_j']$, then the
obtained $r$-graph is still $\C F_{M_2}$-free (by Lemma~\ref{lm:FmFree}) and has
exactly $|G_1|+|A_1|-|B_1| \ge |G_1| + \frac{9}{10}\,|A_1|$ edges.
Therefore, $|G_1| + \frac{9}{10}\,|A_1| \le \mathrm{ex}(n_1,\mathcal{F}_{M_2})$, which implies that
\begin{align*}
|A_1| + |B_1|
\le \frac{11}{10} \cdot \frac{10}{9}\left(\mathrm{ex}(n_1,\mathcal{F}_{M_2}) - |G_1|\right)
\le 4r\delta_1 \binom{n}{r}.
\end{align*}
Now, extend the partition $V_1',\dots,V_{m_i}'$ arbitrarily to $V(G)$, for example, by defining
\begin{align*}
V_j \coloneqq  
\begin{cases}
V_1'\cup Z, & \quad{j = 1}, \\
V_j',         & \quad{j\neq 1}. 
\end{cases}
\end{align*}
Then simple calculations show that $|V_j| \in [\beta n/4, (1-\lambda/4r)n_1]$ for all $j\in [m_i]$, and
the $r$-graph $G' \coloneqq  G \setminus \left(\bigcup_{j\in R_i}G[V_j]\right)$ satisfies  
\begin{align*}
|G'\triangle \blow{E_i}{V_1,\dots,V_{m_i}}| 
\le |Z| \binom{n-1}{r-1} + |A_1| + |B_1|
\le \varepsilon n^r. 
\end{align*}
This completes the proof of Lemma~\ref{LEMMA:stability-level-1}. 
\epf

Now we are ready to prove Part~\ref{it:b} of Theorem~\ref{THM:mixed-pattern}.

\bpf[Proof of Theorem~\ref{THM:mixed-pattern}\ref{it:b}.]
Let $M_0$ and $c_0$ be returned by Lemma~\ref{lm:Max2}.
Let $M_1$ be returned by Lemma~\ref{lm:edit} for $c_0$, and let $M_2$ be returned by Lemma~\ref{LEMMA:stability-level-1}. 
Let us show that $M\coloneqq \max\{ M_0,M_1,M_2\}$ works in Theorem~\ref{THM:mixed-pattern}\ref{it:b}.

Take any $\varepsilon >0$.  
Let $\ell \in \mathbb{N}$ be a sufficiently large integer such that, in particular, $\left(1-\frac{\lambda}{4r}\right)^{\ell} \ll \varepsilon$. 
Next, choose sufficiently small positive constants 
$\delta_{\ell}\gg \dots\gg \delta_1\gg \delta$.
Let $n$ be sufficiently large. 
Let $G$ be an $\mathcal{F}_{M}$-free $r$-graph on $n$ vertices with $|G|\ge (1-\delta)\mathrm{ex}(n, \mathcal{F}_M)$. 
By Lemma~\ref{LEMMA:stability-level-1}, there exist $i\in I$ and a partition $V(G) = V_1 \cup \cdots \cup V_{m_i}$ with 
$|V_j| \in [\beta n/4, (1-\lambda/4r)n]$ such that
$G'\coloneqq  G\setminus \left(\bigcup_{j\in R_i}G[V_j]\right)$ satisfies $|G'\triangle \blow{E_i}{V_1,\dots,V_{m_i}}| \le \delta_1 n^r$.

Let $\hat{G} \coloneqq  \blow{E_i}{V_1,\dots,V_{m_i}} \cup \left(\bigcup_{j\in R_i}G[V_j]\right)$. 
Then Lemma~\ref{lm:FmFree} implies that $\hat{G}$ is still $\mathcal{F}_{M}$-free,
and our argument above shows that $|\hat{G}| \ge |G| - \delta_1 n^r \ge \mathrm{ex}(n,\mathcal{F}_M) - 2\delta_1 n^r$.

Now take any $j\in R_i$. 
Note that we have $|G[V_j]| \ge \mathrm{ex}(|V_j|, \mathcal{F}_{M}) - 2\delta_1 n^r$, 
since otherwise, by replacing $G[V_j]$ in $\hat{G}$ by a maximum $\mathcal{F}_{M}$-free $r$-graph on $V_j$ we would get an $r$-graph which is $\mathcal{F}_{M}$-free (by Lemma~\ref{lm:FmFree}) and has more than $\mathrm{ex}(n,\mathcal{F}_M)$ edges, a contradiction.
Since $|V_{j}| \ge \beta n/4$ is sufficiently large and $\delta_1 \ll \delta_2$, there exist, by Lemma~\ref{LEMMA:stability-level-1} again, an index $i' \in I$ and a partition $V_j = V_{j,1}' \cup \cdots \cup V_{j,m_{i'}}'$ such that
$|V_{j,k}'| \in [(\beta/4)^2 n, (1-\lambda/4r)^2n]$ for all $k\in [m_{i'}]$ and
$G_j'\coloneqq  G[V_j]\setminus \left(\bigcup_{k\in R_{i'}}G[V_{j,k}]\right)$ satisfies $|G_j'\triangle \blow{E_{i'}}{V_{j,1},\dots,V_{j,m_{i'}}}| \le \delta_2 |V_j|^r$. 
Summing over all $j\in R_i$, we get 
\begin{align*}
\sum_{j\in R_i}|G_j'\triangle \blow{E_{i'}}{V_{j,1},\dots,V_{j,m_{i'}}}| 
\le \delta_2 \sum_{j\in R_i}|V_j|^r
\le \delta_2 \left(\sum_{j\in R_i}|V_j|\right)^r
\le \delta_2 n^r. 
\end{align*}
Repeating this argument until the $\ell$-th level, we see that $G$ can be made into a $P_I$-mixing construction by removing and adding at most 
\begin{align*}
\sum_{i=1}^{\ell}\delta_i n^r + \sum_{\B i}\binom{|V_{\B i}|}{r}
& \le \sum_{i=1}^{\ell}\delta_i n^r + \frac{n}{\left(1-\frac{\lambda}{4r}\right)^{\ell}n}\binom{\left(1-\frac{\lambda}{4r}\right)^{\ell} n}{r} \\
& \le \left(\sum_{i=1}^{\ell}\delta_i + \left(1-\frac{\lambda}{4r}\right)^{\ell}\right)n^r
\le \varepsilon n^r
\end{align*}
edges.
Here the second summation is over all legal (with respect to $G$) vectors of length~$\ell$. 
This completes the proof of Theorem~\ref{THM:mixed-pattern}\ref{it:b}.
\epf

\section{Finite $r$-graph families with rich extremal Tur\'an constructions}\label{se:rich}

In this section, we prove Theorem~\ref{th:ExpMany}. We need some preliminaries first.

Recall that a \textbf{Steiner triple system} on a set $V$ is a $3$-graph $D$ with vertex set $V$ such that every pair of distinct elements of $V$ is covered by exactly one edge of $D$. Let $\STS_t$ be the set of all Steiner triple systems on~$[t]$. It is known that such a design $D$ exists (and thus $\STS_t\not=\emptyset$)  if and only if the residue of $t\ge 3$ modulo 6 is $1$ or $3$, a result that was proved by  Kirkman~\cite{Kirkman47} already in 1847.

We will need the following result, which is a special case of~\cite[Proposition~2.2]{LMR1}.

\begin{lemma}[\cite{LMR1}]
\label{lm:Prop2.2} Let $t\ge 55$ and $D\in \STS_t$ be arbitrary.
Then, for every $(x_1,\dots,x_t)\in \IS_t$, it holds that
 \begin{align*}
 \lambda_{\OO D}(x_1,\ldots,x_t)\le \lambda_{\OO D}\left(\frac1t,\dots,\frac1t\right)- \frac23 \sum_{i=1}^t \left(x_i-\frac1t\right)^2.
 \end{align*}
 (Recall that $\OO D=\binom{[t]}{3}\setminus D$ denotes the complement of $D$.)
\end{lemma}

The above lemma implies, in particular, that the uniform weight is the unique $\OO D$-optimal vector. Also,
note that
$
|\OO D|=\binom{t}{3}- \frac13 \binom{t}{2}=\frac{t(t-1)(t-3)}6
$
and thus
\begin{align*}
\lambda_{\OO D}=\lambda_{\OO D}\left(\frac1t,\dots,\frac1t\right)=\frac{3!\, |\OO D|}{t^3}=\frac{(t-1)(t-3)}{t^2}.
\end{align*}
Here, $\lambda_{\OO D}:=\lambda_{(t,\OO D,\emptyset)}$ denotes the maximum of the Lagrange polynomial $\lambda_{\OO D}(\B x)$ 
of the 3-graph $\OO D$ over $\B x\in \IS_t$. 

For $D\in\STS_t$, let $P_D\coloneqq (t,\OO D,[t])$ be the pattern where every vertex of the complementary $3$-graph $\OO D$ is made recursive. If we take the uniform blow-ups of $\OO D$ at all levels when making a $P_D$-construction, then the obtained limiting edge density $\lambda_1$ satisfies the relation $\lambda_1=\lambda_{\OO D}+\lambda_1 t (1/t)^3$ which gives that
\begin{align*}
\lambda_{P_D}\ge \lambda_1= \frac{\lambda_{\OO D}}{1-1/t^2}=\frac{(t-1)(t-3)}{t^2-1}=\frac{t-3}{t+1}.
\end{align*}

\begin{lemma}
\label{lm:OptHRec} Let $t$ be sufficiently large and let $D\in STS_t$. Let $f(\B x)\coloneqq \lambda_{\OO D}(\B x)+\lambda_1 \sum_{i=1}^t x_i^3$.
Then $f(\B x)\le \lambda_1$ for every $\B x\in\IS_t^*$ with equality if and only if $\B x$ is the uniform vector $(\frac1t,\ldots,\frac1t)$.

It follows that $\lambda_{P_D}=\frac{t-3}{t+1}$ and  the set of $P_D$-optimal vectors consists only of the uniform vector $(\frac1t,\ldots,\frac1t)\in\IS_t^*$. \end{lemma}

\bpf Take any $\B x\in\IS_t^*$. In order to prove that $f(\B x)\le \lambda_1$, we split the argument into two cases.

First, suppose that $\max\{x_i\mid i\in [t]\}\le 1/2$.
Here, we have by Lemma~\ref{lm:Prop2.2} that
\begin{align*}
f(\B x)\le \lambda_{\OO D}\left(\frac1t,\dots,\frac1t\right) -\frac23\sum_{i=1}^t \left(x_i-\frac1t\right)^2+\lambda_1 \sum_{i=1}^t x_i^3 =\lambda_{\OO D}+\sum_{i=1}^t g(x_i),
\end{align*}
where $g(x)\coloneqq -\frac23 (x-\frac1t)^2+\lambda_1 x^3$ for $x\in [0,\frac12]$.

 The second derivative of the cubic polynomial $g$ has zero at $x_0\coloneqq \frac{2(t+1)}{9(t-3)}$. We have that  $x_0>1/t$. While this can be checked to hold for every $t$, it is trivial  here since $t$ was assumed to be sufficiently large. Thus function $g$ is concave on $[0,\frac1t]$. Unfortunately, it is not concave on the whole interval $[0,1/2]$ so we consider a different function $g^*$ which equals $g$ on $[0,\frac1t]$ and whose graph on $[\frac1t,\frac12]$ is the line tangent to the graph of $g$ at $1/t$, that is, we set
 \begin{align*}
  g^*(x)\coloneqq g\left(\frac1t\right)+g'\left(\frac1t\right) \left(x-\frac1t\right)
  ,\quad\mbox{for $x\in [1/t,1/2]$}.
     \end{align*}
  By above,  $g^*$ is a concave function on $[0,1/2]$. Also, we have that $g\le g^*$ on this interval. Indeed, since the second derivative $g''$ changes sign from negative to positive at $x_0>1/t$, it is enough to check only that $g(1/2)\le g^*(1/2)$ and routine calculations give that
  \begin{align*}
  g^*(1/2)-g(1/2)=\frac{(t-2)^2(t^2+t+36)}{24t^3(t+1)}>0,
  \end{align*}
  as desired. Thus, since $g^*$ is a concave majorant of $g$, we have that
  \begin{align*}
   \frac1t \sum_{i=1}^t g(x_i)\le \frac1t\sum_{i=1}^t g^*(x_i)\le g^*\left(\frac{x_1+\cdots+x_t}t\right) = g^*(1/t)=g(1/t).
   \end{align*}
 This gives that $f(\B x)\le \lambda_{\OO D}+tg(1/t)=\lambda_1$. Moreover, if we have equality then $x_1=\cdots=x_t=1/t$ (because $g^*$ is strictly concave on $[0,1/t]$).

Thus it remains to consider the case when some $x_i$ is strictly larger than $1/2$. Without loss of generality, assume that $x_1> 1/2$. Here, we can bound
\begin{align*}
 f(\B x)\le h(x_1)\coloneqq 3!\, x_1(1-x_1)^2\frac{t-3}{2(t-1)} 
 +\frac{3!\,\binom{t-1}{3}}{(t-1)^3}\, (1-x_1)^3
 + \lambda_1 \left(x_1^3+(1-x_1)^3\right),
 \end{align*}
 where the stated three terms come from the following arguments. The first term accounts for the triples containing $x_1$ in the Lagrange polynomial~$\lambda_{\OO D}(\B x)$. The link graph $L_D(1)$ is just a perfect matching $M$ on $\{2,\dots,t\}$ (because $D$ is a Steiner triple system) and receives total weight $1-x_1$. As it is well-known (see e.g.\ \cite{MS65}),
 the Largrangian of a graph is maximised by putting the uniform weight on a maximum clique which, for the complement $L_{\OO D}(1)$ of a perfect matching, has size $s\coloneqq (t-1)/2$. Thus $(1-x_1)^{-2}\sum_{ij\in L_D(1)} x_ix_j\le  \binom{s}{2}/s^2=\frac{s-1}{2s}=\frac{t-3}{2(t-1)}$, giving the first term. The second term just upper bounds the Lagrangian of ${\OO D}-1={\OO D}[\{2,\ldots,t\}]$ by the Largrangian of the complete $3$-graph on $t-1$ vertices, scaling the result by the cube of the total weight $1-x_1$. The third term uses the fact that the sum of cubes of non-negative entries with sum $1-x_1$ is maximised when we put all weight on a single element.

The coefficient at $x^3$ in the cubic polynomial $h(x)$ is $\frac{2t^2-10t+9}{(t-1)^2}>0$. Also, the derivative of $h$ has two roots, which are $\pm1/\sqrt2+o(1)$ as $t\to\infty$. Thus (since $t$ is large),
the function $h$, when restricted to the interval  on $[1/2,1]$, first decreases and then increases. So,
in order to show that $f(\B x)\le \lambda_1$,
it is enough to check the inequality $h(x)\le \lambda_1$ only at the points $x=1/2$ and $x=1$. There, the values of $h(x)-\lambda_1$ are respectively $-\frac{2 t^3-20 t^2+47 t-27}{8 (t-1)^2
   (t+1)}<0$
and~$0$. Thus $f(\B x_1)\le \lambda_1$ also in Case~2. Furthermore, the equality can only be attained if $x_1=1$ (and $x_2=\cdots=x_t=0$); however, we exclude standard basis vectors from $\IS_t^*$. This proves the first part of the lemma.

Using the proved inequality $f\le \lambda_1$, one can show by a simple induction on the number of levels that every $P$-construction on $n\to\infty$ vertices has edge density at most $\lambda_1+o(1)$. Thus $\lambda_{P_D}=\lambda_1$. Also, the set of $P_D$-extremal vectors, which by definition consists of $\B x\in\IS_t^*$ with $f(\B x)=\lambda_1$, is exactly as claimed.\epf

In the rest of the this section, whenever we have any $I\subseteq \STS_t$, a set of Steiner triple systems on~$[t]$, we denote $P_I\coloneqq \{P_D\mid D\in I\}$. (Recall that $P_D=(t,\OO D,[t])$.) Also, let us call a partition $V=V_1\cup\cdots\cup V_t$ \textbf{balanced} if for all $i,j\in [t]$ we have $\left|\,|V_i|-|V_j|\,\right|\le 1$.

\begin{lemma}
\label{lm:OptHRecExact} Let $t$ be sufficiently large and take any non-empty $I\subseteq \STS_t$. Then there is $n_0=n_0(I)$ such that  every maximum $P_I$-mixing construction $G$ on $n\ge n_0$ vertices has its bottom partition balanced.
\end{lemma}

\bpf View $n$ as tending to $\infty$ and take  any maximum $P_I$-mixing construction $G$ on~$[n]$. Let $G$ have the base $P_D$ and the bottom partition $[n]=V_1\cup\cdots\cup V_t$. Let $v_i\coloneqq |V_i|$ and $e_i\coloneqq |G[V_i]|$ for $i\in [t]$. By Lemma~\ref{lm:OptHRec}, we have $v_i=(1/t+o(1))n$ for each $i\in [t]$.

Suppose on the contrary that some two part sizes differ by more than~1, say, $v_1\ge v_2+2$. Let $u$ (resp.\ $w$) be a vertex of minimum degree in $G[V_1]$ (resp.\ maximum degree in $G[V_2]$). Let $G'$ be obtained from $G$ by removing $u$ and adding a double $w'$ of $w$. Thus $G'$ is also a $P_I$-mixing construction, with the bottom parts $V_1'\coloneqq V_1\setminus\{u\}$, $V_2'\coloneqq V_2\cup \{w'\}$, and $V_i'\coloneqq V_i$ for $3\le i\le t$. By the maximality of $G$, we have
$
|G|\ge |G'|.
$

Let us estimate $|G'|-|G|$. The contribution from the edges that intersect both the first part and the second part is
\begin{align*}
\left((v_1-1)(v_2+1)-v_1v_2\right) \sum_{i\colon \{1,2,i\}\in {\OO D}} |V_i| = (v_1-v_2-1) \left(\frac{t-3}{t}+o(1)\right)n,
\end{align*}
 where the equality uses Lemma~\ref{lm:OptHRec} and the fact that there are exactly $t-3$ triples containing the pair~$\{1,2\}$. By the maximality of $G$, each part $V_i$ spans a maximum $P_I$-mixing construction; thus $e_i=\Lambda_{P_I}(v_i)$. Since $v_1\ge v_2$, we have by Lemma~\ref{lm:monotone} that $\Lambda_{P_I}(v_1)/\binom{v_1}{3}\le \Lambda_{P_I}(v_2)/\binom{v_2}{3}$. Thus the degree of $u$ in $G[V_1]$ is at most
\begin{align*}
\frac{3e_1}{v_1}=\frac{3\Lambda_{P_I}(v_1)}{v_1}\le
\frac{3\Lambda_{P_I}(v_2)\, \binom{v_1}{3}}{v_1\binom{v_2}{3}}
= \frac{3e_2}{v_2}\cdot \frac{(v_1-1)(v_1-2)}{(v_2-1)(v_2-2)}.
\end{align*}
 The degree of $w'$ in $G'$ is at least the degree of $w$ in $G$ which in turn is at least $3e_2/v_2$. Thus the contribution to $|G'|-|G|$ of edges inside the first or second part is at least
\begin{align*}
\frac{3e_2}{v_2}-\frac{3e_1}{v_1} \ge \frac{3e_2}{v_2}\left(1-\frac{(v_1-1)(v_1-2)}{(v_2-1)(v_2-2)}\right) = -\frac{3e_2}{v_2}\cdot \frac{(v_1-v_2)(v_1+v_2-3)}{(v_2-1)(v_2-2)}.
\end{align*}
 This is, by $e_2=(\lambda_{P_D}+o(1))\binom{v_2}{3}=\left(\frac{t-3}{t+1}+o(1)\right)\binom{n/t}{3}$,
 \begin{align*}
(-1+o(1))\,\frac{3\,\frac{t-3}{t+1}\,\binom{n/t}{3}}{n/t}\cdot \frac{(v_1-v_2)\,2n/t}{(n/t)^2}=-\frac{(v_1-v_2)(t-3)n}{t(t+1)} +o((v_1-v_2)n).
\end{align*}
 By putting both contributions together and using $v_1-v_2-1\ge (v_1-v_2)/2$, we obtain that
\begin{align*}
0&\ge |G'|-|G|\ \ge\ (v_1-v_2-1)\, \frac{t-3}{t}n-\frac{(v_1-v_2)(t-3)n}{t(t+1)}+
o((v_1-v_2)n)\\
&\ge (v_1-v_2)\,\frac{(t-3)n}{t}\left(\frac12 - \frac1{t+1}\right) +
o((v_1-v_2)n)\ >\ 0,
\end{align*}
 which is the desired contradiction.\epf

For every $D\in\STS_t$, we construct a parameter $F(D)$ of much lower complexity than $D$ that nonetheless contains enough information  to compute the sizes of maximum balanced $\OO D$-blowups of all large orders. More precisely, we
do the following for every $q\in\{0,\ldots,t-1\}$. Let $\ell\to\infty$. For every $q$-set $Q\subseteq [t]$, consider a $\OO D$-blowup $G$ on $[t\ell+q]$ with partition $V_1\cup\cdots\cup V_t$, where
\begin{align*}
|V_i|=\left\{ \begin{array}{ll}
 \ell+1,&\mbox{if $i\in Q$,}\\
 \ell,&\mbox{if $i\in [t]\setminus Q$}.
\end{array}\right.
\end{align*}
 Thus, the $q$ larger parts are exactly those specified by $Q$. Clearly, the size of $G$ is
\begin{align*}
 p_{D,Q}(\ell)\coloneqq (\ell+1)^3 t_{D,Q,3} + \ell(\ell+1)^2 t_{D,Q,2}+ \ell^2(\ell+1) t_{D,Q,1}+\ell^3\left(|\OO D|-t_{D,Q,3}-t_{D,Q,2}-t_{D,Q,1}\right),
\end{align*}
 where, for $i\in [3]$, we let $t_{D,Q,i}$ denote the number of triples in $\OO D$ that have exactly $i$ vertices in $Q$. This is a polynomial function of~$\ell$. If we take another $q$-set $Q'$ then the polynomial $p_{D,Q}(\ell)-p_{D,Q'}(\ell)$ does not change sign for large $\ell$ (possibly being the zero polynomial). Since there are finitely many choices of $Q$ (namely, $\binom{t}{q}$ choices), there are $Q_{D,q}\in \binom{[t]}{q}$ and $\ell_0$
such that
 \begin{align}\label{eq:QD}
 p_{D,Q_{D,q}}(\ell)\ge p_{D,Q'}(\ell),\quad\mbox{for all $Q'\in \binom{[t]}{q}$ and  $\ell\ge \ell_0$}.
\end{align}
 Fix one such $Q_{D,q}$ for each $q\in \{0,\dots,t-1\}$ and define \begin{align}\label{eq:F(D)}
 F(D)\coloneqq ((t_{D,Q_{D,q},i})_{i=1}^3)_{q=0}^{t-1}.
\end{align}
The number of possible values of $F$ is upper bounded by, say, $t^{9t}$ because each individual $t_{D,Q,i}$ assumes at most $\binom{t}{3}+1\le t^3$ possible values.

\hide{
We claim that the sequence $F(D)\coloneqq ((t_{D,Q_{D,q},i})_{i=1}^3)_{q=0}^{m-1}$ determines the maximum size of a balanced $D$-blowup $G$ of large order~$n$. Indeed, let $n$ to sufficiently large so that, in particular,~\eqref{eq:QD} holds
 with respect to  $\ell\coloneqq \lfloor n/t\rfloor$, $q\coloneqq n-t\ell$, and $Q\coloneqq Q_{D,q}$.
Let $V_1\cup\cdots\cup\ldots V_t$ be a balanced partition that maximises $|\blow{D}{V_1,\ldots,V_t}|$. Define $Q'\coloneqq \{i\in [t]\colon |V_i|=\ell+1\}$. Thus $Q'$ is the index set of the larger parts and, of course, its size is $q$, the residue of $n$ modulo~$t$. By our choice of $Q_{D,q}$ and since $n$ is large, the inequality in~\eqref{eq:QD} holds. It must be in fact equality by the maximality of $\blow{D}{V_1,\ldots,V_t}$, proving the claim.
}

\begin{lemma}
\label{lm:F(D)} Let $t$ be sufficiently large. Let $I\subseteq \STS_t$ be any subset on which the above function $F$ is constant. Then there is $n_0$ such that for all $n\ge n_0$ and all $D\in I$ there is a maximum  $P_I$-mixing construction $G$ on $[n]$ with the base~$P_D$.
\end{lemma}

\bpf Choose $n_0$ sufficiently large, in particular so that, for every $n\ge n_0$, \eqref{eq:QD} holds for every $D\in I$ with respect to $\ell\coloneqq \lfloor n/t\rfloor$ and $q\coloneqq n-t\ell$.

Take any $n\ge n_0$ and $D\in I$. We have to exhibit a maximum $P_I$-mixing construction $G$ with base~$P_D$. Let $\ell\coloneqq \lfloor n/t\rfloor$, $q\coloneqq n-t\ell$, and $Q\coloneqq Q_{D,q}$.
Take a balanced partition $V_1\cup\cdots\cup V_t$ of $[n]$ with $|V_j|=\ell+1$ exactly for $j\in Q$. Let $G$ be obtained from $\blow{{\OO D}}{V_1,\ldots,V_t}$ by adding for every $j\in [t]$ a maximum $P_I$-mixing construction on~$V_j$.

Let us show that the $P_I$-mixing construction $G$ is maximum.
Note that the size of the graph $G$ is
\begin{align}\label{eq:F(D)Size}
 |G|=p_{D,Q}(\ell) + q\, \Lambda_{P_I}(\ell+1)+(t-q)\,\Lambda_{P_I}(\ell).
\end{align}
Let $G'$ be any maximum $P_I$-mixing construction on $[n]$. Let $G'$ have the base $D'$ and the bottom partition $V_1'\cup\dots\cup V_t'$. This partition must be balanced by Lemma~\ref{lm:OptHRecExact}, since $n$ is large. Let $Q'\in \binom{[t]}{q}$ consist of the indices of parts of size $\ell+1$. Clearly, $|Q'|=q$. By the maximality of $G'$, every part $V_s'$, $s\in [t]$, induces a maximum $P_I$-mixing construction.
Thus the obvious analogue of~\eqref{eq:F(D)Size} holds for $G'$ as well. Let $Q''\coloneqq Q_{D',q}$. Since $\ell\ge \lfloor n_0/t\rfloor$ is large, we have
$p_{D',Q'}(\ell)\le p_{D',Q''}(\ell)$. Since $F(D)=F(D')$, the polynomials $p_{D',Q''}$ and $p_{D,Q}$ are the same. Putting all together, we obtain
\begin{align*}
 |G'|-|G|=p_{D',Q'}(\ell)-p_{D,Q}(\ell)\le  p_{D',Q''}(\ell)-p_{D,Q}(\ell)=0.
\end{align*}
 Thus $|G|$ is indeed a maximum $P_I$-mixing construction.
\epf

Let us remark that Lemma~\ref{lm:F(D)} need not hold for small $n$ when it is in general possible that some of the patterns $P_D$ for $D\in I$ cannot be the base in a maximum $P_I$-mixing construction.

\bpf[Proof of Theorem~\ref{th:ExpMany}.]
Keevash~\cite[Theorem 2.2]{Keevash18} proved that if $t\to\infty$ is 1 or 3 modulo 6 then the number of Steiner triples systems on $[t]$ is $(t/\me^2+o(1))^{t^2/6}$. Note that the function $F$ assumes at most $t^{9t}$ values while each isomorphism class of $\STS_t$ has at most $t!$ elements. Thus we can fix a sufficiently large $t$ and a subset $I\subseteq \STS_t$ consisting of non-isomorphic 3-graphs such that $F$ is constant on $I$ while
\begin{align*}
|I|\ge \frac{(t/\me^2+o(1))^{t^2/6}}{t!\, t^{9t}}> t!.
\end{align*}
 Let $\C F$ be the finite family of $3$-graphs returned by Theorem~\ref{THM:mixed-pattern}; thus maximum $\C F$-free $3$-graphs are exactly maximum $P_I$-mixing constructions. Let us show that this family $\C F$ satisfies both parts of Theorem~\ref{th:ExpMany}. Let $n_0$ be sufficiently large.

Given any $n\in \mathbb{N}$, let $V\coloneqq [n]$ and consider the family $\C P$ of 3-graphs $G$ that can be recursively constructed as follows. If the current vertex set $V$ has less than $n_0$ vertices, put any maximum $P_I$-mixing construction on $V$ and stop. So suppose that $|V|\ge n_0$. Pick any $D\in I$. Let $G'$ be a maximum $P_I$-mixing construction on $V$ with the base $P_D$, which exists by Lemma~\ref{lm:F(D)}. Let $V=V_1\cup\cdots\cup V_t$ be the bottom partition of~$G'$. Let $G$ be obtained by taking all bottom edges of $G'$ and adding for each $i\in [t]$ a $3$-graph on $V_i$ that can be recursively constructed by the above procedure.

Let us observe some easy properties of any obtained $G\in\C P$. By definition, $G$ is a $P_I$-mixing construction. In fact, it is a maximum one, which can be shown by induction on the number of vertices: the initial $P_I$-mixing construction $G'$ is maximum and when we ``erase'' edges in $G'[V_i]$, we add back the same number of edges by induction.

\hide{
As the base, pick any $P_D$ for $D\in I$. Let $\ell\coloneqq \lfloor n/t\rfloor$, $q\coloneqq n-t\ell$ and $Q\coloneqq Q_{D,q}$. Take a balanced partition $V\coloneqq V_1\cup\cdots\cup V_t$ so that the parts of size $\ell+1$ are exactly $V_i$ with $i\in Q$. Take the blow up $\blow{\OO D}{V_1,\ldots,V_t}$ and apply recursion inside each set~$V_i$. (Thus, for example, if each of the $t$ level-$1$ sets has at least $n_0$ elements, then we have $|I|^t$ choices of bases at this level.)
}

Let the \textbf{size-$m$ truncated tree} $\B T_{G}\SizeTruncated{m}$ of $P_I$-mixing construction $G$ be the subtree of $\B T_{G}$ where we keep only those nodes that corresponds to parts in $\BV{G}$ of size at least~$m$. Of course, if a node is not included into $\B T_{G}\SizeTruncated{m}$ then all its descendants are not, so $\B T_{G}\SizeTruncated{m}$ is a subtree of $\B T_{G}$.

In the rest of Section~\ref{se:rich}, we would need to work with \textbf{unordered} trees where the $t$ children of a non-leaf vertex are not ordered and the first part $\B i$ of each node $(\B i, x)$, that can be used to order children, is erased (but we keep the second part $x$). For these objects we will use the (non-bold) symbol~$T$ instead of~$\B T$.
In particular, the \textbf{size-$n_0$ truncated tree} $T_{G}\SizeTruncated{n_0}$ is the unordered subtree of the  tree of the $P_I$-mixing construction $G$, where we keep only those nodes that corresponds to parts of size at least~$n_0$. Suppose now that $G\in\Sigma P_I$ is maximum. Then every branch of $T_{G}\SizeTruncated{n_0}$ has length at least, say, $\log_{t+1}(n/n_0)$, because when we subdivide any set $V$ with at least $n_0$ vertices, its partition is balanced by Lemma~\ref{lm:OptHRecExact}  and thus the smallest part size is at least $\lfloor |V|/t\rfloor\ge |V|/(t+1)$. Thus any resulting tree has at least $s\coloneqq t^{\log_{t+1}(n/n_0)-1}$ non-leaf vertices. Since the children of every non-leaf 
can be labelled by any element of $I$ (as long as we use the same element for all children), we have at least $|I|^s$ choices here.
The number of ways that an isomorphic copy of any feasible $I$-labelled rooted tree $T$ can be generated as above, rather roughly, is most $(t!)^s$. Thus there are at least $(|I|/t!)^s$ non-isomorphic (unordered $I$-labelled rooted) trees obtainable this way. This is exponential in $n$ (since $n_0$ is fixed); thus the first part of Theorem~\ref{th:ExpMany} follows from the following claim.

\bcl{IsoCP}{If $G,G'\in \C P$ are isomorphic $3$-graphs then their $n_0$-truncated unordered trees $T_{G}\SizeTruncated{n_0}$ and $T_{G'}\SizeTruncated{n_0}$ are isomorphic.}

\bcpf We use induction on $n$, the number of vertices in $G$. If $n<n_0$, then the truncated trees of $G$ and $G'$ are empty and thus the conclusion vacuously holds. Suppose $n\ge n_0$ and that the identity map on $[n]$ gives an isomorphism between some $G,G'\in\C P$.

By Lemmas~\ref{lm:regular} and~\ref{lm:MaxRigid}, the maximum $P_I$-mixing construction $G$ is rigid. (In fact, the proof of Lemma~\ref{lm:MaxRigid} simplifies greatly in this case and every $P_I$-mixing construction $G$ with all bottom parts non-empty is rigid since we can identify bottom edges of $G$ as precisely those that do not contain the symmetric difference of some two distinct edges of $G$.)

 Thus, $G$ and $G'$ have the same base $P_D$ and the isomorphism between their bottom edge-sets, $\blow{\OO D}{V_1,\ldots,V_t}\subseteq G$ and $\blow{\OO D}{V_1',\ldots,V_t'}\subseteq G'$, comes from an automorphism $h$ of $D$, that is, $V_{h(i)}'=V_i$ for every $i\in [t]$. Now, using the automorphism $h$ relabel the bottom parts of $G'$ so that $V_i'=V_i$ for $i\in [t]$ and apply induction to each pair $G[V_i],G'[V_i]\in \C P$ of isomorphic (in fact, equal) $3$-graphs.\ecpf

We now turn to the second part of the theorem. Take any $s$. Let $\ell$ satisfy $(|I|/t!)^{\ell-1}\ge s$. Fix  sufficiently small $\varepsilon>0$ and let $n\to\infty$. It is enough to find $s$ maximum $P_I$-mixing constructions $\mathcal{G}_{1}, \ldots, \mathcal{G}_{s}$ on $[n]$ so that every two are at edit distance at least $\varepsilon^3 n^3$, as this will demonstrate that $\xi(\mathcal{F}) \ge s$. 
Indeed, suppose on the contrary that $\xi(\mathcal{F}) \le s-1$. Then there exist $s-1$ $3$-graphs $\mathcal{H}_{1}, \ldots, \mathcal{H}_{s-1}$ on $[n]$ such that, for each $i \in [s]$, there is some $j \in [s-1]$ where the edit distance between $\mathcal{G}_i$ and $\mathcal{H}_j$ is at most $\varepsilon^3 n^3/3$. 
By the Pigeonhole Principle, there must be two $3$-graphs, say $\mathcal{G}_{i_1}$ and $\mathcal{G}_{i_2}$, that are both within edit distance of at most $\varepsilon^3 n^3/3$ to the same $\mathcal{H}_{j}$. 
However, by the triangle inequality, this implies that the edit distance between $\mathcal{G}_{i_1}$ and $\mathcal{G}_{i_2}$ is at most $2 \varepsilon^3 n^3/3 < \varepsilon^3 n^3$, a contradiction. 

Call two $I$-labelled unordered trees $T$ and $T'$ \textbf{isomorphic up to level $m$} if their first $m$ levels span isomorphic trees. For example, two trees are isomorphic  up to level 1 if and only if the children of the root are labelled by the same element of $I$ in the both trees. (Recall that the root always get label $\emptyset$ while all children of a node always get the same label.)
For convenience, let us agree that every two trees are isomorphic up to level~0.

By our choice of $\ell$, there are $s$ trees
that are pairwise non-isomorphic up to level~$\ell$. Since $n$ is sufficiently large, we can assume by Lemma~\ref{lm:F(D)} that for each of these trees $T$ there is a maximum $P_I$-mixing construction $G$ on $[n]$ whose unordered tree is isomorphic to $T$ up to level~$\ell$, that is, $T_G\TruncatedLevel{\ell}\cong T$.
Thus it is enough to show that any two of the obtained $3$-graphs, say $G$ and $G'$ with unordered trees $T$ and $T'$ respectively where $T'$ and $T'$ are non-isomorphic up to level $\ell$, are at the edit distance at least~$\varepsilon^3 n^3$. Suppose that this is false and the identity map exhibits this, that is, $|G\bigtriangleup G'|< \varepsilon^3 n^3$. Let $\B V$ and $\B V'$ be the partition structures of $G$ and $G'$ respectively. Also, for a tree $T$ and a sequence $(i_1,\dots,i_m)$, let $T[i_1,\dots,i_m]$ be the subtree formed by the vertex with the first coordinate $(i_1,\dots,i_m)$ and all its descendants in $T$.

\bcl{FarInEdit}{For every $m\le \ell$, there are sequences $(i_1,\ldots,i_m),(i_1',\ldots,i_m')\in [t]^m$ such that $T[i_1,\dots,i_m]$ and $T'[i_1',\ldots,i_m']$ are non-isomorphic up to level $\ell-m$ and
\begin{align*}
 \left|V_{i_1,\ldots,i_m}\bigtriangleup V'_{i_1',\ldots,i_m'}\right|\le 2(t-1)m\cdot \varepsilon n.
 \end{align*}
 }

\bcpf We use induction on $m$ with the base case $m=0$ being satisfied by the empty sequences. Let $m\ge1$ and suppose that we already constructed some sequences $(i_1,\ldots,i_{m-1}),(i_1',\ldots,i_{m-1}')\in [t]^{m-1}$ that satisfy the claim for~$m-1$. Let $U\coloneqq V_{i_1,\ldots,i_{m-1}}$ and $U'\coloneqq V'_{i_1',\ldots,i_{m-1}'}$. For $j\in [t]$, let $U_j\coloneqq V_{i_1,\ldots,i_{m-1},j}$ and $U'_j\coloneqq V'_{i_1',\ldots,i_{m-1}',j}$. In other words, $U=U_1\cup\dots\cup U_t$ and $U'\coloneqq U'_1\cup\dots\cup U'_t$ are the bottom partitions of the  $P_I$-mixing constructions $G[U]$ and $G'[U']$ respectively. These constructions are maximum; thus, by Lemma~\ref{lm:OptHRecExact}, their bottom partitions are balanced and, by simple induction on $m$, each part has size $n/t^{m}+O(1)$.

For each $i\in [t]$, there cannot be distinct $j,k\in [t]$ with each of $A\coloneqq U_i\cap U'_j$ and $B\coloneqq U_i\cap U'_k$ having at least $\varepsilon n$ elements. Indeed, otherwise the co-degree of each pair $(a,b)\in A\times B$ is at most $|U_i|-2=n/t^{m}+O(1)$ in $G[U]$ as all such edges lie inside $U_i$, and at least $(t-3)n/t^{m}+O(1)$ in $G'[U']$ (since the bottom pattern of $G'[U']$ has $t-3$ triples containing the pair $\{i,j\}$). This is impossible, since then
 \begin{align}\label{eq:Edit1}
 |G\bigtriangleup G'|\ge |G'\setminus G|\ge |A|\cdot |B|\cdot \left(\frac{(t-4)n}{t^{m}}+O(1)\right) \ge \varepsilon^3 n^3.
 \end{align}

 Thus there is $h:[t]\to [t]$ with $|U_i\cap U'_j|<\varepsilon n$ for each $j\in [t]$ different from~$h(i)$.
 Since $\e\ll 1/t^\ell\le 1/t^m$, we have
 \begin{align}\label{eq:U'H}
  |U_i\cap U'_{h(i)}|\ge |U_i|-(t-1)\varepsilon n -|U'\setminus U| \ge \frac{n}{t^{m}}-2(t-1)m\cdot \varepsilon n+O(1),
 \end{align}
 which is strictly larger than half of $|U'_{h(i)}|$.
 Thus $h$ is an injection and, by the equal finite cardinality of its domain and its range, a bijection.


The map $h$ has to be an isomorphism between the base patterns $D$ and $D'$ of $G[U]$ and $G'[U']$ respectively. Indeed,
if $h$ does not preserve the adjacency for at least one triple then \eqref{eq:U'H} implies that, say, $|G\bigtriangleup G'|\ge (n/(2t^{m}))^3>\varepsilon^3 n^3$, a contradiction. 

Since the trees $T[i_1,\dots,i_{m-1}]$ and $T'[i_1',\ldots,i_{m-1}']$ (which are the trees of $G[U]$ and $G'[U']$) are non-isomorphic up to level $\ell-m+1$, there must be $i$ such that, letting $i_m\coloneqq i$ and $i_m'\coloneqq h(i)$,
the trees $T[i_1,\dots,i_{m}]$ and $T'[i_1',\ldots,i_{m}']$ are non-isomorphic up to level $\ell-m$. Also, by~\eqref{eq:U'H} (and its version when the roles of $G$ and $G'$ are exchanged), we have that
 \begin{align*}
 \left|V_{i_1,\dots,i_{m}} \bigtriangleup V'_{i_1',\ldots,i_{m}'}\right|\le |U\bigtriangleup U'| +2(t-1)\cdot \varepsilon n
 \le 2(t-1)m\cdot \varepsilon n.
 \end{align*}
Thus the obtained sequences $(i_1,\ldots,i_m),(i_1',\ldots,i_m')\in [t]^m$ satisfy the claim.\ecpf

Finally, a desired contradiction follows by taking $m\coloneqq \ell$ in the above claim (as every two trees are isomorphic up to level 0), finishing the proof of the second part of Theorem~\ref{th:ExpMany}.\epf

\section{Feasible region}\label{SEC:feasible-region}
In this section we prove Theorem~\ref{THM:feasibe-region}. We need some auxiliary results first.

\begin{lemma}\label{LEMMA:minimal-shadow-complete}
Suppose that $P = (m,E,R)$ is a minimal pattern.
Then every pair $\{i,j\} \in \binom{[m]}{2}$ is contained in some multiset in $E$.
\end{lemma}
\bpf
Suppose to the contrary that there exists a pair in $[m]$ that is not contained in any multiset in $E$. By relabeling the vertices in $P$, we may assume that this pair is $\{1,2\}$. 
Let $\lambda\coloneqq \lambda_P$. For $\B z\in \I S_m$, let $f(\B z)\coloneqq \lambda_{E}(\B z)+\lambda\sum_{t\in R}z_t^{r}$. Note that no term in $f(\B z)$ contains the product $z_1z_2$.
So we can rewrite $f(\B z)$ as
\begin{align*}
    f(\B z) 
    = \sum_{i=1}^{r}\left( \alpha_i z_1^{i} + \beta_i z_2^{i} \right) + \gamma,
\end{align*}
where $\alpha_i, \beta_i, \gamma \ge 0$ all depend only on $z_3, \ldots, z_{m}$.

Let $\B x\in \C X$ be an optimal vector for $P$. 
By symmetry, we may assume that $\sum_{i=1}^{r} \alpha_i x_1^{i-1} \ge \sum_{i=1}^{r} \beta_i x_2^{i-1}$. 
Let $\B y \coloneqq (x_1+x_2, 0, x_3, \ldots, x_m) \in \mathbb{S}_{m}$. 
Notice that 
\begin{align*}
    f(\B y) - f(\B x)
    & = \sum_{i=1}^{r} \alpha_i (x_1 + x_2)^{i} - \sum_{i=1}^{r}\left( \alpha_i x_1^{i} + \beta_i x_2^{i} \right) \\
    & = (x_1+x_2) \sum_{i=1}^{r} \alpha_i (x_1 + x_2)^{i-1} - \left( x_1 \sum_{i=1}^{r} \alpha_i x_1^{i-1} + x_2 \sum_{i=1}^{r} \beta_i x_2^{i-1} \right) \\
    & \ge (x_1+x_2) \sum_{i=1}^{r} \alpha_i (x_1 + x_2)^{i-1} - \left( x_1 \sum_{i=1}^{r} \alpha_i x_1^{i-1} + x_2 \sum_{i=1}^{r} \alpha_i x_1^{i-1} \right) \\
    & = (x_1+x_2) \sum_{i=1}^{r} \left(\alpha_i (x_1 + x_2)^{i-1} - \alpha_i x_1^{i-1} \right)
    \ge 0.
\end{align*}
Since $\B x$ is an optimal vector for $P$, we must have $f(\B y) - f(\B x) = 0$, implying that $\B y$ is an optimal vector for $P$ or that  $\B y=\B e_1$ is the first standard basis vector. In either case, one can easily derive that $\lambda_P=\lambda_{P-2}$,  contradicting the minimality of $P$.
\epf

For every $i\in I$ and pattern $P_i = (m_i,E_i,R_i)$ we let $S_i \subseteq [m_i]\setminus R_i$ be the collection of $j\in [m_i]\setminus R_i$ such that
for every blowup $\blow{E_i}{V_1, \ldots, V_{m_i}}$ of $E_i$ the shadow of $\blow{E_i}{V_1, \ldots, V_{m_i}}$ has no edge inside~$V_j$.

\begin{lemma}\label{LEMMA:stability-shadow}
For every  finite collection $P_I$ of minimal $3$-graph patterns with the same Lagrangian $\lambda \in (0,1)$, there is an integer $M$ such that for every $\varepsilon>0$ there exist $\delta >0$ and $N_0 >0$ such that the following holds for all $n\ge N_0$.
If $G$ is an $\mathcal{F}_{M}$-free $3$-graph on $[n]$ with at least $(\lambda-\delta) \binom{n}{3}$ edges,
then there exists a $P_I$-mixing construction $H$ on $[n]$ with base $P_i$ and bottom partition $V(G) = V_1\cup \cdots \cup V_{m_i}$ for some $i\in I$ such that
$\delta(H) \ge (\lambda-\varepsilon)\binom{n-1}{2}$,
$|H \bigtriangleup G|\le \varepsilon\binom{n}{3}$, $(|V_1|/n,\dots,|V_{m_i}|/n)$ is $\varepsilon$-close to a $P_i$-optimal vector, and
$\sum_{j\in S_i}|(\partial G)[V_j]| \le \varepsilon \binom{n}{2}$.

\end{lemma}
\bpf Many steps of this proof are similar to the analogous parts of the proof of Lemma~\ref{lm:Max2}. So we may be brief when the appropriate adaptation is straightforward,

Let $\ell_0$ be the constant returned by Lemma~\ref{lm:tight} and then let $M$ be sufficiently large. Given $\varepsilon>0$, choose sufficiently small positive constants in this order:  $\delta_4\gg \delta_3\gg \delta_2\gg \delta_1\gg \delta$. Let $n\to\infty$ and let $G$ be any $\C F_M$-free $r$-graph on $V\coloneqq [n]$ with at least $(\lambda-\delta)\binom{n}{3}$ edges.

Due to Theorem~\ref{THM:mixed-pattern}\ref{it:b}, we may assume that
$|G\bigtriangleup H|\le \delta_1 \binom{n}{3}$ for some $P_{I}$-mixing construction~$H$.
Let $H$ have the partition structure $\B V$, the base $P_i$ and the bottom partition $V_1,\ldots,V_{m_i}$.
By modifying $H$ as in the argument leading to~\eqref{eq:MinDegStage1}, we can further assume that
\begin{align}\label{eq:HShadow}
\delta(H) \ge (\lambda-\delta_2)\binom{n-1}{2}\quad\mbox{and}\quad |G\bigtriangleup H|\le \delta_2 \binom{n}{3}.
\end{align}
By Lemma~\ref{lm:OmegaN}\ref{it:OmegaNx}, the bottom part ratios are within $\varepsilon$ from a $P_i$-optimal vector. Thus this $H$ satisfies the first three properties stated in the lemma.
The rest of the proof is dedicated to proving the remaining property that the total shadow of $G$ inside the parts $V_j$ with $j\in S_i$ is ``small''.

Let
\begin{align}
    Z_{1} \coloneqq  \left\{v\in V\colon d_{G}(v)\le (\lambda-\delta_3)\binom{n-1}{2}\right\}\quad\text{and}\quad  G_1 \coloneqq  G-Z_1. \notag
\end{align}

By~\eqref{eq:HShadow}, we have the following inequalities.

\bcl{Z1}{It holds that
\begin{align*}
|Z_1| &\le \frac{3|G\bigtriangleup H|}{(\delta_3-\delta_2)\binom{n-1}{2}}\ \le\
\delta_3 n,\\
|G_1| &\ge |G|-|Z_1|\binom{n-1}{2}\ge (\lambda-\delta)\binom{n}{3}-\delta_3 n \binom{n-1}{2}\ \ge\ \lambda\binom{n}{3} - \delta_3 n^3,\\
 \delta(G_1) &\ge  (\lambda-\delta_3)\binom{n-1}{2} -|Z_1| n\ \ge\ (\lambda-4\delta_3)\binom{n}{2}.\ecpf
 \end{align*}
}

Let $V' \coloneqq  V\setminus Z_1$,
and let $H_1$ be the induced subgraph of $H$ on~$V'$.
Clearly, we have $|G_1\bigtriangleup H_1|\le |G\bigtriangleup H|\le \delta_2 \binom{n}{3}$.

Define
\begin{align}
    Z_{2} \coloneqq  \left\{v\in V'\colon |L_{G_1}(v)\cap L_{H_1}(v)|\le (\lambda-\delta_4)\binom{n}{2}\right\}
    \quad\text{and}\quad G_2 \coloneqq  G_1-Z_2. \notag
\end{align}

Similarly to Claim~\ref{cl:Z1}, we have the following.

\bcl{Z2}{We have $|Z_2| \le \delta_4 n$,  $|G_2| \ge \lambda\binom{n}{3} - \delta_4 n^{3}$ and $\delta(G_2) \ge (\lambda-4\delta_4)\binom{n}{2}$.\ecpf}

Let $V'' \coloneqq  V'\setminus Z_2$ and $V_{\B i}'' \coloneqq  V_{\B i}\cap V''$ for every legal sequence $\B i$.
Let $H_2$ be the induced subgraph of $H$ on $V''$.
Similarly to above, we have $|G_2 \bigtriangleup H_2|\le |G\bigtriangleup H|\le \delta_2 n^{3}$ and
 \begin{align}\label{eq:delta(H2)}
 \delta(H_2)
 \ge \delta(H)-\left(|Z_1|+|Z_2|\right)\binom{n-1}{2}
 \ge (\lambda-3\delta_4)\binom{n}{2}.
 \end{align}
%
%
In addition, it follows from the definition of $Z_2$ that for every $v \in V''$,
\begin{align}\label{eq:delta-H2-G2-intersect}
    \left| L_{G_2}(v) \cap L_{H_2}(v) \right|
    \ge (\lambda - \delta_4) \binom{n}{2} - |Z_2| (n-2) 
    \ge (\lambda - 3 \delta_4) \binom{n}{2}. 
\end{align}

Take any $j \in S_i$. Our next aim is to show that $\partial G \cap \binom{V_j''}{2} = \emptyset$.
Suppose to the contrary that there exists $D= \{u_1,u_2,u_3\}$ in $G$ with $u_1,u_2\in V_j''$.

Let $H_3$ be obtained from $H_2$ by removing edges contained in $H_2[V''_{\B i}]$ for all legal sequences $\B i$ of length at least~$\ell_0$.
Let $\B T\coloneqq \B T_{H_3}$. By the min-degree property of $H$ and Lemma~\ref{lm:OmegaN}\ref{it:OmegaNVi}, each part of $\BV{H}$ at level at most $\ell_0$ has at least $(\beta/2)^{\ell_0}n$ vertices. This is much larger than $|Z_1\cup Z_2|\le (\delta_3+\delta_2)n$, the number of the removed vertices when building $H_3$ from $H$. Thus, in particular, $\B T_{H_3}=\B T_{H}\TruncatedLevel{\ell_0}$.
Note that if the tree $\B T$ is extendible then it is maximal up to level~$\ell_0$. Indeed, each involved recursive part $V_{\B i}$ of $H$ has at least $(\beta/2)^{\ell_0}n$ vertices and thus must span some edges, for otherwise the edge density of the $P_I$-mixing construction $H$, which is $\rho(H)\ge \rho(G)-\delta_2\ge \lambda-\delta-\delta_2$, can be increased by \begin{align*}
(\lambda-o(1))\binom{|V_{\B i}|}{3}/\binom{n}{3}\ge \lambda(\beta/2)^{3\ell_0}-o(1)\gg \delta+\delta_2,
\end{align*}
thus jumping above
the maximum density $\lambda+o(1)$, a contradiction.
Let $F$ be the rigid $P_{I}$-mixing construction with $\B T_F=\B T$, returned by Lemma~\ref{lm:tight} or~\ref{LEMMA:non-extendable-partition}.
Since $M\gg \ell_0$, we can assume that $F$ has at most $M-1$ vertices.

By relabelling $V(F)$ we can assume that $u_1,u_2\in V(F)$ and these vertices have the same branch in $F$ and $H$ (which is just the single-element sequence $(j)$) while $u_3\not\in V(F)$. Let $F^D$ be obtained from $F$ by adding the edge $D$ (and the new vertex~$u_3$). Let $\B W$ be the partition structure of $F$.

\bcl{FDShadow}{$F^D \in \mathcal{F}_{M}$}

\bcpf Suppose to the contrary that there exists an embedding $f$ of $F^D$ into some $P_{I}$-mixing construction $Q$ with base $P_t$ and bottom partition $U_1,\dots,U_{m_t}$. We can assume that $f(V(F^D))$ intersects at least two parts~$U_s$. First, suppose that already the image of $V(F)\subseteq V(F^D)$ under $f$ intersects at least two bottom parts of~$Q$. By the rigidity of $F$, we have that $t=i$ and,
by relabelling the parts $U_s$, we can assume that $W_s\subseteq U_s$ for every $s\in [m_i]$. Thus $f(\{u_1,u_2\})\subseteq U_{j}$. By $j\in S_i$,
no edge of $Q$ can cover the pair $\{u_1,u_2\}$, a contradiction to $f$ mapping $D$ to an edge of~$Q$. Thus, suppose that $f(V(F))\subseteq U_s$ for some $s\in [m_t]$. Since $|F|>0$, we have $s\in R_t$.
By $V(F^D)=V(F)\cup\{u_3\}$, it holds that $f(u_3)\in U_{k}$ for some $k\not=s$. However, then the profile of the edge $f(D)$ in $Q$ is $\multiset{\rep{s}{2},k}$,  contradicting by Lemma~\ref{lm:density=1} our assumption that $\lambda_{P_t}=\lambda<1$.\ecpf

Now, we apply the familiar argument where we consider all maps $f:V(F^D)\to V(H_3)$ such that $f$ is the identity on $D$ and, for every vertex $u$ of $F^D$ different from $u_3$, it holds that $\branch{F}{u}=\branch{H_3}{f(u)}$, that is, the branch of $f(u)$ in $H_3$ is the same as the branch of $u$ in~$F$. We have by Lemma~\ref{lm:OmegaN}\ref{it:OmegaNVi} that for every legal sequence $\B i$ of length at most $\ell_0$
we have $|V''_{\B i}| \ge (\beta/2)^{\ell_0}|V''|$. Thus the number of choices of $f$ is at least, say,
$((\beta/2)^{\ell_0}n/3)^{v(F)-2}$. By Claim~\ref{cl:FDShadow}, each choice of $f$ reveals a missing edge $Y
\in H_3\setminus G_2$. It is impossible that $Y$ is disjoint from $\{u_1,u_2\}$ for at least half of the choices, for otherwise the size of $H_3\setminus G_2$ is too large, since every such $Y$ is overcounted at most $n^{v(F)-5}$ times. Thus at least half of the time we have $Y\cap\{u_1,u_2\}\not=\emptyset$. By $j\in S_i$, each such $Y$ intersects $\{u_1,u_2\}$ in exactly one vertex. Note that, for each $j\in \{1,2\}$, by~\eqref{eq:delta-H2-G2-intersect} and Lemma~\ref{LEMMA:pattern-max-degree}, we have that
\begin{align*}
    |L_{H_3}(u_j)\setminus L_{G_2}(u_j)|
    & \le |L_{H_2}(u_j)\setminus L_{G_2}(u_j)| \\
    & = |L_{H_2}(u_j)| - |L_{G_2}(u_j) \cap L_{H_2}(u_j)| \\
    & \le  \Delta(H) - |L_{G_2}(u_j) \cap L_{H_2}(u_j)| 
     \le (\lambda+\delta_3)\binom{n}{2} - (\lambda- 3\delta_4)\binom{n}{2}
    \le 4\delta_4 n^2.
\end{align*}
Thus the total number of such choices of $f$ can be upper bounded by $2 \times 4 \delta_4 n^2 = 8 \delta_4 n^2$, the number of choices of $Y$, times the trivial upper bound $n^{v(F)-4}$ on the number of extensions to the remaining vertices of $V(F)$. This contradicts that $\delta_4$ was sufficiently small depending on $\beta$ and~$\ell_0$.

We have shown that, for every $j\in S_i$, the set $V_j''$ spans no edges in $\partial G$, that is, every edge of $\partial G$ in $V_{j}$ must contain at least one vertex from $Z_1\cup Z_2$.
Thus $\sum_{j\in S_i}|(\partial G) [V_j]| \le |Z_1\cup Z_2|n \le \varepsilon \binom{n-1}{2}$.
This proves Lemma~\ref{LEMMA:stability-shadow}.
\epf

Now we are ready to prove Theorem~\ref{THM:feasibe-region}.

\bpf[Proof of Theorem~\ref{THM:feasibe-region}.]
Again, we can assume that the assumptions and terminology of Section~\ref{se:assumptions} apply to~$P_I$.

Let $M$ be a sufficiently large integer, in particular, such that $M$ satisfies Lemma~\ref{LEMMA:stability-shadow} and the family $\C F_M$ satisfies Theorem~\ref{THM:mixed-pattern}. Let us show that $\C F\coloneqq \C F_{M}$ satisfies
Theorem~\ref{THM:feasibe-region}. Theorem~\ref{THM:mixed-pattern} implies that $\pi(\C F_{M})=\pi(\C F_\infty)=\lambda$, from which it easily follow that $M(\C F_\infty)\subseteq M(\C F_{M})$. Let us show the converse implication.

Fix any $x$ in $M(\mathcal{F}_{M})$ and $\e>0$. We have to approximate the point $(x,\lambda)\in \Omega(\C F_M)$ by an element of $\Omega(\C F_\infty)$ within $\e$ in, say, the supremum norm.
Let $\left(G_{n}\right)_{n=1}^{\infty}$ be a sequence of $\mathcal{F}_{M}$-free $3$-graphs that realizes $(x, \lambda)$. Let $v_n\coloneqq v(G_n)$ for $n\ge 1$. By passing to
a subsequence if necessary, we can assume that the sequence $(v_n)_{n\in \I N}$ is
strictly increasing.
Since $\lim_{n\to \infty}\rho(G_n) = \pi(\mathcal{F}_{M})$, Lemma~\ref{LEMMA:stability-shadow} applies for all large $n$ and returns a $P_{I}$-mixing construction $H_n$ on $V(G_n)$. In particular, it holds that $\delta(H_n) \ge (\lambda-o(1))\binom{v_n-1}{2}$
and $|G_n \triangle H_n| =o(v_n^3)$ as $n\to \infty$.
Since $H_n$ is a $P_{I}$-mixing construction and $\lim_{n\to \infty}\rho(H_n) = \lim_{n\to \infty}\rho(G_n) = \pi(\mathcal{F}_{M})$,
we have $\lim_{n\to \infty}\rho(\partial H_n) \in M(\C F_\infty)$.

Let $\ell$ be a sufficiently large integer. Define $\delta_{\ell+1} \coloneqq  (1-\lambda/(2r))^{\ell+1}/2$ and then let $\delta_{\ell}\gg \dots\gg \delta_1$ be sufficiently small positive constants, chosen in this order. Let $n\to\infty$.
Let $\B V$ be the partition structure of~$H_n$.
By choosing $n$ sufficiently large we can assume that $\delta(H_n) \ge (\lambda-\delta_1)\binom{v_n-1}{2}$,
$|G_n \triangle H_n| \le \delta_1 \binom{v_n}{3}$, and $|G_n| \ge (\lambda-\delta_1) \binom{v_n}{3}$.

We call edges in $\partial G_n \setminus \partial H_n$ \textbf{bad shadows}.
By Lemma~\ref{LEMMA:minimal-shadow-complete}, every edge $e \in \partial G_n \setminus \partial H_n$ is contained entirely inside some bottom part of $H_n$. Let $\mathrm{br}_{\B V}(e)$ denote the (unique, non-empty) maximal sequence $\B i$ such that $e\subseteq V_{\B i}$. 

Every bad shadow with branch of length at least $\ell+1$ lies inside some part $V_{\B i}$ with $|\B i|\ge \ell+1$, whose size is at most $(1-\lambda/(2r))^{\ell+1}v_n$ by Lemma~\ref{lm:OmegaN}\ref{it:OmegaNVi}. By the Hand-Shaking Lemma, the number of such pairs is at most
$\frac12\, (1-\lambda/(2r))^{\ell+1}v_n\cdot v_n=
\delta_{\ell+1} v_n^2$.

Recall that each part of height at most $\ell$ in $H_n$ has size at least $(\beta/2)^\ell v_n$ by Lemma~\ref{lm:OmegaN}\ref{it:OmegaNVi}. In particular, by Lemma~\ref{LEMMA:minimal-shadow-complete}, the collection of all bad shadows whose branch has length $1$ is exactly the set $\bigcup_{j\in S_{b(H_n)}}\partial G_n [V_j]$,
and by Lemma~\ref{LEMMA:stability-shadow}, we have that
 \begin{align*}\sum_{j\in S_{b(H_n)}}|\partial G_n [V_j]| \le \delta_2 v_n^2.\end{align*}
Now, repeat the following for every recursive index $j_1\in R_{b(H_n)}$.
Since $\delta(H_n) \ge (\lambda-\delta_1)\binom{v_n-1}{2}$, we have by Lemma~\ref{lm:OmegaN} that $|V_{j_1}| \ge \beta v_n/2$ and $\delta(H_n[V_{j_1}]) \ge (\lambda-\delta_2)\binom{|V_{j_1}|-1}{2}$.
Since $|G_n[V_{j_1}] \triangle H_n[V_{j_1}]|\le |G_n \triangle H_n| \le \delta_1 \binom{v_n}{3}$, we have
\begin{align*}
|G_n[V_{j_1}]| \ge |H_n[V_{j_1}]| - \delta_1 \binom{v_n}{3} \ge (\lambda-2\delta_2)\binom{|V_{j_1}|}{3}.
\end{align*}
So applying Lemma~\ref{LEMMA:stability-shadow} to $H_n[V_{j_1}]$, we obtain that, for example,
\begin{align*}
 \sum_{j\in S_{b(H_n[V_{j_1}])}}|\partial G_n [V_{j_1, j}]| \le \frac1m\,{\delta_3 \binom{v_n}{2}},
\end{align*}
where $m \coloneqq  \max\{m_k\colon k\in I\}$.
Therefore,
\begin{align*}
\sum_{j_1\in R_{b(H_n)}}\sum_{j\in S_{b(H_n[V_{j_1}])}}|\partial G_n [V_{j_1, j}]| \le m\frac{\delta_3 \binom{v_n}{2}}{m} = \delta_3 \binom{v_n}{2},
\end{align*}
 that is, the number of bad shadows whose branch has length $2$ is at most~$\delta_3 \binom{v_n}{2}$.
Repeating this argument we can show that the the number of bad shadows whose length of branch is $h \le \ell$ is at most $\delta_{h+1} \binom{v_n}{2}$.
Therefore, the total number of bad shadows is bounded by
\begin{align}
\delta_2 \binom{v_n}{2} + \delta_3 \binom{v_n}{2}  + \cdots + \delta_{\ell}\binom{v_n}{2} + \delta_{\ell+1}\binom{v_n}{2} \le 2 \delta_{\ell+1} \binom{v_n}{2}. \notag
\end{align}
Therefore, $|\partial G_n \setminus \partial H_n|/\binom{v_n}{2} < \e/2$. On the other hand, every pair $xy\in\partial H_n$ at level at most $\ell$ is covered by at least $(\beta/2)^{\ell}v_n$ triples in $H_n$. This observation and our previous estimate of the total number of pairs at levels at least $\ell+1$ give that
 \begin{align*}
 |\partial H_n\setminus \partial G_n|\le \delta_{\ell+1}v_n^2 + \frac{3|G_n\bigtriangleup H_n|}{(\beta/2)^{\ell}v_n}< \frac{\e}2\, \binom{v_n}{2}.
\end{align*}
Since $\e>0$ was arbitrary, we conclude that $x\in M(\C F_M)$.
This proves Theorem~\ref{THM:feasibe-region}.
\epf
\section{Proof of Corollary~\ref{CORO:feasible-function-Cantor}}\label{SEC:examples}
In this section we prove Corollary~\ref{CORO:feasible-function-Cantor} by applying Theorem~\ref{THM:feasibe-region} to two specific patterns.

Consider the following two patterns $P_1 = (5,K_{5}^{3}, \{1\})$ and $P_2 = (7,B_{5,3}, \{1\})$, where the $3$-graph $B_{5,3}$ was defined in Section~\ref{example}.

\begin{lemma}\label{LEMMA:Lagrangian-P1-P2}
The following statements hold.
\begin{enumerate}[label=(\alph*)]
\item\label{it:Lagrangian-B53-K53} We have $\lambda_{P_1} = \lambda_{P_2} = \lambda\coloneqq  \frac{3\left(\sqrt{7}-2\right)}{4} \approx 0.484313$.
\item\label{it:optimal-vector-K53} For $\B x\in\IS_5^*$,
we have $\lambda_{K_{5}^{3}}(\B x)+ \lambda x_1^3 = \lambda$ if and only if
\begin{align}\label{equ:Lagrange-maxmizer-1}
x_1 = \frac{\sqrt{7}-2}{3}, \quad\text{and}\quad
x_2  = \cdots = x_5 = \frac{5-\sqrt{7}}{12}.
\end{align}
\item\label{it:optimal-vector-B53} For $\B y\in\IS_7^*$,
we have $\lambda_{B_{5,3}}(\B y)+ \lambda y_1^3 = \lambda$ if and only if
\begin{align}\label{equ:Lagrange-maxmizer-2}
y_1 = \frac{\sqrt{7}-2}{3},\quad
y_2 = y_3 = \frac{5-\sqrt{7}}{12}, \quad\text{and}\quad
y_4 = \cdots = y_7 = \frac{5-\sqrt{7}}{24}.
\end{align}
\end{enumerate}
\end{lemma}

\noindent\textbf{Remark.}
Parts~\ref{it:optimal-vector-K53} and~\ref{it:optimal-vector-B53} imply that $P_1$ and $P_2$ are minimal.

\bpf
For $\B x\in \mathbb{S}_{5}^*$ and $\B y\in \mathbb{S}_{7}^*$, 
let
\begin{align}
g_1(\B x) := \frac{\lambda_{K_{5}^3}(\B x)}{1-x_1^3} \quad\text{and}\quad
g_2(\B y) := \frac{\lambda_{B_{5,3}}(\B y)}{1-y_1^3}. \notag
\end{align}
It follows from the AM-GM Inequality that
\begin{align}
g_1(\B x)
 = \frac{\lambda_{K_{5}^3}(\B x)}{1-x_1^3}
& = \frac{6\left(x_1\sum_{ij\in \binom{[2,5]}{2}}x_ix_j + \sum_{ijk\in \binom{[2,5]}{3}}x_ix_jx_k\right)}{1-x_1^3} \notag\\
& \le \frac{6\left(x_1\binom{4}{2}\left(\frac{1-x_1}{4}\right)^2 + \binom{4}{3}\left(\frac{1-x_1}{4}\right)^3\right)}{1-x_1^3} \notag\\
& = \frac{3(1-x_1)(1+5x_1)}{8(1+x_1+x_1^2)}
\le \frac{3(\sqrt{7}-2)}{4}, \notag
\end{align}
where the equality holds if and only if~(\ref{equ:Lagrange-maxmizer-1}) holds.

Similarly, for $g_2(\B y)$ we have
\begin{align}
\lambda_{B_{5,3}}(\B y)
& = 6y_1\left(y_2y_3+(y_2+y_3)(y_4+y_5+y_6+y_7)+(y_4+y_5)(y_6+y_7)\right) \notag\\
& \quad  +6\left( y_2y_3(y_4+y_5+y_6+y_7)+ y_2(y_4+y_6)(y_5+y_7) + y_3(y_4+y_7)(y_5+y_6)\right).  \notag
\end{align}
Notice that
\begin{align}
& y_2y_3+(y_2+y_3)(y_4+y_5+y_6+y_7)+(y_4+y_5)(y_6+y_7) \notag\\
& = y_2y_3+y_2(y_4+y_5)+y_2(y_6+y_7)+ y_3(y_4+y_5)+y_3(y_6+y_7)+(y_4+y_5)(y_6+y_7) \notag\\
& = \sigma_2(y_2, y_3, y_4+y_5, y_6+y_7)
 \le \binom{4}{2}\left(\frac{y_2+ \cdots + y_7}{4}\right)^2
 = \frac{3}{8} \left(1-y_1\right)^2, \notag
\end{align}
where
\begin{align*}
\sigma_k(x_1,\dots,x_s)\coloneqq \sum_{X \in \binom{[s]}{k}}\prod_{i\in X} x_i
\end{align*} is the \textbf{$i$-th symmetric polynomial}, and the last inequality follows from the Maclaurin Inequality (see e.g.~\cite[Theorem~11.2]{Cvetkovski12i}) that $\sigma_k(x_1,\dots,x_s)\le \binom{s}{k} (x_1+\dots+x_s)^k/s^k$ for any non-negative~$x_i$'s and $k\in\I N$.

On the other hand, it follows from the AM-GM and Maclaurin Inequalities that
\begin{align}
& y_2y_3(y_4+y_5+y_6+y_7)+ y_2(y_4+y_6)(y_5+y_7) + y_3(y_4+y_7)(y_5+y_6) \notag\\
& \le y_2y_3(y_4+y_5+y_6+y_7) + y_2 \left(\frac{y_4+y_6+y_5+y_7}{2}\right)^2+ y_3 \left(\frac{y_4+y_7+y_5+y_6}{2}\right)^2 \notag\\
& = \sigma_3(y_2, y_3, (y_4+y_5+y_6+y_7)/2, (y_4+y_5+y_6+y_7)/2)
\le \binom{4}{3}\left(\frac{y_2+ \cdots + y_7}{4}\right)^3
= \frac{1}{16}(1-y_1)^3. \notag
\end{align}
Therefore, $\lambda_{B_{5,3}}(\B y) \le 6 \left(\frac{3}{8} y_1\left(1-y_1\right)^2 + \frac{1}{16}(1-y_1)^3\right)$.
Similar to $g_1(\B x)$ we obtain
\begin{align}
g_2(\B y)
 = \frac{\lambda_{B_{5,3}}(\B y)}{1-y_1^3}
 \le \frac{6 \left(\frac{3}{8} y_1\left(1-y_1\right)^2 + \frac{1}{16}(1-y_1)^3\right)}{1-y_1^3}
 \le \frac{3(\sqrt{7}-2)}{4}, \notag
\end{align}
and, as it is easy to check, equality holds if and only if (\ref{equ:Lagrange-maxmizer-2}) holds.
This gives all claims of the lemma.
\epf

For the next lemma we need the following definitions.
Given a collection $\{P_i = (m_i, E_i, R_i) \colon i\in I\}$ of $r$-graph patterns,
we define a family $\Sigma^k P_I$ recursively for every integer $k\ge 0$ in the following way. (Note that it is different from the family $\Sigma_k P_I$ that appeared in the proof of Lemma~\ref{LEMMA:lambda-sigma}.)

\begin{definition}[The $k$-th $P_I$-mixing construction]
Let
\begin{align*}
    \Sigma^0 P_I \coloneqq  \bigcup_{i\in I} \left\{ H \colon \text{$H$ is a $P_i$-construction} \right\}.
\end{align*}
For every integer $k\ge 1$ an $r$-graph $H$ is a \textbf{$k$-th $P_I$-mixing construction} if there exist $i\in I$ and a partition $V(H) = V_1 \cup \cdots \cup V_{m_i}$ such that $H$ can be obtained from the blowup $\blow{E_i}{V_1, \ldots, V_{m_i}}$ by adding, for each $j\in R_i$, an arbitrary $(k-1)$-th $P_I$-mixing construction on $V_j$.
Let $\Sigma^k P_I$ denote the collection of all $k$-th $P_I$-mixing constructions.
\end{definition}

Informally speaking, in a $k$-th $P_I$-mixing construction we have to fix $i\in I$ and are allowed to use only pattern $P_i$ on all levels larger than $k$. 
It easy to see that $\Sigma^k P_I \subseteq \Sigma^{k'} P_I$ for all $k'\ge k$ and,  if $\lambda_{P_i}=\lambda$ for each $i\in I$, then the maximum asymptotic density attainable by $\Sigma^k P_I$ for each $k\ge 0$ is $\lambda$. Moreover,
\begin{align*}
    \bigcup_{k\ge 0}\Sigma^k P_I = \Sigma P_I.
\end{align*}

For every integer $k \ge 0$ let $M_{\Sigma^k P_I}$ (resp. $M_{\Sigma P_I}$) be the collection of points $x\in [0,1]$ such that there exists a sequence $\left(H_n\right)_{n=1}^{\infty}$ of $r$-graphs in $\Sigma^k P_I$ (resp. $\Sigma P_I$) with
\begin{align*}
\lim_{n\to \infty}v(H_n) = \infty, \quad
\lim_{n\to \infty}\rho(H_n) = \lambda_{P_I}, \quad\text{and}\quad
\lim_{n\to \infty}\rho(\partial H_n) = x.
\end{align*}

It is easy to observe that the set $M_{\Sigma P_I}$ is the closure of $\bigcup_{k\ge 0} M_{\Sigma^k P_I}$, i.e.
\begin{align*}
    M_{\Sigma P_I}
    = \overline{\bigcup_{k\ge 0} M_{\Sigma^k P_I}},
\end{align*}
which proves Observation~\ref{OBS:shadow-density-mixing-construction}.

We will need the following theorem (which is a rather special case of e.g.\ \cite[Theorems (1) and (3)]{Hutchinson81}) for determining the Hausdorff dimension of a self-similar set.

\begin{theorem}
\label{THM:open-set-condition}
Suppose that $m\ge 2$ and,
for each $i\in [m]$, $\psi_i(x)\coloneqq r_i (x-x_i)+ x_i$ a linear map with $r_i,x_i\in\I R$ and $|r_i|<1$. 
Additionally, suppose that this collection of maps $\{\psi_1, \ldots, \psi_{m}\}$ satisfies the \textbf{open set condition}, namely that
there exists a non-empty open set $V$ such that
\begin{enumerate}
    \item $\bigcup_{i\in[m]}\psi_i(V) \subseteq V$, and
    \item the sets in $\{\psi_{i}(V) \colon i\in [m]\}$ are pairwise disjoint.
\end{enumerate}
Then the Hausdorff dimension $d$ of the (unique) bounded closed non-empty set $A$ that satisfies $A = \bigcup_{i\in[m]}\psi_i(A)$ is the (unique) solution of the equation
\begin{align*}
    \sum_{i\in [m]}|r_i|^d = 1.
\end{align*}
\end{theorem}

The following lemma determines the Hausdorff dimension of the set $M_{\Sigma\{P_1, P_2\}}$.

\begin{lemma}\label{LEMMA:Huasdorff}
We have $M_{\Sigma^0\{P_1,P_2\}} = \left\{\frac{6-\sqrt{7}}{4},\frac{22-3 \sqrt{7}}{16}\right\}$,
and for every integer $k\ge 0$ we have
\begin{align}\label{equ:rec-shadow}
M_{\Sigma^{k+1} \{P_1, P_2\}}  =
                        \left\{ \psi_1(x) \colon x\in M_{\Sigma^{k} \{P_1, P_2\}} \right\} \cup
                        \left\{ \psi_2(x) \colon x\in M_{\Sigma^{k} \{P_1, P_2\}} \right\},
\end{align}
where we define, for every $x\in\I R$,
\begin{align*}
\psi_1(x) 
& \coloneqq  \frac{11-4\sqrt{7}}{9} \left(x - \frac{6-\sqrt{7}}{4} \right) + \frac{6-\sqrt{7}}{4} \quad\text{and}\quad \\
\psi_2(x) 
& \coloneqq  \frac{11-4\sqrt{7}}{9} \left(x - \frac{22-3\sqrt{7}}{16} \right) + \frac{22-3\sqrt{7}}{16}. \notag
\end{align*}
Moreover, the Hausdorff dimension of $M_{\Sigma\{P_1, P_2\}}$ is
$\delta\coloneqq \frac{\log 2}{\log \left(4 \sqrt{7}+11\right)} \approx 0.225641$.
\end{lemma}
\bpf
Lemma~\ref{LEMMA:Lagrangian-P1-P2} implies that $M_{\Sigma^0{\{P_1,P_2\}}} = \{a,b\}$, where $a\coloneqq \frac{6-\sqrt{7}}{4}$ and $b\coloneqq \frac{22-3 \sqrt{7}}{16}$.
Here, $a$ is obtained by solving the equation 
\begin{align*}
    1-  \left( \frac{\sqrt{7}-2}{3} \right)^2 - 4\left(\frac{5-\sqrt{7}}{12}\right)^2 
    + \left(\frac{\sqrt{7}-2}{3}\right)^2 x =x, 
\end{align*}
and $b$ is obtained by solving the equation 
\begin{align*}
    1-  \left( \frac{\sqrt{7}-2}{3} \right)^2  - 2\left(\frac{5-\sqrt{7}}{12}\right)^2 - 4\left(\frac{5-\sqrt{7}}{24}\right)^2 
    + \left(\frac{\sqrt{7}-2}{3}\right)^2 x = x.  
\end{align*}
Let $k \ge 1$. 
It follows from the definition of $M_{\Sigma^{k+1} \{P_1, P_2\}}$ that $\alpha \in M_{\Sigma^{k+1} \{P_1, P_2\}}$ if and only if there exists $\beta \in M_{\Sigma^{k} \{P_1, P_2\}}$ such that
\begin{align*}
    \alpha 
     = 1-  \left( \frac{\sqrt{7}-2}{3} \right)^2 - 4\left(\frac{5-\sqrt{7}}{12}\right)^2  
    + \left(\frac{\sqrt{7}-2}{3}\right)^2 \beta 
     = \frac{13\sqrt{7}-20}{18} + \frac{11-4\sqrt{7}}{9} \beta 
    = \psi_{1}(\beta),
\end{align*}
or 
\begin{align*}
    \alpha 
    & = 1-  \left( \frac{\sqrt{7}-2}{3} \right)^2 - 2\left(\frac{5-\sqrt{7}}{12}\right)^2 - 4\left(\frac{5-\sqrt{7}}{24}\right)^2 
    + \left(\frac{\sqrt{7}-2}{3}\right)^2 \beta \\
    & = \frac{47\sqrt{7} - 64}{72} + \frac{11-4\sqrt{7}}{9} \beta  
    = \psi_{2}(\beta).
\end{align*}
This proves~(\ref{equ:rec-shadow}).

Next, we prove that the Hausdorff dimension of $M_{\Sigma \{P_1, P_2\}}$ is~$\delta$.
Let $A \coloneqq  M_{\Sigma \{P_1, P_2\}}$ and $V \coloneqq  (a,b)$.
Note that $A$ is a bounded closed set that, by~\eqref{equ:rec-shadow}, satisfies $A=\psi_1(A)\cup \psi_2(A)$. Also,
routine calculations show that $\psi_1(V)=(a,c)$ and $\psi_2(V)=(d,b)$, where $c \coloneqq  \frac{166-17\sqrt{7}}{144}$ is less than $d \coloneqq  \frac{124-23\sqrt7}{72}$. Thus
the open set $V$ and the maps $\psi_1$ and $\psi_2$ satisfy the Open Set Condition, i.e.\
\begin{align*}
    \psi_1(V) \cup \psi_2(V) \subseteq V \quad\text{and}\quad
    \psi_1(V) \cap \psi_2(V) = \emptyset.
\end{align*}

On the other hand, recall that $A\subseteq [0,1]$ is a closed set satisfying $A = \psi_1(A) \cup \psi_2(A)$.
So, by Theorem~\ref{THM:open-set-condition}, the Hausdorff dimension of $A$ is the unique solution $x$ to
\begin{align*}
    \left(\frac{11-4\sqrt{7}}{9}\right)^{x} + \left(\frac{11-4\sqrt{7}}{9}\right)^{x} =1,
\end{align*}
which is $\frac{\log 2}{\log \left(4 \sqrt{7}+11\right)}$.
\epf

Now Corollary~\ref{CORO:feasible-function-Cantor} is an easy consequence of Theorem~\ref{THM:feasibe-region} and Lemma~\ref{LEMMA:Huasdorff}.

\section{Concluding remarks}\label{se:conclusion}

Since our paper is quite long, we restricted applications of Theorem~\ref{THM:mixed-pattern} to $3$-graphs only. Some of these results extends to general $r$ while such an extension for others seems quite challenging.
For example, our proof of Theorem~\ref{THM:feasibe-region} also works for $r$-graph patterns and the $(r-2)$-th shadow (when we consider pairs of vertices covered by edges). However, we do not know if the analogue of Theorem~\ref{THM:feasibe-region} holds for the $(r-i)$-th shadow when $i\ge 3$.

On the other hand, we tried to present a fairly general version of Theorem~\ref{THM:mixed-pattern} (for example, allowing $P_I$ to contain patterns with different Lagrangians) in case it may be useful for some other applications.

In some rather special cases of Theorem~\ref{THM:mixed-pattern}, it may be possible to drop the constraint that each $P_i\in P$ is minimal.
For example, it is shown in~\cite{HLLMZ22} that for every (not necessarily minimal) pattern $(m,H,\emptyset$) where $H$ consists of simple $r$-sets, there exists a finite forbidden family
whose extremal Tur\'an constructions are exactly maximum blowups of $H$.
However, we do not know if this is true in general, since the minimality condition is crucially used in the proof of Theorem~\ref{THM:mixed-pattern}.

\section*{Acknowledgements}
We would like to thank  J{\'o}zsef Balogh,  Felix Christian Clemen, Victor Falgas-Ravry, and Dhruv Mubayi for helpful discussions and/or comments.
We are very grateful to a referee for carefully reading the manuscript and for providing numerous suggestions that significantly improved the presentation. 

\renewcommand{\baselinestretch}{0.98}\small
\bibliography{mixpattern}

\begin{thebibliography}{10}

\bibitem{BT11}
R.~Baber and J.~Talbot.
\newblock Hypergraphs do jump.
\newblock {\em Combin. Probab. Comput.}, 20(2):161--171, 2011.

\bibitem{BaloghClemenLuo}
J.~Balogh, F.~C. Clemen, and H.~Luo.
\newblock Non-degenerate hypergraphs with exponentially many extremal
  constructions.
\newblock E-print arxiv:2208.00652, 2022.

\bibitem{BR83}
W.~G. Brown.
\newblock On an open problem of {P}aul {T}ur{\'a}n concerning 3-graphs.
\newblock In {\em Studies in pure mathematics}, pages 91--93. Springer, 1983.

\bibitem{Cvetkovski12i}
Z.~Cvetkovski.
\newblock {\em Inequalities}.
\newblock Springer, Heidelberg, 2012.
\newblock Theorems, techniques and selected problems.

\bibitem{E64}
P.~Erd\H{o}s.
\newblock On extremal problems of graphs and generalized graphs.
\newblock {\em Israel J. Math.}, 2:183--190, 1964.

\bibitem{Erdos67a}
P.~Erd{\H{o}}s.
\newblock Some recent results on extremal problems in graph theory. {R}esults.
\newblock In {\em {Theory of Graphs (Internat. Sympos., Rome}, 1966)}, pages
  117--123 (English); pp. 124--130 (French). Gordon and Breach, New York, 1967.

\bibitem{Erdos81c}
P.~Erd{\H o}s.
\newblock On the combinatorial problems {I} would most like to see solved.
\newblock {\em Combinatorica}, 1:25--42, 1981.

\bibitem{ErdosSimonovits66}
P.~Erd{\H{o}}s and M.~Simonovits.
\newblock A limit theorem in graph theory.
\newblock {\em Stud.\ Sci.\ Math.\ Hungar.}, 1:51--57, 1966.

\bibitem{ErdosStone46}
P.~Erd{\H o}s and A.~H. Stone.
\newblock On the structure of linear graphs.
\newblock {\em Bull.\ Amer.\ Math.\ Soc.}, 52:1087--1091, 1946.

\bibitem{Fonderflaass88}
D.~G. Fon-der Flaass.
\newblock A method for constructing {$(3,4)$}-graphs.
\newblock {\em Mat.\ Zametki}, 44:546--550, 1988.

\bibitem{FF84}
P.~{Frankl} and Z.~{F\"uredi}.
\newblock {An exact result for 3-graphs}.
\newblock {\em {Discrete Math.}}, 50:323--328, 1984.

\bibitem{FranklPengRodlTalbot07}
P.~Frankl, Y.~Peng, V.~R{\"o}dl, and J.~Talbot.
\newblock A note on the jumping constant conjecture of {E}rd{\H o}s.
\newblock {\em J.\ Combin.\ Theory\ {(B)}}, 97:204--216, 2007.

\bibitem{FR84}
P.~Frankl and V.~R\"{o}dl.
\newblock Hypergraphs do not jump.
\newblock {\em Combinatorica}, 4(2-3):149--159, 1984.

\bibitem{Froh08}
A.~Frohmader.
\newblock More constructions for {T}ur\'{a}n's {$(3,4)$}-conjecture.
\newblock {\em Electron. J. Combin.}, 15(1):Research Paper 137, 23, 2008.

\bibitem{HLLMZ22}
J.~Hou, H.~Li, X.~Liu, D.~Mubayi, and Y.~Zhang.
\newblock Hypergraphs with infinitely many extremal constructions.
\newblock {\em Discrete Anal.}, pages Paper No. 18, 34, 2023.

\bibitem{Hutchinson81}
J.~E. Hutchinson.
\newblock Fractals and self-similarity.
\newblock {\em Indiana Univ. Math. J.}, 30(5):713--747, 1981.

\bibitem{KA68}
G.~Katona.
\newblock A theorem of finite sets.
\newblock In {\em Theory of graphs ({P}roc. {C}olloq., {T}ihany, 1966)}, pages
  187--207, 1968.

\bibitem{KatonaNemetzSimonovits64}
G.~O.~H. Katona, T.~Nemetz, and M.~Simonovits.
\newblock On a graph problem of {Tur\'an} {(In Hungarian)}.
\newblock {\em Mat.\ Fiz.\ Lapok}, 15:228--238, 1964.

\bibitem{KE11}
P.~Keevash.
\newblock Hypergraph {T}ur\'{a}n problems.
\newblock In {\em Surveys in combinatorics 2011}, volume 392 of {\em London
  Math. Soc. Lecture Note Ser.}, pages 83--139. Cambridge Univ. Press,
  Cambridge, 2011.

\bibitem{Keevash18}
P.~Keevash.
\newblock Counting designs.
\newblock {\em J.\ Europ.\ Math.\ Soc}, 20:903--927, 2018.

\bibitem{Kirkman47}
T.~P. Kirkman.
\newblock On a problem in combinations.
\newblock {\em Cambridge and Dublin Math.\ J.}, 2:191--204, 1847.

\bibitem{KO82}
A.~V. {Kostochka}.
\newblock {A class of constructions for {T}ur{\'a}n's (3,4)-problem.}
\newblock {\em {Combinatorica}}, 2:187--192, 1982.

\bibitem{KR63}
J.~B. Kruskal.
\newblock The number of simplices in a complex.
\newblock In {\em Mathematical optimization techniques}, pages 251--278. Univ.
  of California Press, Berkeley, Calif., 1963.

\bibitem{Liu20a}
X.~Liu.
\newblock Cancellative hypergraphs and {S}teiner triple systems.
\newblock {\em J. Combin. Theory Ser. B}, 167:303--337, 2024.

\bibitem{LM1}
X.~Liu and D.~Mubayi.
\newblock The feasible region of hypergraphs.
\newblock {\em J. Comb. Theory, Ser. B}, 148:23--59, 2021.

\bibitem{LM22}
X.~Liu and D.~Mubayi.
\newblock A hypergraph {T}ur{\'a}n problem with no stability.
\newblock {\em Combinatorica}, pages 1--30, 2022.

\bibitem{LMR1}
X.~Liu, D.~Mubayi, and C.~Reiher.
\newblock Hypergraphs with many extremal configurations.
\newblock {\em Israel. J. Math.}
\newblock to appear.

\bibitem{LMR3}
X.~Liu, D.~Mubayi, and C.~Reiher.
\newblock Hypergraphs with many extremal configurations {II}.
\newblock In preparation, 2022.

\bibitem{MS65}
T.~S. Motzkin and E.~G. Straus.
\newblock Maxima for graphs and a new proof of a theorem of {T}ur\'{a}n.
\newblock {\em Canadian J. Math.}, 17:533--540, 1965.

\bibitem{NorinYepremyan17}
S.~Norin and L.~Yepremyan.
\newblock Tur{\'a}n number of generalized triangles.
\newblock {\em J.\ Combin.\ Theory\ {(A)}}, 146:312--343, 2017.

\bibitem{PI14}
O.~Pikhurko.
\newblock On possible {T}ur\'{a}n densities.
\newblock {\em Israel J. Math.}, 201(1):415--454, 2014.

\bibitem{Pikhurko15}
O.~Pikhurko.
\newblock The maximal length of a gap between {$r$}-graph {Tur\'an} densities.
\newblock {\em Electronic J.\ Combin.}, 22:7pp., 2015.

\bibitem{PikhurkoSliacanTyros19}
O.~Pikhurko, J.~{Sliac\v an}, and K.~Tyros.
\newblock Strong forms of stability from flag algebra calculations.
\newblock {\em J.\ Combin.\ Theory\ {(B)}}, 135:129--178, 2019.

\bibitem{RS09}
V.~R\"{o}dl and M.~Schacht.
\newblock Generalizations of the removal lemma.
\newblock {\em Combinatorica}, 29(4):467--501, 2009.

\bibitem{Sido95}
A.~Sidorenko.
\newblock What we know and what we do not know about {T}ur\'{a}n numbers.
\newblock {\em Graphs Combin.}, 11(2):179--199, 1995.

\bibitem{SI68}
M.~Simonovits.
\newblock A method for solving extremal problems in graph theory, stability
  problems.
\newblock In {\em Theory of {G}raphs ({P}roc. {C}olloq., {T}ihany, 1966)},
  pages 279--319. Academic Press, New York, 1968.

\bibitem{Turan41}
P.~{Tur\'an}.
\newblock On an extremal problem in graph theory (in {Hungarian}).
\newblock {\em Mat.\ Fiz.\ Lapok}, 48:436--452, 1941.

\bibitem{YanPeng21}
Z.~Yan and Y.~Peng.
\newblock Non-jumping {T}ur\'an densities of hypergraphs.
\newblock {\em Discrete Math.}, 346(1):Paper No. 113195, 11, 2023.

\end{thebibliography}

\end{document}